		\documentclass[a4paper]{book}

\usepackage{latexsym}
\usepackage{amsfonts}
\usepackage[reqno]{amsmath}
\usepackage{amsthm}
\usepackage{amssymb}

\usepackage{bbm}
\usepackage{mathrsfs}
\usepackage{cases}

\usepackage{pstricks}
\usepackage{pst-eucl}
\usepackage{pstricks-add}
\usepackage{graphicx}			%
\usepackage{calc}			%

\usepackage{ifthen}
\usepackage{multicol}
\usepackage[nottoc]{tocbibind}

\usepackage[utf8]{inputenc}
\usepackage[T1]{fontenc}

\usepackage{textcomp}
\usepackage{lmodern}
\usepackage[a4paper]{geometry}

\usepackage[ps2pdf]{hyperref}
\hypersetup{colorlinks=true,linkcolor=blue}
  \usepackage[numbers,sort&compress]{natbib}

\newcounter{numtho}

\newtheoremstyle{mes_tho}%
		{9pt}{9pt}%
		{\itshape}%
		{}%
		{\bfseries}{.}%
		{\newline}%
		{}%

\theoremstyle{remark}

\theoremstyle{mes_tho}	\newtheorem{lemma}[numtho]{Lemma}
			\newtheorem{theorem}[numtho]{Theorem}
			\newtheorem{remark}[numtho]{Remark}
			\newtheorem{corollary}[numtho]{Corollary}
			\newtheorem{proposition}[numtho]{Proposition}

\newcommand{\arxiv}[2][]{%
\ifthenelse{\equal{#1}{}}{%
	\href{http://www.arxiv.org/abs/#2}{{\tt arXiv:#2}}%
			}
			{%
	\href{http://www.arxiv.org/abs/#2}{{\tt arXiv:#2}[#1]}%
}%
}				%

\newcommand{\defe}[2]{\textbf{#1}\index{#2}}
\newcommand{\tq}{\text{ st }}

\DeclareMathOperator{\SO}{SO}
\DeclareMathOperator{\so}{\mathfrak{so}}
\DeclareMathOperator{\ad}{ad}
\DeclareMathOperator{\Ad}{Ad}
\DeclareMathOperator{\AD}{\textbf{Ad}}
\DeclareMathOperator{\Adh}{Adh}
\DeclareMathOperator{\Int}{Int}
\DeclareMathOperator{\id}{id}
\DeclareMathOperator{\pr}{\texttt{pr}}
\DeclareMathOperator{\tr}{Tr}

\newcommand{\hH}{\mathscr{H}}
\newcommand{\hS}{\mathscr{S}}			%
\newcommand{\hF}{\mathscr{F}}			%

\newcommand{\sA}{\mathcal{A}}
\newcommand{\sG}{\mathcal{G}}
\newcommand{\sH}{\mathcal{H}}			%
\newcommand{\sK}{\mathcal{K}}			%
\newcommand{\sN}{\mathcal{N}}
\newcommand{\sP}{\mathcal{P}}
\newcommand{\sQ}{\mathcal{Q}}
\newcommand{\sR}{\mathcal{R}}
\newcommand{\sS}{\mathcal{S}}
\newcommand{\sZ}{\mathcal{Z}}

\newcommand{\mO}{\mathcal{O}}
\newcommand{\mZ}{\mathcal{Z}}

\newcommand{\lG}{\mathfrak{g}}

\newcommand{\eB}{\mathbbm{B}}			%

\newcommand{\eR}{\mathbbm{R}}

\newcommand{\mtu}{\mathbbm{1}}  			%

\newcommand{\Dsddb}[4]{\frac{d}{d#2}\Big[#1\Big]_{#3=#4}}
\newcommand{\Dsdd}[3]{ \Dsddb{#1}{#2}{#2}{#3}   }

\newcommand\abstractname{Abstract}
\makeatletter
  \newenvironment{abstract}{%
      \if@twocolumn
        \section*{\abstractname}%
      \else
        \small
        \begin{center}%
          {\bfseries \abstractname\vspace{-.5em}\vspace{\z@}}%
        \end{center}%
        \quotation
      \fi}
      {\if@twocolumn\else\endquotation\fi}
\makeatother

\makeindex

\begin{document}

  \begin{center}%
    {\LARGE BTZ black hole from the structure of $\so(2,n)$ \par}%
    \vskip 3em%
    {\large
    Laurent Claessens}
      \vskip 1.5em%
    {\large \today}%
  \end{center}

		\begin{abstract}
	In this paper, we study the relevant structure of the algebra $\so(2,n)$ which makes the BTZ black hole possible in the anti de Sitter space $AdS=SO(2,n)/SO(1,n)$. We pay a particular attention on the reductive Lie algebra structures of $\so(2,n)$ and we study how this structure evolves when one increases the dimension. 

	As in \cite{lcTNAdS} and \cite{BTZ_horizon}, we define the singularity as the closed orbits of the Iwasawa subgroup of the isometry group of anti de Sitter, but here, we insist on an alternative (closely related to the original conception of the BTZ black hole) way to describe the singularity as the loci where the norm of fundamental vector vanishes. We provide a manageable Lie algebra oriented formula which describes the singularity and we use it in order to derive the existence of a black hole and to give a geometric description of the horizon. We also define a coherent structure of black hole on $AdS_2$.

	This paper contains a ``short'' and a ``long'' version of the text. In the short version, only the main results are exposed and the proofs are reduced to the most important steps in order to be easier to follow the developments. The long version contains all the intermediate steps and computations for the sake of completeness.

\end{abstract}
\tableofcontents

\chapter{Short version}
\section{Structure of the algebra}
\label{SHORTSecProgressRidMatrices}

\subsection{The Iwasawa component}

Our study of $AdS_l=\SO(2,l-1)/\SO(1,l-1)$ will be based on the properties of the algebra $\sG=\so(2,l-1)$ endowed with a Cartan involution $\theta$ and an Iwasawa decomposition $\sG=\sA\oplus\sN\oplus\sK$. In this section we want to underline the most relevant facts for our purpose. The part we are mainly interested in is the Iwasawa component $\sA\oplus\sN$ where
\begin{subequations}		\label{SHORTEqLeANEnDimAlg}
\begin{align}
	\sN&=\{ X^{k}_{+0},X^{k}_{0+},X_{++},X_{+-} \}\\
	\sA&=\{ J_1, J_2\},
\end{align}
\end{subequations}
where $k$ runs %
from $3$ to $l-1$. The commutator table is
\begin{subequations}  \label{SHORTEqTableSOIwa}
	\begin{align}
		[X_{0+}^{k},X_{+0}^{k'}]&=\delta_{kk'}X_{++}		&[X_{0+}^{k},X_{+-}]&=2X_{+0}^{k}\\
		[ J_1,X_{+0}^{k}]&=X_{+0}^{k}				&[ J_2,X_{0+}^{k}]&=X_{0+}^{k}\\
		[ J_1,X_{+-}]&=X_{+-}					&[ J_2,X_{+-}]&=-X_{+-}\\
		[ J_1,X_{++}]&=X_{++}					&[ J_2,X_{++}]&=X_{++}.
	\end{align}
\end{subequations}
We see that the Iwasawa algebra belongs to the class of $j$-algebras whose Pyatetskii-Shapiro decomposition is 
\begin{equation}
	\sA\oplus\sN=(\sA_1\oplus_{\ad}\sZ_1)\oplus_{\ad}\big( \sA_2\oplus_{\ad}(V\oplus \sZ_2) \big),
\end{equation}
with
\begin{subequations}
	\begin{align}
		\sA_1&=\langle H_1\rangle		&\sA_2&=\langle H_2\rangle\\
		\sZ_1&=\langle X_{+-}\rangle		&\sZ_2&=\langle X_{++}\rangle\\
		&					&V&=\langle X_{0+}^{k},X_{+0}^{k}\rangle_{k\geq 3}
	\end{align}
\end{subequations}
where 
\begin{equation}
	\begin{aligned}[]
		H_1&=J_1-J_2\\
		H_2&=J_1+J_2.
	\end{aligned}
\end{equation}
\subsection{Dimensional slices}
\label{SHORTSubSecDimensionalSlices}

Since the elements of $\mZ_{\sK}(\sA)$ is a part of $\sH$ (see later), they will have almost no importance in the remaining\footnote{We will however need them in the computation of the coefficients \eqref{SHORTEqCoefsabcBE}.}. The most important part is
\begin{equation}		\label{SHORTEqDecomDimQlipm}
	\sA\oplus\sN\oplus\bar\sN=\underbrace{\langle J_1,J_2,X_{\pm,\pm}\rangle}_{\text{for every dimension}}\oplus\underbrace{\langle X_{0\pm}^{4},X_{\pm 0}^{4}\rangle}_{\text{for $\so(2,\geq 3)$}}\oplus\ldots\oplus\underbrace{\langle X_{0\pm}^l,X_{\pm 0}^{l}\rangle}_{\text{for $\so(2, l-1)$}}.
\end{equation}
We use the following notations in order to make more clear how does the algebra evolve when one increases the dimension:
\begin{equation}
	\begin{aligned}[]
		\sN_2&=\langle X_{+-},X_{++}\rangle,		&\sN_k&=\langle X^k_{0+},X_{+0}^k\rangle	\\
		\bar\sN_2&=\langle X_{-+},X_{--} \rangle, 	&\bar\sN_k&=\langle X_{0-}^k,X_{-0}^k\rangle\\
		\tilde\sN_2&=\langle \sN_2,\bar\sN_2\rangle,	&\tilde\sN_k&=\langle \sN_k,\bar\sN_k\rangle
	\end{aligned}
\end{equation}
for $k\geq 3$. The relations are
\begin{equation}		\label{SHORTEqsCommWithtsNDeuxkA}
	\begin{aligned}[]
		[\tilde\sN_2,\tilde\sN_2]&\subseteq\sA
		&[\tilde\sN_2,\tilde\sN_k]&\subseteq\tilde\sN_k\\
		[\tilde\sN_k,\tilde\sN_{k}]&\subseteq\sA\oplus\tilde\sN_2
		&[\tilde\sN_k,\tilde\sN_{k'}]&\subset \mZ_{\sK}(\sA) \\
	\end{aligned}
\end{equation}
As a consequence of the splitting and these commutation relations,
\begin{equation}		\label{SHORTEqDecompGPourKillOrtho}
	\sG=\mZ_{\sK}(\sA)\oplus\sA\oplus\tilde\sN_2\bigoplus_{k\geq 3}\tilde\sN_k
\end{equation}
is a Killing-orthogonal decomposition of $\so(2,l-1)$.

\subsection{Reductive decomposition}
\label{SHORTSubSecReductiveDecompQ}

Let $\sQ$ be the following vector subspace of $\sG$:
\begin{equation}		\label{SHORTEqDecQEspacesCools}	%
	\sQ=\big\langle \mZ(\sK),J_2,[\mZ(\sK),J_1],(X_{0+}^k)_\sP\big\rangle_{k\geq 3}.
\end{equation}
Then we choice a subalgebra $\sH$ of $\sG$ which, as vector space, is a complementary of $\sQ$. In that choice, we require that there exists an involutive automorphism $\sigma\colon \sG\to \sG$ such that
\begin{equation}
	\sigma=(\id)_{\sH}\oplus(-\id)_{\sQ}.
\end{equation}
In that case the decomposition $\sG=\sH\oplus\sQ$ is reductive, i.e. $[\sQ,\sQ]\subset\sH$ and $[\sQ,\sH]\subset\sQ$.

\begin{remark}

The space $\mZ(\sK)$ is given by the structure of the compact part of $\so(2,n)$, the elements $(X_{0+}^k)_{\sP}$ are defined from the root space structure of $\so(2,n)$ and the Cartan involution. The elements $J_1$ and $J_2$ are a basis of $\sA$. However, we need to know $\sH$ in order to distinguish $J_1$ from $J_2$ that are respectively generators of $\sA_{\sH}$ and $\sA_{\sQ}$.

Thus the basis 
\eqref{SHORTEqDecQEspacesCools} %
is given in a way almost independent of the choice of $\sH$.
\end{remark}
\section{Some properties}

In  this section we list some important properties which can be derived from the given commutators relation. We refer to \cite{btz_so} for detailed proofs.

The space $\sQ$ can be set in a nice way if we consider the following basis:
\begin{subequations}		\label{SHORTEqBasQQziPlusMieux}%
	\begin{align}
		q_0&=(X_{++})_{\mZ(\sK)}&\in \sK\cap\sQ\cap\tilde\sN_2					\\
		q_1&=J_2&\in\sQ\cap\sA\\							
		q_2&=-[J_1,q_0]	&\in\sP\cap\sQ\cap\tilde\sN_2	\label{SHORTsubEqJUnQZQDeux}		\\
		q_k&=(X^k_{0+})_{\sP}&\in\sP\cap\sQ\cap\tilde\sN_k.
	\end{align}
\end{subequations}
These elements correspond to the expression \eqref{SHORTEqDecQEspacesCools}.
\begin{corollary}		\label{SHORTCorQdansPetK}
	We have $q_0\in\sK$ and $q_i\in\sP$ if $i\neq 0$ and the set $\{ q_0,q_1,\ldots,q_l \}$ is a basis of $\sQ$. Moreover, we have $\sQ\cap\tilde\sN_k=\langle q_k\rangle$.
\end{corollary}

The basis $\{ q_i \}$ of $\sQ$ allows to write down a nice basis for $\sH$:
\begin{equation}			\label{SHORTAlignPremDefHH}
	\begin{aligned}[]
		J_1&		&r_k&=[J_2,q_k]			\\
		p_1&=[q_0,q_1]	&p_k&=[q_0,q_k]			\\
		s_1&=[J_1,p_1]	&s_k&=[J_1,p_k].
	\end{aligned}
\end{equation}
The $\so(2,l-1)$ part of the algebra is generated by the elements
\begin{equation}
	[q_i,q_j]=-\frac{1}{ 4 }\big( [X_{0+}^i,X_{0-}^j]+[X_{0-}^i,X_{0+}^j] \big)=r_{ij}.
\end{equation}%

The vectors $q_i$ have the property to be intertwined by some elements of $\sH$. Namely, the adjoint action of the elements 
\begin{subequations}\label{SHORTAlignDefMagicIntert}
	\begin{align}
		X_1=p_1&=-[J_2,q_0]	&\in\sP\cap\sH\cap\tilde\sN_2\\
		X_2=s_1&=[J_1,X_1]	&\in\sK\cap\sH\cap\tilde\sN_2\\
		X_k=-r_k&=-[J_2,q_k]	&\in\sK\cap\sH\cap\tilde\sN_k				\label{SHORTEqDefXkCommeComm}
	\end{align}
\end{subequations}%
intertwines the $q_i$'s in the sense of the following proposition.
\begin{proposition}[Intertwining properties]			\label{SHORTXUnALaTwistingSuperCool}
	The elements defined by equation \eqref{SHORTAlignDefMagicIntert} satisfy
	\begin{multicols}{2}
		\begin{subequations}				\label{SHORTEqCalculBBBJUnUnNirme}
			\begin{align}
				\ad(J_1)q_0&=-q_2		\label{SHORTEqCalculBBBJUnUnNirmeA}\\
				\ad(J_1)q_2&=-q_0.		\label{SHORTEqCalculBBBJUnUnNirmeB}
			\end{align}
		\end{subequations}
		\begin{subequations}				\label{SHORTEqSubEqbXUnqZero}
			\begin{align}
				\ad(X_1)q_1&=q_0		\label{SHORTSubEqbXZeroqUn}\\
				\ad(X_1)q_0&=q_1,		\label{SHORTSubEqbXUnqZero}
			\end{align}
		\end{subequations}
		\begin{subequations}
			\begin{align}
				\ad(X_2)q_2&=q_1		\label{SHORTSubEqXdeuxQdeuxa}\\
				\ad(X_2)q_1&=-q_2		\label{SHORTSubEqXdeuxQun}
			\end{align}
		\end{subequations}
		\begin{subequations}				\label{SHORTEqSubEqbXkqZero}
			\begin{align}
				\ad(X_k)q_k&=-q_1.		\label{SHORTSubEqbXIMoinsqZero}\\
				\ad(X_k)q_1&=q_k		\label{SHORTSubEqbXkQunQk}
			\end{align}
		\end{subequations}
	\end{multicols}
\end{proposition}

This circumstance allows us to prove many results  on $q_1=J_2$ and propagate them to the others by adjoint action. The following result is proven using the intertwining elements.
\begin{theorem}			\label{SHORTThoAdESqqq}
	We have
	\begin{equation}
		\ad(q_i)^2q_j=q_j
	\end{equation}
	if $i\neq j$ and $i\neq 0$. If $i=0$, we have
	\begin{equation}
		\ad(q_0)^2q_j=-q_j.
	\end{equation}
\end{theorem}

Being the tangent space of $AdS$, the space $\sQ$ is of a particular importance.
We know that the directions of light like geodesics are given by elements in $\sQ$ which have a vanishing norm\cite{These}. These elements are exactly the ones which are nilpotent. We are thus led to study the norm of the basis vectors $q_i$ as well as the order of the nilpotent elements. 
The following three properties will be central in the deduction of the causal structure:%
\begin{enumerate}
	\item
		the elements $\{ q_i \}_{i=0,\cdots,l-1}$ are Killing-orthogonal,
	\item
		if $E$ is nilpotent in $\sQ$, then $\ad(E)^3=0$,
	\item
		up to renormalization, the elements in $\sQ$ that have vanishing norm are of the form
		\begin{equation}
			E=q_0+\sum_{i=1}^{l-1}w_iq_i
		\end{equation}
		where $w\in S^{l-2}$.
\end{enumerate}

From now, it will be convenient to work with the following Killing-orthogonal basis of $\sG$:
\begin{equation}		\label{SHORTEqSuperBaseeB}
	\eB=\{J_1,J_2, q_0,q_2,p_1,s_1,q_k,p_k,r_k,s_k \}_{k=3,\ldots,l-1}
\end{equation}

From a computational point of view, we will need the following exponentials in order to describe the geometry of the geodesics.
The action of $e^{xq_0}$ on $\sA$ is
\begin{subequations}			\label{SHORTSubEqsAdxqzJJ}
	\begin{align}
		e^{\ad(xq_0)}J_1=\cos(x)J_1+\sin(x)q_2\\
		e^{\ad(xq_0)}J_2=\sin(x)p_1+\cos(x)q_1,
	\end{align}
\end{subequations}
on $\tilde\sN_2$ we have
\begin{subequations}					\label{SHORTSubEqsexpxzqpsdzuu}
	\begin{align}
		e^{\ad(xq_0)}q_2&=\cos(x)q_2-\sin(x)J_1\\
		e^{\ad(xq_0)}q_0&=q_0\\
		e^{\ad(xq_0)}p_1&=\cos(x)p_1-\sin(x)q_1 		\label{SHORTEqAdqzpUn}\\
		e^{\ad(xq_0)}s_1&=s_1,
	\end{align}
\end{subequations}
and the action on $\tilde\sN_k$ is
\begin{subequations}		\label{SHORTEqExpAdqkpk}
	\begin{align}
		e^{\ad(xq_0)}q_k&=\cos(x)q_k+\sin(x)p_k\\
		e^{\ad(xq_0)}p_k&=\cos(x)p_k-\sin(x)q_k\\
		e^{\ad(xq_0)}p_k&=\cos(x)p_k-\sin(x)q_k.
	\end{align}
\end{subequations}
\section{The causally singular structure}
\label{SHORTSecBlacHole}

\subsection{Closed orbits}

The singularity in $AdS_l$ is defined as the closed orbits of $AN$ and $A\bar N$ in $G/H$. This subsection is intended to identify them.

\begin{proposition}		\label{SHORTPropCartanExtExpo}
	The Cartan involution $\theta\colon \sG\to \sG$ is an inner automorphism, namely it is given by $\theta=\Ad(k_{\theta})$ where $k_{\theta}= e^{\pi q_0}$.
\end{proposition}

One check the proposition setting $x=\pi$ in the exponentials given earlier.
In the following results, we use the fact that the group $K$ splits into the commuting product $K=\SO(2)\times\SO(l-1)$.
\begin{lemma}		\label{SHORTLemExistxTqansAbarN}
	For every $an\in AN$, there exists $x\in\mathopen[ 0 , 2\pi [$ such that $[an e^{xq_0}]\in[A\bar N]$.
\end{lemma}

\begin{proof}
	Let $k\in K$ such that $ank\in A\bar N$. The element $k$ decomposes into $k=st$ with $s= e^{xq_0}\in \SO(2)$ and $t\in\SO(n)\subset H$. Thus $[ans]\in[A\bar N]$.
\end{proof}

\begin{lemma}		\label{SHORTLemansse}
	If $[an]=[s]$ with $s\in\SO(2)$, then $s=e$.
\end{lemma}

\begin{proof}
	The assumption implies that there exists a $h\in H$ such that $an=sh$. Using the $KAN$ decomposition of $H$, such a $h$ can be written under the form $h=ta'n'$ with $t\in\SO(n)$. Thus we have $an=sta'n'$. By unicity of the decomposition $kan$, we must have $st=e$, and then $s=e$.
\end{proof}

Using the fact that $k_{\theta}ANk_{\theta}\subset A\bar N$, we also prove the following.
\begin{lemma}		\label{SHORTLemANksk}
	If $an\in AN$ and if $[ank_{\theta}]=[s]$ with $s\in\SO(2)$, then $s=k_{\theta}$.
\end{lemma}

\begin{theorem}		\label{SHORTThoOrbitesOuverttes}
	The closed orbits of $AN$ in $AdS_l$ are $[AN]$ and $[AN k_{\theta}]$ where $k_{\theta}$ is the element of $K$ such that $\theta=\Ad(k_{\theta})$. The closed orbits of $A\bar N$ are $[A\bar N]$ and $[A\bar N k_{\theta}]$. The other orbits are open.
\end{theorem}

\begin{proof}
	Let us deal with the $AN$-orbits in order to fix the ideas. First, remark that each orbit of $AN$ pass trough $[SO(2)]$. Indeed, each $[ank]$ is in the same orbit as $[k]$ with $k\in K=\SO(2)\otimes\SO(n)$. Since $\SO(n)\subset H$, we have $[k]=[s]$ for some $s\in\SO(2)$.

	We are thus going to study openness of the $AN$-orbit of elements of the form $[e^{x q_0}]$ because these elements are ``classifying'' the orbits. Using the isomorphism $  dL_{g^{-1}}\colon T_{[g]}(G/H)\to \sQ$, we know that a set $\{ X_1,\ldots X_l \}$ of vectors in $T_{[ e^{x q_0}]}AdS_l$ is a basis if and only if the set $\{ dL_{ e^{-xq_0}}X_i \}_{i=1,\ldots l}$ is a basis of $\sQ$. We are thus going to study the elements 
	\begin{equation}
		\begin{aligned}[]
			dL_{ e^{-xq_0}}X^*_{[ e^{xq_0}]}	&=dL_{ e^{-xq_0}}\Dsdd{ \pi\big(  e^{-t X} e^{xq_0} \big) }{t}{0}\\
								&=\Dsdd{ \pi\big(  \AD( e^{-xq_0} ) e^{-tX} \big) }{t}{0}\\
								&=-\pr_{\sQ} e^{-\ad(xq_0)}X
		\end{aligned}
	\end{equation}
	when $X$ runs over the elements of $\sA\oplus \sN$.

	Taking the projection over $\sQ$ of the exponentials given by equations \eqref{SHORTSubEqsAdxqzJJ}--\eqref{SHORTEqExpAdqkpk}, we see
	that an orbit trough $[ e^{xq_0}]$ is open if and only if $\sin(x)\neq 0$. It remains to study the orbits of $[ e^{\pi q_0}]$ and $[e]$. Lemma \ref{SHORTLemansse} shows that these two orbits are disjoint.

	Let us now prove that $[AN]$ is closed. A point outside $\pi(AN)$ reads $\pi(ans)$ where $s$ is an elements of $\SO(2)$ which is not the identity. Let $\mO$ be an open neighborhood of $ans$ in $G$ such that every element of $\mO$ read $a'n's't'$ with $s'\neq e$. The set $\pi(\mO)$ is then an open neighborhood of $\pi(ans)$ which does not intersect $[AN]$. This proves that the complementary of $[AN]$ is open. The same holds for the orbit $[A\bar N]$.
	
	The orbit $[ANk_{\theta}]$ and $[A\bar Nk_{\theta}]$ are closed too because $ANk_{\theta}=k_{\theta}A\bar N$.

\end{proof}

\subsection{Vanishing norm criterion}

In the preceding section, we defined the singularity by means of the action of an Iwasawa group. We are now going to give an alternative way of describing the singularity, by means of the norm of a fundamental vector of the action. This ``new'' way of describing the singularity is, in fact, much more similar to the original BTZ black hole where the singularity was created by identifications along the integral curves of a Killing vector field\cite{these_Detournay}. The vector $J_1$ in theorem \ref{SHORTThosSequivJzero} plays here the role of that ``old'' Killing vector field.

Discrete identifications along the integral curves of $J_1$ would produce the causally singular space which is at the basis of our black hole.

What we will prove is the
\begin{theorem}		\label{SHORTThosSequivJzero}
	We have $\sS\equiv \| J_1^* \|_{[g]}=\| \pr_{\sQ}\Ad(g^{-1})J_1 \|=0$.
\end{theorem}

Thanks to this theorem, our strategy will be to compute $\| \pr_{\sQ}\Ad(g^{-1})J_1 \|$ in order to determine if $[g]$ belongs to the singularity or not. The proof will be decomposed in three steps. The first step is to obtain a manageable expression for $\| J_1^* \|$.

\begin{lemma}		\label{SHORTLemExpressionCoolNormJUn}
	We have $\| (J_1^*)_{[g]} \|=\| \pr_{\sQ}\Ad(g^{-1})J_1 \|$ for every $[g]\in AdS_l$.
\end{lemma}

\begin{proof}
	This is a direct application of the formula for the Killing-induced norm on an homogeneous space\cite{Kerin}.
\end{proof}

\begin{proposition}		\label{SHORTPropPtpsSjzero}
	If $p\in\hS$, then $\| J_1^* \|_p=0$.
\end{proposition}

\begin{proof}
	We are going to prove that $\pr_{\sQ}\Ad(g^{-1})J_1$ is a light like vector in $\sQ$ when $g$ belongs to $[AN]$ or $[A\bar N]$. A general element of $AN$ reads $g=a^{-1}n^{-1}$ with $a\in A$ and $n\in N$. Since $\Ad(a)J_1=J_1$, we have $\Ad(g^{-1})J_1=\Ad(n)J_1$. We are going to study the development
	\begin{equation}
		\Ad( e^{Z})J_1= e^{\ad(Z)}J_1=J_1+\ad(Z)J_1+\frac{ 1 }{2}\ad(Z)^2J_1+\ldots
	\end{equation}
	where $Z=\ln(n)\in\sN$. The series is finite because $Z$ is nilpotent  and begins by $J_1$ while all other terms belong to $\sN$. Notice that the same remains true if one replace $\sN$ by $\bar \sN$ everywhere. 
	
	Moreover, $\Ad( e^{Z})J_1$ has no $X_{0+}$-component (no $X_{0-}$-component in the case of $Z\in\bar\sN$) because $[X_{0+},J_1]=0$, so that the term $[Z,J_1]$ is a combination of $X_{+0}$, $X_{++}$ and $X_{+-}$.  Since the action of $\ad(X_{+\pm})$ on such a combination is always zero, the next terms are produced by action of $\ad(X_{0+})$ on a combination of $X_{+0}$, $X_{++}$ and $X_{+-}$. Thus we have
\begin{equation}		\label{SHORTEqAdanJUnabck}
	\Ad( e^{Z})J_1=J_1+aX_{++}+bX_{+-}+c_kX_{+0}^k
\end{equation}
for some\footnote{One can show that every combinations of these elements are possible, but that point is of no importance here.} constants $a$, $b$ and $c_k$.

The projection of $\Ad( e^{Z})J_1$ on $\sQ$ is made of a combination of the projections of $X_{+0}$, $X_{++}$ and $X_{+-}$.  The conclusion is that $\pr_{\sQ}\big( e^{\ad(Z)}J_1\big)$ is a multiple of $q_0+q_2$, which is light like. The conclusion still holds with $\bar\sN$, but we get a multiple of $q_0-q_2$ instead of $q_0+q_2$.

	Now we have $\Ad(k_{\theta})J_1=J_1$ and $\Ad(k_{\theta})(q_0\pm q_2)=-(q_0\pm q_2)$, so that the same proof holds for the closed orbits $[ANk_{\theta}]$ and $[A\bar N k_{\theta}]$.
\end{proof}

\begin{remark}		\label{SHORTRemANANbarYapas}
	The coefficients $a$, $b$ and $c_k$ in equation \eqref{SHORTEqAdanJUnabck} are continuous functions of the starting point $an\in AN$. More precisely, they are polynomials in the coefficients of $X_{++}$, $X_{+-}$, $X_{0+}$ and $X_{0+}$ in $Z$. The vector $\pr_{\sQ}\Ad(g^{-1})J_1=(a+b)(q_0+q_1)$ is thus a continuous function of the point $[g]\in[AN]$.

	If $[g]\in[AN]\cap[A\bar N]$, then $\pr_{\sQ}\Ad(g^{-1})J_1$ has to vanish as it is a multiple of $q_0+q_1$ and of $q_0-q_1$ in the same time. We conclude that in each neighborhood in $[AN]$ of an element of $[AN]$, there is an element which does not belong to $[A\bar N]$.
\end{remark}

\begin{proposition}
	If $\| J_1^* \|_p=0$, then $p\in\sS$.
\end{proposition}

\begin{proof}
	As before we are looking at a point $[g]=[(an)^{-1}s^{-1}]$ with $s= e^{xq_0}$. The norm $\| J_1^* \|$ vanishes if
	\begin{equation}
		\| \pr_{\sQ} \Ad( e^{xq_0})\Ad(an)J_1 \|=0.
	\end{equation}
	We already argued in the proof of proposition \ref{SHORTPropPtpsSjzero} that $\Ad(an)J_1$ is equal to $J_1$ plus a linear combination of $X_{++}$, $X_{+-}$ and $X_{+0}$. A straightforward computation yields
	\begin{equation}		\label{SHORTEqprQexpxqzXanroots}
		\begin{aligned}[]
			\pr_{\sQ} e^{\ad(xq_0)}(J_1&+aX_{++}+bX_{+-}+\sum_k c_kX^k_{+0})\\
								&=(a+b)q_0+(a-b)\sin(x)q_1+\big( \sin(x)-(a+b)\cos(x) \big)q_2+\sum _k c_k\sin(x)q_k.
		\end{aligned}
	\end{equation}
	The norm of this vector, as function of $x$, is given by
	\begin{equation}
		n(x)=(a+b)\sin(2x)+(4ab-C^2-1)\big( 1-\cos(2x) \big),
	\end{equation}
	where $C^2=\sum_kc_k^2$. Using the variables with $u=a+b$ and $v=(1+c^2-4ab)/2$,
	\begin{equation}
		n(x)=u\sin(2x)+v\cos(2x)-v.
	\end{equation}
	Following $u=0$ or $u\neq 0$, the graph of that function vanishes two or four times between $0$ and $2\pi$. Points of $[AN]$ are divided into two parts: the \emph{red points} which correspond to $u\neq 0$, and the \emph{blue points} which correspond to $u=0$. By continuity, the red part is open.

	Let $P=an\in AN$. 
	We consider the unique  $x_0$, $x_1$, $x_2$ and $x_3$ in $\mathopen[ 0 , 2\pi [$ such that 
\begin{subequations}
	\begin{align}
		[P e^{x_0q_0}]&\in[AN]\\
		[P e^{x_1q_0}]&\in[ANk_{\theta}]\\
		[P e^{x_2q_0}]&\in[A\bar N]\\
		[P e^{x_3q_0}]&\in[A\bar Nk_{\theta}].
	\end{align}
\end{subequations}
They satisfy $x_0=0$, $x_1=\pi$ and $x_3=x_2+\pi$ modulo $2\pi$. Now, we divide $[AN]$ into two parts. The elements of $[AN]\cap [A\bar N]$ and $[AN]\cap[A\bar Nk_{\theta}]$ are said to be of \emph{type I}, while the other are said to be of \emph{type II}. We are going to prove that type I points are exactly blue points, while type II points are the red ones.

If $P$ is a point of type II, we know that the $x_i$ are four different numbers so that the norm function $n_P(x)$ vanishes \emph{at least} four times on the interval $\mathopen[ 0 , 2\pi [$, each of them corresponding to a point in the singularity. But our division of $[AN]$ into red and blue points shows that $n_P(x)$ can vanish \emph{at most} four times. We conclude that a point of type II is automatically red, and that the four roots of $n_P(x)$ correspond to the four values $x_i$ for which $P e^{x_iq_0}$ belongs to the singularity. The proposition is thus proved for points of type II.

Let now $[P]$ be of type I (say $P\in [AN]\cap [A\bar N]$) and let us show that $P$ is blue. We consider a sequence of points $P_k$ of type II which converges to $P$ (see remark \ref{SHORTRemANANbarYapas}). We already argued that $P_k$ is red, so that $x_0(P_k)\neq x_2(P_k)$ and $x_1(P_k)\neq x_3(P_k)$, but
\begin{subequations}
	\begin{align}
		x_0(P_k)-x_2(P_k)\to 0\\
		x_1(P_k)-x_3(P_k)\to 0.
	\end{align}
\end{subequations}
The continuity of $n_Q(x)$ with respect to both $x\in\mathopen[ 0 , 2\pi [$ and $Q\in[AN]$ implies that $P$ has to be blue, and then $n_P(x)$ vanishes for exactly two values of $x$ which correspond to $P e^{xq_0}\in\sS$.

Let us now prove that everything is done. We begin by points of type I. Let $P$ be of type $I$ and say $P\in[AN]\cap[A\bar N]$. The curve $n_P(x)$ vanishes exactly two times in $\mathopen[ 0 , 2\pi [$. Now, if $P e^{x_1 q_0}\in[ANk_{\theta}]$, thus $x_1 = \pi$ and we also have $P e^{x_1q_0}\in[A\bar Nk_{\theta}]$, but $P$ does not belong to $[ANk_{\theta}]$, which proves that $n_P(x)$ vanishes \emph{at least} two times which correspond to the points $P e^{xq_0}$ that are in the singularity. Since the curve vanishes in fact exactly two times, we conclude that $n_P(x)$ vanishes if and only if $P e^{xq_0}$ belongs to the singularity.

If we consider a point $P$ of type II, we know that the values of $x_i$ are four different numbers, so that the curve $n_P(x)$ vanishes \emph{at least} four times, corresponding to the points $P e^{xq_0}$ in the singularity. Since the curve is in fact red, it vanishes \emph{exactly} four times in $\mathopen[ 0 , 2\pi [$ and we conclude that the curve $n_P(x)$ vanishes if and only if $P e^{xq_0}$ belongs to the singularity.

The conclusion follows from the fact that 
\begin{equation}
	AdS_l=\Big\{ [P e^{xq_0}] \tq \text{$P$ is of type I or II and $x\in\mathopen[ 0 , 2\pi [$} \Big\}.
\end{equation}

\end{proof}
Proof of theorem \ref{SHORTThosSequivJzero} is now complete.

\subsection{Existence of the black hole}
\label{SHORTSubSecExistenceTrouNoir}

We know that the geodesic trough $[g]$ in the direction $X$ is given by
\begin{equation}
	\pi\big( g e^{sX} \big)
\end{equation}
and that a geodesics is light-like when the \defe{direction}{Direction} $X$ is given by a nilpotent element in $\sQ$\cite{lcTNAdS}.
Let us study the geodesics issued from the point $[ e^{-xq_0}]$. They are given by
\begin{equation}
	l^w_x(s)=\pi\big(    e^{-xq_0} e^{sE(w)} \big)
\end{equation}
where $E(w)=q_0+\sum_iw_iq_i$ with $\| w \|=1$. According to our previous work, the point $l^w_x(s)$ belongs to the singularity if and only if 
\begin{equation}		\label{SHORTEqNormAFaireZeroOuPas}
	n_x^w(s)=\left\|   \pr_{\sQ} e^{-\ad(sE(w))} e^{\ad(xq_0)}J_1  \right\|^2=0.
\end{equation}
Using the exponential previously computed, the fact that $\ad(E)^3=0$, and collecting the terms in $\sQ$ we find%
\begin{equation}		\label{SHORTEqnwxBBB}
	\begin{aligned}[]
		n_x^w(s)	&=s^2\sin^2(x)B\big( \ad(E)^2q_2,q_2 \big)\\
				&\quad+s^2\cos^2(x)B\big( \ad(E)J_1,\ad(E) J_1\big)\\
				&\quad-2s\cos(x)\sin(x)B\big( \ad(E)J_1,q_2 \big)\\
				&\quad+\sin^2(x)B(q_2,q_2).
	\end{aligned}
\end{equation}
The problem now reduces to the evaluation of the three Killing products in this expression. %
Using the relations
\begin{equation}
	\begin{aligned}[]
		q_1&\in\sP\cap\sQ\cap\sA	,	&&	[\tilde\sN_k,\tilde\sN_2]\subset\tilde\sN_k\\
		q_2&\in\sP\cap\sQ\cap\tilde\sN_2,	&&	[\tilde\sN_2,\tilde\sN_2]\subset\sA\\
		q_k&\in\sP\cap\sQ\cap\tilde\sN_k,	&&	[\tilde\sN_k,\tilde\sN_k]\subset\sA\oplus\tilde\sN_2,
	\end{aligned}
\end{equation}%
and the fact that $\ad(E)J_1=w_2q_0+q_2$ we find %
\begin{equation}
	\begin{aligned}[]
		\frac{ n_x^w(s) }{ B(q_2,q_2) }=\big( \cos^2(x)-w^2_2 \big)s^2-2\cos(x)\sin(x)s+\sin^2(x).
	\end{aligned}
\end{equation}
We have $n_x^w(s)=0$ when $s$ equals
\begin{equation}
	s_{\pm}=\frac{ \cos(x)\sin(x)\pm| w_2\sin(x) | }{ \cos^2(x)-w_2^2 }.
\end{equation}
If $w_2\sin(x)\geq 0$, we have\footnote{The solutions \eqref{SHORTEqRacinesellTN} were already deduced in \cite{lcTNAdS} in a quite different way.}
\begin{equation}				\label{SHORTEqRacinesellTN}
	\begin{aligned}[]
		s_+	&=\frac{ \sin(x) }{ \cos(x)-w_2 }&\text{and}&&
		s_-	&=\frac{ \sin(x) }{ \cos(x)+w_2 },
	\end{aligned}
\end{equation}
and if $w_2\sin(x)<0$, we have to exchange $s_+$ with $s_-$.

If we consider a point $ e^{xq_0}$ with $\sin(x)>0$ and $\cos(x)<0$, the directions $w$ with $| w_2 |<| \cos(x) |$ escape the singularity as the two roots \eqref{SHORTEqRacinesellTN} are simultaneously negative. Such a point does not belong to the black hole. That proves that the black hole is not the whole space.

If we consider a point $ e^{xq_0}$ with $\sin(x)>0$ and $\cos(x)>0$, we see that for every $w_2$, we have $s_+>0$ or $s_->0$ (or both). That shows that for such a point, every direction intersect the singularity. Thus the black hole is actually larger than only the singularity itself.

The two points with $\sin(x)=0$ belong to the singularity. At the points $\cos(x)=0$, $\sin(x)=\pm1$, we have $s_+=-1/w_2$ and $s_-=1/w_2$. A direction $w$ escapes the singularity only if $w_2=0$ (which is a closed set in the set of $\| w \|=1$). 

\section{Some more computations}
\label{SHORTSecMoreComputations}

As we saw, the use of theorem \ref{SHORTThosSequivJzero} leads us to study the function
\begin{equation}			\label{SHORTEqprSqbignor}
	n_{[g]}^w(s)=\|   \pr_{\sQ} \Ad\big(  e^{-s E(w)} \big)X  \|^2=0
\end{equation}
where $E(w)=q_0+w_1q_1+\ldots+w_{l-1}q_{l-1}$ ($w\in S^{l-2}$) and (see equation \eqref{SHORTEqAdanJUnabck})
\begin{equation}
	X=e^{\ad(xq_0)}\Ad(na)J_1=e^{\ad(xq_0)}\big(  J_1+aX_{++}+bX_{+-}+\sum_{k=3}^{l-1}c_kX_{+0}^k\big).
\end{equation}
A long but straightforward computation shows that

\begin{equation}
	n_{[g]}^w(s)=\|   \pr_{\sQ} \Ad\big(  e^{-s E(w)} \big)X  \|^2=a(E)s^2+b(E)s+c
\end{equation}
where
\begin{subequations}		\label{SHORTSubEqsabcEBBB}
	\begin{align}
		a(E)&=-B\big( \ad(E)X,\sigma\ad(E)X \big)		\label{SHORTEqCoefaEBX}\\
		b(E)&=-2B\big( X_{\sQ},\ad(E)X_{\sH} \big)\\
		c&=B(X_{\sQ},X_{\sQ}).
	\end{align}
\end{subequations}
More precisely,%

\begin{subequations}\label{SHORTEqCoefsabcBE}
	\begin{align}
		\frac{ a(E) }{ B(q_0,q_0) }&=M\Big( w_2^2+\cos(x)w_2+\cos^2(x)\Big)	   \label{SHORTEqCoeffaE}\\
		\frac{ b(E) }{ B(q_0,q_0) }&=-2M\sin(x)\big( w_2+\cos(x)\big))\\
		\frac{ c }{ B(q_0,q_0) }&=M\sin^2(x)
	\end{align}
\end{subequations}
where $M=\big( a^2-b^2-\sum_kc_k^2 \big)$. The important point to notice is that these expressions only depend on the $w_2$-component of the direction. Notice that $c=0$ if and only if the point $[g]$ belongs to the singularity because $s=0$ is a solution of \eqref{SHORTEqprSqbignor} only in the case $[g]\in\hS$.

\begin{remark}	\label{SHORTRemImapoabcE}
	Importance of the coefficients \eqref{SHORTEqCoefsabcBE}. If $v\in\hF_l$, there is a direction $E_0$ in $AdS_l$ which escapes the singularity from $v$. Thus the polynomial $a(E_0)s^2+b(E_0)s+c$ has only non positive roots. From the expressions \eqref{SHORTEqCoefsabcBE}, we see that the polynomial corresponding to $\iota(v)$ is the same, so that the direction $E_0$ escapes the singularity from $\iota(v)$ as well.

	This is the main ingredient of the next section.
\end{remark}

\section{Description of the horizon}
\label{SHORTSecHorizonSansMatrices}

\subsection{Induction on the dimension}

The horizon in $AdS_3$ is already well understood \cite{Keio,BTZ_horizon}. We are not going to discuss it again. We will study how does the causal structure (black hole, free part, horizon) of $AdS_{l}$ includes itself in $AdS_{l+1}$ by the inclusion map
\begin{equation}
	\iota\colon AdS_l\to AdS_{l+1}.
\end{equation}

\begin{lemma}
	The direction $E_0$ in $AdS_l$ escapes the singularity from $v\in AdS_l$ if and only if it escapes the singularity from $\iota(v)$ in $AdS_{l+1}$.
\end{lemma}

\begin{proof}
	The fact for $v$ to escape the singularity in the direction $E_0$ means that the equation
	\begin{equation}
		a_{v}(E_0)s^2+b_{v}(E_0)+c_{v}=0
	\end{equation}
	where the coefficients are given by \eqref{SHORTEqCoefsabcBE} has no positive solutions with respect to $s$. Since these coefficients are the same for $v$ and $\iota(v)$, the equation for $\iota(v)$ is in fact the same and has the same solutions.
\end{proof}

\begin{lemma}		\label{SHORTLemDueiINtlIntlpu}
	We have
	\begin{equation}
		\iota\big( \Int(\hF_l) \big)\subset\Int(\hF_{l+1})
	\end{equation}
	or, equivalently, 
	\begin{equation}
		\Adh(BH_{l+1})\cap\iota(AdS_l)\subset\iota\big( \Adh(BH_l) \big).
	\end{equation}
\end{lemma}

\begin{proof}
	Let $v'\in\Int(\hF_l)$ and $\mO$, an open set of directions in $AdS_l$ that escape the singularity. The coefficient $a_l(E)$ is not constant on $\mO$ because the coefficient $M=a^2-b^2-C^2$ is only zero on the singularity (see equation \eqref{SHORTEqCoeffaE}). Thus we can choose $E_0\in\mO$ such that $a_l(E_0)\neq 0$. We consider $a_{l+1}(E_0)$, the coefficient of $s^2$ for the point $\iota(v')$ in the direction $E_0$. From the expression \eqref{SHORTEqCoeffaE} we know that $a_{l+1}(E_0)=a_l(E_0)$. The coefficients $b(E_0)$ and $c$ are also the same for $v'$ and $\iota(v')$. 
	
	Since $a(E_0)\neq 0$ and $v'\in\Int(\hF_l)$, we have two solutions to the equation $a(E_0)s^2+b(E_0)s+c=0$ and both of these are outside $\eR^+_0$. This conclusion is valid for $v'\in AdS_l$ as well as for $\iota(v')\in AdS_{l+1}$. Then there is a neighborhood of $\iota(v')$ on which the two solutions keep outside $\eR^+_0$. That proves that $\iota(v')\in\Int(\hF_{l+1})$.

	For the second line, suppose that $v\in\iota(AdS_l)$ does not belong to $\iota\big( \Adh(BH_l) \big)$, thus $v\in\iota\big( \Int(\hF_l) \big)\subset\Int(\hF_{l+1})$. In that case $v$ does not belong to $\Adh(BH_{l+1})$.

\end{proof}

\begin{proposition}		\label{SHORTProphFdanshF}
	We have
	\begin{equation}
		\hF_{l+1}\cap\iota(AdS_l)\subset\iota(\hF_l)
	\end{equation}
\end{proposition}

\begin{proof}
	If $v=\iota(v')\in\hF_{l+1}$, there is a direction $E_0$ in $AdS_{l+1}$ which escape the singularity from $v$. That direction is given by a vector $(w_1,\cdots v_l)\in S^{l}$. Since the coefficients $a(E)$, $b(E)$ and $c$ do only depend on $w_2$, a direction $(w'_1,\cdots,w'_{l-1},0)$ with $w'_2=w_2$ escapes the singularity from $v'$. This proves that $v'\in\hF_l$.
\end{proof}

\begin{lemma}		\label{SHORTLemHiH}
	We have
	\begin{equation}
		\hH_{l+1}\cap\iota(AdS_l)\subset\iota(\hH_{l}).
	\end{equation}
\end{lemma}

\begin{proof}
	First,
	\begin{equation}
		v\in\hH_{l+1}\cap\iota(AdS_l)\subset\hF_{l+1}\cap\iota(AdS_l)\subset\iota(\hF_l)
	\end{equation}
	from proposition \ref{SHORTProphFdanshF}. Now, let's take $v'\in\hF_l$ such that $v=\iota(v')$. We have to prove that $v'\in\hH_l$. Let us suppose that $v'\in\Int(\hF_l)$, so $v\in\Int(\hF_{l+1})$ because of lemma \ref{SHORTLemDueiINtlIntlpu}. This is in contradiction with the fact that $v\in\hH_{l+1}$.
\end{proof}

\begin{corollary}		\label{SHORTCorDeuxTrucsBHhH}
	We have
	\begin{multicols}{2}
		\begin{enumerate}

		\item
			$\iota(\hS_l)\subset\hS_{l+1}$,
		\item
			$\iota(\hF_l)\subset\hF_{l+1}$,
		\item
			$\iota(BH_l)\subset BH_{l+1}$,
		\item
			$\iota(\hH_l)\subset \hH_{l+1}$.

		\end{enumerate}
	\end{multicols}
\end{corollary}

\begin{proof}
	We have $\Ad\big( \iota(g^{-1}) \big)J_1=\Ad(g^{-1})J_1$, so that the condition of theorem \eqref{SHORTThosSequivJzero} is invariant under $\iota$. Thus one immediately has $\iota(\hS_l)\subset\hS_{l+1}$ and $\iota(\hF_l)\subset\hF_{l+1}$.

	An element $v$ which does not belong to $BH_{l+1}$ belongs to $\hF_{l+1}$, but if $v$ belongs to $\iota(AdS_l)\cap\hF_{l+1}$, it belongs to $\iota(\hF_l)$ by proposition \ref{SHORTProphFdanshF} and then does not belong to $\iota(BH_l)$. Thus $\iota(BH_{l})\subset BH_{l+1}$.

	Now if $v'\in\hH_l$, let us consider $\mO$, a neighborhood of $v=\iota(v')$ in $AdS_{l+1}$. The set $\iota^{-1}\big( \mO\cap\iota(AdS_l) \big)$ contains a neighborhood $\mO'$ of $v'$ in $AdS_{l}$. Since $v'\in\hH_l$, there is $\bar v\in\mO'$ such that $\bar v\in BH_l$. Thus $\iota(\bar v)\in\mO$ belongs to $BH_{l+1}$ by the first item.

\end{proof}

Let $X\in\sG_{l+1}$ such that $[X,J_1]=0$, and let $R$ be the group generated by $X$. The following results are intended to show that such a group can be used in order to transport the causal structure from $AdS_l$ to $AdS_{l+1}$

The key ingredient will be the fact that, since $X$ commutes with $J_1$, we have
\begin{equation}		\label{SHORTEqAdOkSurJun}
	\Ad\big( (g e^{sE})^{-1} \big)J_1=\Ad\big( ( e^{\alpha X}g e^{sE})^{-1} \big)J_1.
\end{equation}

\begin{lemma}		\label{SHORTLemRSsubsetS}
	A group $R$ as described above preserves the causal structure in the sense that
	\begin{multicols}{2}
	\begin{enumerate}

		\item
			$R\cdot\hS\subset\hS$
		\item
			$R\cdot BH\subset BH$
		\item
			$R\cdot \hF\subset\hF$
		\item
			$R\cdot \hH\subset\hH$.

	\end{enumerate}
	\end{multicols}
\end{lemma}

\begin{proof}
	From equation \eqref{SHORTEqAdOkSurJun}, we deduce that a direction $E_0$ will escape the singularity from the point $[g]$ is and only if it escapes the singularity from the points $r[g]$ for every $r\in R$. The first three points follow.

	Now let $v\in\hH$ and $r\in R$ and let us prove that $r\cdot v\in BH$. By the third point, $r\cdot v\in\hF$. Let now $\mO$ be a neighborhood of $r\cdot v$. The set $r^{-1}\cdot \mO$ is a neighborhood of $v$ and we can consider $\bar v\in BH\cap r^{-1}\cdot\mO$. By the second point, $r\cdot \bar v$ is a point of the black hole in $\mO$.
\end{proof}

\begin{remark}		\label{SHORTRemdqnqRSlsubsetSlpu}
	Combining corollary \ref{SHORTCorDeuxTrucsBHhH} and lemma \ref{SHORTLemRSsubsetS} we have
	\begin{multicols}{2}
	\begin{enumerate}
		\item
			$ R\cdot\iota(\hS_l)\subset\hS_{l+1}$
		\item
			$ R\cdot\iota(\hF_l)\subset\hF_{l+1}$
		\item
			$ R\cdot\iota(BH_{l})\subset BH_{l+1}$
		\item
			$ R\cdot\iota(\hH_{l})\subset\hH_{l+1}$.
	\end{enumerate}
	\end{multicols}
\end{remark}

\begin{theorem}		\label{SHORTThoCausalPasseParR}
	If moreover the one parameter group $R$ has the property to generate $AdS_{l+1}$ (in the sense that $R\cdot\iota(AdS_l)=AdS_{l+1}$), then we have
	\begin{multicols}{2}
	\begin{enumerate}
		\item
			$ R\cdot\iota(\hS_l)=\hS_{l+1}$
		\item
			$ R\cdot\iota(\hF_l)=\hF_{l+1}$
		\item\label{SHORTItemStrucalpb}
			$ R\cdot \iota(BH_{l})=BH_{l+1}$
		\item
			$ R\cdot\iota(\hH_{l})=\hH_{l+1}$.
	\end{enumerate}
	\end{multicols}
\end{theorem}

\begin{proof}
	The inclusions in the direct sense are already done in the remark \ref{SHORTRemdqnqRSlsubsetSlpu}.

	Let $r= e^{\alpha X}$ be an element of $R$. Since, by assumption, we have $[X,J_1]=0$, the action of $r$ leaves invariant the condition of theorem \eqref{SHORTThosSequivJzero}:
	\begin{equation}		\label{SHORTEqalpharjnagitpas}
		\Ad\big( (g e^{sE})^{-1} \big)J_1=\Ad\big( ( e^{\alpha X}g e^{sE})^{-1} \big)J_1.
	\end{equation}
	\begin{enumerate}
		\item
			Let $[g]\in\hS_{l+1}$, there exists a $r\in R$ such that $r[g]\in \iota(AdS_l)$. There exists an element $g'\in G_l$ such that $rg=\iota(g')$. Now $[g']\in\hS_l$ because
			\begin{equation}
				\Ad(g'^{-1})J_1=\Ad(\iota(g'^{-1}))J_1=\Ad\big(  (rg)^{-1} \big)J_1=\Ad(g^{-1})J_1,
			\end{equation}
			but by assumption the norm of the projection on $\sQ$ of the right hand side is zero.

		\item
			If $v$ is free in $AdS_{l+1}$, there is a direction $E_0$ escaping the singularity from $v$ and an element $r\in R$ such that $v'=r\cdot v\in\iota(AdS_l)$. The point $v'$ is also free in $AdS_{l+1}$ as the direction $E_0$ works for $r\cdot v$ as well as for $v$. Thus by proposition \ref{SHORTProphFdanshF} we have
			\begin{equation}
				v'\in\hF_{l+1}\cap\iota(AdS_l)\subset\iota(\hF_l)
			\end{equation}
			and $v\in R\cdot\iota(\hF_l)$.
		\item
			If $v\in BH_{l+1}$, the point $r\cdot v\in\iota(AdS_l)$ also belongs to $BH_{l+1}$. If $r\cdot v=\iota(v')$, then $v'\in BH_l$ from if $v'\in\hF_l$, then $\iota(v')\in\hF_{l+1}$.
		\item
			If $v\in\hH_{l+1}$, there exists a $r\in R$ such that $v'=r\cdot v\in\iota(AdS_l)$ and moreover, $v'$ belongs to the horizon in $AdS_{l+1}$ since the horizon is invariant under $R$. Thus $v'$ belongs to $\iota(AdS_l)\cap\hH_{l+1}\subset \iota(\hH_l)$ by lemma \ref{SHORTLemHiH}.  Now, $v\in R\cdot\iota(\hH_l)$.

	\end{enumerate}
\end{proof}

\subsection{Examples of surjective groups}

Theorem \ref{SHORTThoCausalPasseParR} describes the causal structure in $AdS_l$ by induction on the dimension provided that one knows a group $R$ such that $AdS_{l+1}=R\cdot \iota(AdS_l)$. Can one provide examples of such groups? The following proposition provides a one.

\begin{proposition}		\label{SHORTPropSurjectif}
	If $R$ is the one parameter subgroup of $\SO(2,l)$ generated by $r_{l,l+1}$, then we have
	\begin{equation}
		R\cdot \iota(AdS_l)= AdS_{l+1}.
	\end{equation}
\end{proposition}

\begin{proof}
	If one realises $AdS_n$ as the set of vectors of length $1$ in $\eR^{2,n-1}$, $AdS_l$ is included in $AdS_{l+1}$ as the set of vectors with vanishing last component and the element $r_{l,l+1}$ is the rotation in the plane of the two last coordinates. In that case, we have to solve
	\begin{equation}
		e^{\alpha r_{l,l+1}}\begin{pmatrix}
			u'	\\ 
			t'	\\ 
			x'_1	\\ 
			\vdots	\\ 
			x'_{l-2}	\\ 
			0	
		\end{pmatrix}=
		\begin{pmatrix}
			u	\\ 
			t	\\ 
			x_1	\\ 
			\vdots	\\ 
			x_{l-2}	\\ 
			x_{l-1}	
		\end{pmatrix}
	\end{equation}
	with respect to $\alpha$, $u'$, $t'$ and $x'_i$. There are of course exactly two solution if $x_{l-2}^2+x_{l-1}^2\neq 0$.
\end{proof}

In fact, many others are available, as the one showed at the end of \cite{BTZ_horizon}. In fact, since, in the embedding of $AdS$ in $\eR^{2,n}$, the singularity is given by $t^2-y^2=0$, almost every group which leaves invariant the combination $t^2-y^2$ can be used to propagate the causal structure. One can found lot of them for example by looking at the matrices given in \cite{These}.

\subsection{Backward induction}

Using proposition \ref{SHORTProphFdanshF}, lemma \ref{SHORTLemHiH}, corollary \ref{SHORTCorDeuxTrucsBHhH} and the fact that the norm of $J_1^*$ is the same in $AdS_l$ as in $\iota(AdS_l)\subset AdS_{l+1}$, we have
\begin{enumerate}
	\item
		$\iota(\hF_l)=\hF_{l+1}\cap\iota(AdS_l)$,
	\item
		$\iota(\hH_l)=\hH_{l+1}\cap\iota(AdS_l)$,
	\item 
		$\iota(\hS_l)=\hS_{l+1}\cap\iota(AdS_l)$.
\end{enumerate}
These equalities hold for $l\geq 3$. For $l=2$ we can take the latter as a definition and set
\begin{equation}
	\hS_2=\{ v\in AdS_2\tq \iota(v)\in\hS_3 \}.
\end{equation}
The Iwasawa decomposition of $\SO(2,1)$ is given by $\sA=\langle J_2\rangle$, $\sN=\langle X_+\rangle$, $\sK=\langle q_0\rangle$ where $X_+=p_1-q_0$. Notice that we \emph{do not} have $\iota(X_+)=X_{++}$. Instead we have $X_+=\frac{ 1 }{2}(X_{++}+X_{-+})$. Thus $\hS_2$ is not given by the closed orbits of $AN$ in $AdS_2$.

The light-like directions are given by the two vectors $E=q_0\pm q_1$. In order to determine if the point $[e^{-\alpha J_2}e^{-aX_+} e^{-xq_0}]$ belongs to the black hole, we follow the same way as in section \ref{SHORTSubSecExistenceTrouNoir}: we compute the norm
\begin{equation}
	\big\| \pr_{\sQ}  e^{-s\ad(E)}e^{x\ad(q_0)} e^{a\ad(X_+)} e^{\alpha\ad(J_2)}J_1 \big\|
\end{equation}
and we see under which conditions it vanishes. Small computations show that, for $E=q_0+q_1$,

\begin{equation}
	\pr_{\sQ} e^{s\ad(E)}X=\Big( \big( \cos(x)-a\sin(x)+a \big)s+(\sin(x)+a\cos(x)) \Big)q_2.
\end{equation}
Its norm vanishes for the value of $s$ given by
\begin{equation}
	s^+=\frac{ a\cos(x)+\sin(x) }{ \big( \sin(x)-1 \big)a-\cos(x) }.
\end{equation}
The same computation with $E=q_0-q_1$ provides the value
\begin{equation}
	s^-=\frac{ a\cos(x)+\sin(x) }{ \big( \sin(x)+1 \big)a-\cos(x) }.
\end{equation}

For small enough $a$, the sign of $s^+$ and $s^-$ are both given the sign of $-\tan(x)$ that can be either positive or negative. Thus there is an open set of points in $AdS_2$ which intersect the singularity in every direction and an open set of points which escape the singularity.

\chapter{Long version}
\section{Introduction}
\label{LONGSecSumStructExist}

\subsection{Anti de Sitter space and the BTZ black hole}

The anti de Sitter space (hereafter abbreviated by $AdS$, or $AdS_l$ when we refer to a precise dimension) is a static solution to the Einstein's equations that describes a universe without mass. It is widely studied in different context in mathematics as well as in physics.

The BTZ black hole, initially introduced in \cite{BTZ_un,BTZ_deux} and then described and extended in various ways \cite{HolstPeldan,Aminneborg,Madden}, is an example of black hole structure which does not derives from a metric singularity. 

The structure of the BTZ black hole as we consider it here grown from the papers \cite{BTZB_deux,Keio} in the case of $AdS_3$. Our dimensional generalization was first performed in \cite{lcTNAdS}. See also \cite{These} for for a longer review. Our point of view insists on the homogeneous space structure and the action of Iwasawa groups. One of the motivation in going that way is to embed the study of BTZ black hole into the noncommutative geometry and singleton physics \cite{BTZ_WZW,articleBVCS}. 

\subsection{The way we describe the BTZ black hole}

We look at the anti de Sitter space as the homogeneous space
\begin{equation}
	AdS_l=\frac{ \SO(2,l-1) }{ \SO(1,l-1) }=G/H.
\end{equation}
We denote by $\sG=\so(2,l-1)$ and $\sH=\so(1,l-1)$ the Lie algebras and by $\pi$ the projection $G\to G/H$. The class of $g$ will be written $[g]$ or $\pi(g)$. We choose an involutive automorphism $\sigma\colon \sG\to \sG$ which fixes elements of $\sH$, and we call $\sQ$ the eigenspace of eigenvalue $-1$ of $\sigma$. Thus we have the reductive decomposition
\begin{equation}		\label{LONGEqIntroRedDecompHQLieAlg}
	\sG=\sH\oplus\sQ.
\end{equation}
The compact part of $\SO(2,l-1)$ decomposes into $K=\SO(2)\times\SO(l-1)$.

Let $\theta$ be a Cartan involution which commutes with $\sigma$, and consider the corresponding Cartan decomposition
\begin{equation}
	\sG=\sK\oplus\sP,
\end{equation}
where $\sK$ is the $+1$ eigenspace of $\theta$ and $\sP$ is the $-1$ eigenspace. A maximal abelian algebra $\sA$ in $\sP$ has dimension two and one can choose a basis $\{ J_1,J_2 \}$ of $\sA$ in such a way that $J_1\in\sH$ and $J_2\in\sQ$.

Now we consider an Iwasawa decomposition
\begin{equation}
	\sG=\sK\oplus\sA\oplus\sN,
\end{equation}
and we denote by $\sR$ the Iwasawa component $\sR=\sA\oplus\sN$. We are also going to use the algebra $\bar\sN=\theta\sN$ and the corresponding Iwasawa component $\bar\sR=\sA\oplus\bar\sN$.

The Iwasawa groups $R=AN$ and $\bar R=A\bar N$ are naturally acting on anti de Sitter by $r[g]=[rg]$. It turns out that each of these two action has exactly two closed orbits, regardless to the dimension we are looking at. The first one is the orbit of the identity and the second one is the orbit of $[k_{\theta}]$ where $k_{\theta}$ is the element which generates the Cartan involution at the group level: $\AD(k_{\theta})=\theta$. In a suitable choice of matrix representation, the element $k_{\theta}$ is the block-diagonal element which is $-\mtu$ on $\SO(2)$ and $\mtu$ on $\SO(l-1)$. The $A\bar N$-orbits of $\mtu$ and $k_{\theta}$ are also closed. Moreover we have
\begin{equation}
	\begin{aligned}[]
		[A\bar N k_{\theta}]&=[k_{\theta}AN]\\
		[AN k_{\theta}]&=[k_{\theta}A\bar N]
	\end{aligned}
\end{equation}
because $A$ is invariant under $\AD(k_{\theta})$ and, by definition, $\Ad(k_{\theta})N=\bar N$. We define as \defe{singular}{Singular point} the points of the closed orbits of $AN$ and $A\bar N$ in $AdS$.

The Killing form of $\SO(2,l-1)$ induces a Lorentzian metric on $AdS$. The sign of the squared norm of a vector thus divides the vectors into three classes:
\begin{equation}
	\begin{aligned}[]
		\| X \|^2&>0&\rightarrow&&\text{time like,}\\
		\| X \|^2&<0&\rightarrow&&\text{space like,}\\
		\| X \|^2&=0&\rightarrow&&\text{light like.}
	\end{aligned}
\end{equation}
A geodesic is time (reps. space, light) like if its tangent vector is time like (reps. space, light).

If $E_1$ is a nilpotent element in $\sQ$, then every nilpotent in $\sQ$ are given by $\{ \Ad(k)E_1 \}_{k\in\SO(l-1)}$. These elements are also all the light like vectors at the base point. A light like geodesic trough the point $\pi(g)$ in the direction $\Ad(k)E_1$ is given by
\begin{equation}
	\pi(g e^{s\Ad(k)E_1}).
\end{equation}
One say that  points with $s>0$ are in the \defe{future}{Future} of $\pi(g)$ while points with $s<0$ are in the \defe{past}{Past} of $\pi(g)$. 

We say that a point in $AdS_l$ belongs to the \defe{black hole}{Black hole} if all the light like geodesics trough that point intersect the singularity in the future. We call \defe{horizon}{Horizon} the boundary of the set of points in the black hole. One say that there is a (non trivial) black hole structure when the horizon is non empty or, equivalently, when there are some points in the black hole, and some outside.

All these properties can be easily checked using the matrices given in \cite{These,lcTNAdS}. In this optic, I wrote a program using Sage\cite{Sage} which checks all the properties that are shown in this paper. It will be published soon.

As far as notations are concerned, we denote by $X_{\alpha\beta}$ the basis of $\sN$ and $\bar\sN$ corresponding to our choice of Iwasawa decomposition. We have $\ad(J_1)X_{\alpha\beta}=\alpha X_{\alpha\beta}$ and $\ad(J_2)X_{\alpha\beta}=\beta X_{\alpha\beta}$.

\subsection{Organization of the paper}

One of the main goal of this paper is to reorganize all this structure in a coherent way. Then we use it efficiently in order to define the singularity of the BTZ black hole, to prove that one has a genuine black hole in every dimension, and to determine the horizon.

In section \ref{LONGSecProgressRidMatrices}, we list the commutators of $\so(2,n)$ with respect to its root spaces and we organize them in such a way to get a clear idea about the evolution of the structure when the dimension increases. We prove that, when one passes from $\so(2,n)$ to $\so(2,n+1)$, one gets four more vectors in the root spaces and that these are Killing-orthogonal to the vectors existing in $\so(2,n)$ (this is the ``dimensional slice'' described in subsection \ref{LONGSubSecDimensionalSlices}).

We give in subsection \ref{LONGSubSecReductiveDecompQ} an original way to describe the space $\sQ$ without reference to $\sH$. The space $\sQ$ is usually described as a complementary of $\sH$. Here we show that it can be described by means of the root spaces and the Cartan involution $\theta$. The space $\sH$ is then described as $\sH=[\sQ,\sQ]$. In some sense, we describe the quotient space $AdS=G/H$ directly by its tangent space $\sQ$ without passing trough the definition of $H$. Of course, the knowledge of $\sH$ will be of crucial importance later.

The subsection \ref{LONGSubSecPropRedDecompQH} is devoted to the proof of many properties of the decompositions $\sG=\sH\oplus\sQ$ and $\sG=\sK\oplus\sP$.

The first important result is the proposition \ref{LONGXUnALaTwistingSuperCool} that shows that the elements of $\sQ$ are $\ad$-conjugate to each others: there exist elements of the adjoint group which are intertwining the elements of $\sQ$. We also provide an orthonormal basis $\{ q_i \}$ of $\sQ$, we compute the norm of these elements and we identify the nilpotent vectors in $\sQ$ (these are the light-like vectors). In the same time, we prove that the space $G/H$ is Lorentzian.

The second central result is the fact that nilpotent elements in $\sQ$ are of order two: if $E\in\sQ$ is nilpotent, then $\ad(E)^3=0$. That result will be used in a crucial way in the proof of the black hole existence, as well as in the study of its properties.

In section \ref{LONGSecBlacHole}, we define and study the structure of the BTZ black hole in the anti de Sitter space. First we identify the closed orbits of the Iwasawa group (theorem \ref{LONGThoOrbitesOuverttes}) and we define them as singular. In a second time, we provide an alternative description the singularity: theorem \ref{LONGThosSequivJzero} shows that the singularity can be described as the loci of points at which a fundamental vector field has vanishing norm. We also provide in lemma \ref{LONGLemExpressionCoolNormJUn} a convenient way to compute that norm on arbitrary point of the space.

We prove, in section \ref{LONGSubSecExistenceTrouNoir}, that our definition of singularity gives rise to a genuine black hole in the sense that there exists points from which some geodesics escape the singularity in the future and there exists some points from which all the geodesics are intersecting the singularity in the future.

In section \ref{LONGSecHorizonSansMatrices}, we provide a geometric description of the horizon (theorem \ref{LONGThoCausalPasseParR}). We show that, if we see $AdS_l$ as a subset of $AdS_{l+1}$, the space $AdS_{l+1}$ is generated by the action on $AdS_l$ by a one parameter group which leaves the singularity invariant. This action thus leaves invariant the whole causal structure and we are able to express the horizon in any dimension from the well known horizon in $AdS_3$.

\section{Structure of the algebra}
\label{LONGSecProgressRidMatrices}

\subsection{The Iwasawa component}

Our study of $AdS_l=\SO(2,l-1)/\SO(1,l-1)$ will be based on the properties of the algebra $\sG=\so(2,l-1)$ endowed with a Cartan involution $\theta$ and an Iwasawa decomposition $\sG=\sA\oplus\sN\oplus\sK$. In this section we want to underline the most relevant facts for our purpose. The part we are mainly interested in is the Iwasawa component $\sA\oplus\sN$ where
\begin{subequations}		\label{LONGEqLeANEnDimAlg}
\begin{align}
	\sN&=\{ X^{k}_{+0},X^{k}_{0+},X_{++},X_{+-} \}\\
	\sA&=\{ J_1, J_2\},
\end{align}
\end{subequations}
where $k$ runs %
\footnote{The ``new'' vectors which appear in $AdS_l$ with respect to $AdS_{l-1}$ are $X_{0\pm}^{l-1}$ and $X_{\pm 0}^{l-1}$. Such an element appears for the first time in $AdS_4$ and is not present when one study $AdS_3$.} %
from $3$ to $l-1$. The commutator table is
\begin{subequations}  \label{LONGEqTableSOIwa}
	\begin{align}
		[X_{0+}^{k},X_{+0}^{k'}]&=\delta_{kk'}X_{++}		&[X_{0+}^{k},X_{+-}]&=2X_{+0}^{k}\\
		[ J_1,X_{+0}^{k}]&=X_{+0}^{k}				&[ J_2,X_{0+}^{k}]&=X_{0+}^{k}\\
		[ J_1,X_{+-}]&=X_{+-}					&[ J_2,X_{+-}]&=-X_{+-}\\
		[ J_1,X_{++}]&=X_{++}					&[ J_2,X_{++}]&=X_{++}.
	\end{align}
\end{subequations}
We see that the Iwasawa algebra belongs to the class of $j$-algebras whose Pyatetskii-Shapiro decomposition is 
\begin{equation}
	\sA\oplus\sN=(\sA_1\oplus_{\ad}\sZ_1)\oplus_{\ad}\big( \sA_2\oplus_{\ad}(V\oplus \sZ_2) \big),
\end{equation}
with
\begin{subequations}
	\begin{align}
		\sA_1&=\langle H_1\rangle		&\sA_2&=\langle H_2\rangle\\
		\sZ_1&=\langle X_{+-}\rangle		&\sZ_2&=\langle X_{++}\rangle\\
		&					&V&=\langle X_{0+}^{k},X_{+0}^{k}\rangle_{k\geq 3}
	\end{align}
\end{subequations}
where 
\begin{equation}
	\begin{aligned}[]
		H_1&=J_1-J_2\\
		H_2&=J_1+J_2.
	\end{aligned}
\end{equation}

The general commutators of such an algebra are
\begin{subequations}		\label{LONGsubEqsGenPySO}
\begin{align}
	[H_1,X_{+-}]		&=2X_{+-}	&	[H_2,X_{0+}^{k}]	&=X_{0+}^{k}				&[H_1,V]	&\subset V	\\		
				&		&	[H_2,X_{+0}^{k}]	&=X_{+0}^{k}				&[X_{+-},V]	&\subset V	\\	
				&		&	[H_2,X_{++}]		&=2X_{++}								\\
				&		&	[X_{0+}^k,X_{+0}^l]	&=\Omega(X_{0+}^k,X_{+0}^l)X_{++}
\end{align}
\end{subequations}
In the case of $\so(2,n)$, we have the following more precise relations:
\begin{subequations}		\label{LONGSubEqsPlusPresPySO}
	\begin{align}
		[H_1,X^k_{0+}]&=-X^k_{0+}\\		
		[X_{+-},X^k_{0+}]&=-2X^k_{+0}
	\end{align}
\end{subequations}
and the link between $\sN$ and $\tilde\sN$ is given by
\begin{subequations}		\label{LONGSubEqsThethaPySO}
	\begin{align}
		[\theta X_{+0}^{k},X_{++}]	&=2X_{0+}^{k}\\
		[\theta X_{0+}^{k},X^k_{0+}]	&=2J_2\\
		[\theta X_{++},X_{++}]		&=4H_2=4(J_1+J_2)\\
		[\theta X_{++},X^{k}_{0+}]	&=2X_{-0}^{k}\\
		[\theta X_{+-},X_{+-}]		&=4H_1=4(J_1-J_2)\\
		[\theta X_{+-},X^{k}_{+0}]	&=2X_{0+}^{k}.
	\end{align}
\end{subequations}
The relations between the higher dimensional root spaces are
\begin{equation}
	\begin{aligned}[]
		[X_{0+}^i,X_{-0}^j]&=-\delta_{ij}X_{-+}\\
		[X_{0+}^i,X_{+0}^j]&=\delta_{ij}X_{++}\\
		[X_{+0}^i,X_{0+}^j]&=-\delta_{ij}X_{++}\\
		[X_{+0}^i,X_{0-}^j]&=\delta_{ij}X_{+-}.
	\end{aligned}
\end{equation}
The space $\mZ_{\sK}(\sA)$ of the elements of $\sK$ which commute with all the elements of $\sA$ is given by the elements $r_{ij}=\frac{ 1 }{2}[X_{0+}^i,X_{0+}^j]$ for every $i,j\geq 3$, $i\neq j$ realise the $\so(l-1)$ algebra. We will say more about them in subsection \ref{LONGSubSec_Thecompactpart}.

We deduce the following relations that will prove useful later
\begin{equation}
	\begin{aligned}[]
		[\theta X_{0+}^k,X_{++}]&=-2X^k_{+0}\\
		[\theta X_{+0}^k,X_{+-}]&=-2X^k_{0-}\\
		[\theta X_{++},X_{+0}^k]&=-2X^k_{0-}\\
		[X_{-+},X^k_{0-}]&=-2X^k_{-0}.
	\end{aligned}
\end{equation}

\subsection{The compact part}
\label{LONGSubSec_Thecompactpart}

The compact part of $\so(2,l-1)$, the algebra $\so(2)\oplus\so(l-1)$ is well known. What is interesting from our point of view is to write the commutation relations between the elements of $\so(l-1)$ and the roots.

We define the following elements that are non vanishing:
\begin{equation}
	[X_{0+}^i,X_{0-}^j]=[X_{0-}^i,X_{0+}^j]=2r_{ij}.
\end{equation}
One immediately has $\theta r_{ij}=r_{ij}$, so that $r_{ij}\in\sK$. We also have $r_{ij}\in\sG_0$ so that $r_{ij}\in\mZ_{\sK}(\sA)$ and they act on the root spaces. The action is given by
\begin{subequations}
	\begin{align}		\label{LONGEqComsRRN}
		[r_{ij},X_{+0}^k]&=0		&\text{if $i$, $j$ and $k$ are different}\\
		[r_{ij},X_{0+}^k]&=0		&\text{if $i$, $j$ and $k$ are different}\\
		[r_{ij},X_{+0}^j]&=X_{+0}^i	&\text{if $i\neq j$}\\
		[r_{ij},X_{0+}^j]&=-X_{0+}^i	&\text{if $i\neq j$}\\
		[r_{ij},X_{\pm\pm}]&=0.			\label{LONGsubeqrXpmpm}
	\end{align}
\end{subequations}
The elements $r_{ij}$ satisfy the algebra of $\so(n)$.

\begin{remark}
	If $\sigma$ is an involutive automorphism which commutes with $\theta$ and such that $\sigma J_1=J_1$, $\sigma J_2=-J_2$, then one has $\sigma r_{ij}=r_{ij}$. We will see later that this fact makes $r_{ij}\in\sH$.
\end{remark}

We know\footnote{See \cite{Berndt} for example.} that $\sG=\sG_0\oplus\sN\oplus\bar\sN$ where $\sG_0=\sA\oplus\mZ_{\sK}(\sA)$. Let us perform a dimension count in order to be sure that the vectors $r_{ij}$ generate $\mZ_{\sK}(\sA)$. When we are working with $AdS_l$, we have
\begin{equation}
	\begin{aligned}[]
		\dim(\sA)&=2\\
		\dim(\tilde\sN_2)&=4\\
		\dim(\bigoplus_{k}\tilde\sN_k)&=4(l-3)\\
		\dim(\langle r_{ij}\rangle)&=\frac{ 1 }{2}(l-4)(l-3).
	\end{aligned}
\end{equation}
The last line comes from the fact that we have the elements $r_{34}$, $r_{35}$,\ldots,$r_{45}$,\ldots The first such element appears in $AdS_5$. Making the sum, we obtain $\frac{ l(l+1) }{ 2 }$, which is the dimension of $\so(2,l-1)$. Thus we have
\begin{equation}
	\sG=\mZ_{\sK}(\sA)\oplus\sA\oplus\tilde\sN_2\oplus\bigoplus_{k}\tilde\sN_k
\end{equation}
where $\mZ_{\sK}(\sA)$ is generated by the elements $r_{ij}$.
\subsection{Dimensional slices}
\label{LONGSubSecDimensionalSlices}

Since the elements of $\mZ_{\sK}(\sA)$ is a part of $\sH$ (see later), they will have almost no importance in the remaining\footnote{We will however need them in the computation of the coefficients \eqref{LONGEqCoefsabcBE}.}. The most important part is
\begin{equation}		\label{LONGEqDecomDimQlipm}
	\sA\oplus\sN\oplus\bar\sN=\underbrace{\langle J_1,J_2,X_{\pm,\pm}\rangle}_{\text{for every dimension}}\oplus\underbrace{\langle X_{0\pm}^{4},X_{\pm 0}^{4}\rangle}_{\text{for $\so(2,\geq 3)$}}\oplus\ldots\oplus\underbrace{\langle X_{0\pm}^l,X_{\pm 0}^{l}\rangle}_{\text{for $\so(2, l-1)$}}.
\end{equation}
We use the following notations in order to make more clear how does the algebra evolve when one increases the dimension:
\begin{equation}
	\begin{aligned}[]
		\sN_2&=\langle X_{+-},X_{++}\rangle,		&\sN_k&=\langle X^k_{0+},X_{+0}^k\rangle	\\
		\bar\sN_2&=\langle X_{-+},X_{--} \rangle, 	&\bar\sN_k&=\langle X_{0-}^k,X_{-0}^k\rangle\\
		\tilde\sN_2&=\langle \sN_2,\bar\sN_2\rangle,	&\tilde\sN_k&=\langle \sN_k,\bar\sN_k\rangle
	\end{aligned}
\end{equation}
for $k\geq 3$. The relations are
\begin{equation}		\label{LONGEqsCommWithtsNDeuxkA}
	\begin{aligned}[]
		[\tilde\sN_2,\tilde\sN_2]&\subseteq\sA
		&[\tilde\sN_2,\tilde\sN_k]&\subseteq\tilde\sN_k\\
		[\tilde\sN_k,\tilde\sN_{k}]&\subseteq\sA\oplus\tilde\sN_2
		&[\tilde\sN_k,\tilde\sN_{k'}]&\subset \mZ_{\sK}(\sA) \\
	\end{aligned}
\end{equation}

As a consequence of the splitting and the commutation relations, we have many Killing-orthogonal subspaces in $\so(2,l-1)$ :
\begin{subequations}
	\begin{align}
		\tilde\sN_k\perp\tilde\sN_{k'}		\label{LONGEqtsnkperprtsnkp}	\\
		\sA\perp\tilde\sN_2			\label{LONGEqAperpNTrois}		\\
		\sA\perp\tilde\sN_k.	\label{LONGEqAperpNk}				\\
		\tilde\sN_2\perp\tilde\sN_k	\label{LONGEqNTtroisperpNk}
	\end{align}
\end{subequations}

In order to check there relations, look at the trace of action of $\ad(X)\circ\ad(Y)$ on the various spaces :
\begin{equation}
	\ad(\sA)\circ\ad(\tilde\sN_2) \colon 
	\begin{cases}
		\sA\to\tilde\sN_2\to\tilde\sN_2\\
		\tilde\sN_2\to\sA\to 0\\
		\tilde\sN_k\to\tilde\sN_k\to\tilde\sN_k\\
		\mZ_{\sK}(\sA)\to\tilde\sN_2\to\tilde\sN_2,
	\end{cases}
\end{equation}
and
\begin{equation}
	\ad(\sA)\circ\ad(\tilde\sN_k)\colon
	\begin{cases}
		\sA\to\tilde\sN_k\to\tilde\sN_k\\
		\tilde\sN_2\to\tilde\sN_k\to\tilde\sN_k\\
		\tilde\sN_k\to\sA\oplus\tilde\sN_2\to\tilde\sN_2\\
		\langle r_{kl}\rangle\to\tilde\sN_l\to\tilde\sN_l
	\end{cases}
\end{equation}
and
\begin{equation}
	\ad(\tilde\sN_2)\circ\ad(\tilde\sN_k)\colon
	\begin{cases}
		\sA\to\tilde\sN_k\to\tilde\sN_k\\
		\tilde\sN_2\to\tilde\sN_k\to\tilde\sN_k\\
		\tilde\sN_k\to\sA\oplus\tilde\sN_2\to\tilde\sN_2\oplus\sA\\
		\langle r_{kl}\rangle\to\tilde\sN_l\to\tilde\sN_l\oplus\sA
	\end{cases}
\end{equation}
and
\begin{equation}
	\ad(J_1)^2|_{\tilde\sN_2}=\ad(J_2)^2|_{\tilde\sN_2}=\id|_{\tilde\sN_2}.
\end{equation}
As a consequence,
\begin{equation}		\label{LONGEqDecompGPourKillOrtho}
	\sG=\mZ_{\sK}(\sA)\oplus\sA\oplus\tilde\sN_2\bigoplus_{k\geq 3}\tilde\sN_k
\end{equation}
is a Killing-orthogonal decomposition of $\so(2,l-1)$.

\subsection{Reductive decomposition}
\label{LONGSubSecReductiveDecompQ}

Let $\sQ$ be the following vector subspace of $\sG$:
\begin{equation}		\label{LONGEqDecQEspacesCools}	%
	\sQ=\big\langle \mZ(\sK),J_2,[\mZ(\sK),J_1],(X_{0+}^k)_\sP\big\rangle_{k\geq 3}.
\end{equation}
Then we choice a subalgebra $\sH$ of $\sG$ which, as vector space, is a complementary of $\sQ$. In that choice, we require that there exists an involutive automorphism $\sigma\colon \sG\to \sG$ such that
\begin{equation}
	\sigma=(\id)_{\sH}\oplus(-\id)_{\sQ}.
\end{equation}
In that case the decomposition $\sG=\sH\oplus\sQ$ is reductive, i.e. $[\sQ,\sQ]\subset\sH$ and $[\sQ,\sH]\subset\sQ$.

From definition \eqref{LONGEqDecQEspacesCools}, it is immediately apparent that one has a basis of $\sQ$ made of elements in $\sK$ and $\sP$, so that one immediately has
\begin{equation}
	[\sigma,\theta]=0.
\end{equation}

If $X\in\sG$, the projections are given by
\begin{equation}		\label{LONGEqProjHQPKsigmatheta}
	\begin{aligned}[]
		X_{\sH}&=\frac{ 1 }{2}(X+\sigma X),		&X_{\sK}&=\frac{ 1 }{2}(X+\theta X),\\
		X_{\sQ}&=\frac{ 1 }{2}(X-\sigma X),		&X_{\sP}&=\frac{ 1 }{2}(X-\theta X).
	\end{aligned}
\end{equation}
In particular $\theta\sH\subset\sH$ since $\theta$ and $\sigma$ commute.

We introduce the following elements of $\sQ$:
\begin{equation}			\label{LONGAlignPRedDefQQ}
	\begin{aligned}[]
		q_0&=(X_{++})_{\sK\sQ}\\
		q_1&=J_2\\
		q_2&=-[J_1,q_0]=-(X_{++})_{\sP\sQ}\\
		q_k&=(X^k_{0+})_{\sP}.
	\end{aligned}
\end{equation}
In order to see that $[J_1,q_0]=(X_{++})_{\sP\sQ}$, just write the $\sP\sQ$-component of the equality $[J_1,X_{++}]=X_{++}$. We will prove later that this is a basis and that each of these elements correspond to one of the spaces listed in \eqref{LONGEqDecQEspacesCools}. 

Since $X_{\sP}=(\theta X-X)/2$, we have
\begin{equation}
	[q_i,q_j]=-\frac{1}{ 4 }\big( [X_{0+}^i,X_{0-}^j]+[X_{0-}^i,X_{0+}^j] \big)=r_{ij}.
\end{equation}%

From equations \eqref{LONGEqAperpNTrois} and \eqref{LONGEqAperpNk}, we have $q_1\perp q_2$ and $q_2\perp q_k$. Using the other perpendicularity relations $\sK\perp\sP$ and \eqref{LONGEqtsnkperprtsnkp}, \eqref{LONGEqAperpNTrois}, \eqref{LONGEqAperpNk}, 
we see that the set $\{ q_i \}_{0\leq i\leq l-1}$ is orthogonal.

The space $\sH$ is defined as generated by the elements
\begin{equation}			\label{LONGAlignPremDefHH}
	\begin{aligned}[]
		J_1&		&r_k&=[J_2,q_k]			\\
		p_1&=[q_0,q_1]	&p_k&=[q_0,q_k]			\\
		s_1&=[J_1,p_1]	&s_k&=[J_1,p_k].
	\end{aligned}
\end{equation}
Elements \eqref{LONGAlignPRedDefQQ} and \eqref{LONGAlignPremDefHH} will be studied in great details later.%

\subsubsection{Remark on the compact part}
\label{LONGSubSubSecRemCompPart}

Elements of $\sK$ are elements of the form $X+\theta X$. A part of the elements inside $\mZ_{\sK}(\sA)$, these elements are of two kinds:
\begin{subequations}
	\begin{align}
		X_{++}+X_{--}\\
		X_{+-}+X_{-+}
	\end{align}
\end{subequations}
on the one hand, and
\begin{subequations}
	\begin{align}
		X^k_{0+}+X^k_{0-}\\
		X^k_{+0}+X^k_{-0}
	\end{align}
\end{subequations}
on the other hand. The first two are commuting, so that $\mZ(\sK)$ is two dimensional when one study $AdS_3$. That correspond to the well known fact that the compact part of $\so(2,2)$ is $\so(2)\oplus\so(2)$ which is abelian. These elements, however, do not commute with the two other. For example, the combination
\begin{equation}
	\big( X_{++}+X_{--}\big)-\big(X_{+-}+X_{-+}\big)
\end{equation}
does not commute with the elements of the second type. Now, one checks that the combination
\begin{equation}		\label{LONGEqZKChoixblabl}
	\big( X_{++}+X_{--}\big)+\big( X_{+-}+X_{-+})
\end{equation}
commutes with all the other, so that it is the generator of $\mZ(\sK)$ for $AdS_{\geq 4}$. This corresponds to the fact that the compact part of $\so(2,n)$ is $\so(2)\oplus\so(n)$. In other terms,
\begin{equation}		\label{LONGEqDecompZKenDeuxSuivantDim}
	\mZ(\sK)=\langle X_{++}+X_{--}+X_{+-}+X_{-+}\rangle\oplus\underbrace{\langle X_{++}+X_{--}-X_{+-}-X_{-+}\rangle}_{\text{only for $AdS_3$}}.
\end{equation}

Notice that, for $AdS_{\geq 4}$, we can define $q_0=(X_{++})_{\mZ(\sK)}$ as $\sK=\so(2)\oplus\so(l-2)$ for $AdS_l$. The case of $AdS_3$ is particular because $\mZ(\sK)$ is of dimension two and we have to set by hand what part of $\mZ(\sK)$ belongs to $\sQ$ (the other part belongs to $\sH$). From what is said around equation \eqref{LONGEqZKChoixblabl}, we know that $q_0$ is a multiple of $X_{++}+X_{--}+X_{+-}+X_{-+}$.

Dimension counting shows that $\dim\sQ=l$ and general theory of homogeneous spaces shows that $\sQ$ has to be seen as the tangent space of the manifold $G/H$.

\subsection{Properties of the reductive decompositions}
\label{LONGSubSecPropRedDecompQH}

We are considering the two reductive decompositions $\sG=\sH\oplus\sQ=\sK\oplus\sP$ in the same time as the root space decomposition \eqref{LONGEqDecompGPourKillOrtho}. We are now giving some properties of them.

We know that $\sK\cap\sQ=\langle q_0\rangle$ belongs to $\tilde\sN_2$. As a consequence, the elements $X^k_{\alpha 0}$ and $X^k_{0\alpha}$ have no $\sK\sQ$-components and 
\begin{equation}
	\begin{aligned}[]
		[J_2,\tilde\sN_k]_{\sP\sH}&=0\\
		[J_2,\tilde\sN_k]_{\sP\sQ}&=0\\
		\pr_{\sP\sQ}X_{\alpha 0}^k&=0\\
		\pr_{\sP\sH}X_{0\alpha }^k&=0.
	\end{aligned}
\end{equation}
Since $X^k_{\alpha 0}$ and $X^k_{0\alpha}$ are not  eigenvectors of $\theta$, they have a non vanishing $\sP$-component. We deduce that 
\begin{equation}				\label{LONGEqXalphazeroaduPH}
	\begin{aligned}[]
		\pr_{\sP\sH}X^k_{\alpha 0}\neq 0\\
		\pr_{\sP\sQ}X^k_{0 \alpha }\neq 0.
	\end{aligned}
\end{equation}

As a consequence of compatibility between $\theta$ and $\sigma$, we have
\begin{equation}
	\begin{aligned}[]
		[J_1,(X_{\alpha\beta})_{\sH}]&=\alpha (X_{\alpha\beta})_{\sH}\\
		[J_1,(X_{\alpha\beta})_{\sQ}]&=\beta (X_{\alpha\beta})_{\sQ}
	\end{aligned}
\end{equation}
and
\begin{equation}		\label{LONGSubEqsJdeuxXalphaneta}
	\begin{aligned}[]
		[J_2,(X_{\alpha\beta})_{\sH}]&=\beta (X_{\alpha\beta})_{\sQ}\\
		[J_2,(X_{\alpha\beta})_{\sQ}]&=\alpha (X_{\alpha\beta})_{\sH}.
	\end{aligned}
\end{equation}
So $X_{\sQ}$ itself is an eigenvector of $\ad(J_1)$. In the same way, we prove that
\begin{equation}		\label{LONGEqJUnXabPK}
	\begin{aligned}[]
		[J_1,(X_{\alpha\beta})_{\sP}]=\alpha(X_{\alpha\beta})_{\sK}\\
		[J_1,(X_{\alpha\beta})_{\sK}]=\alpha(X_{\alpha\beta})_{\sP}
	\end{aligned}
\end{equation}
because $J_1\in\sP$.

\begin{corollary}		\label{LONGCorHPHKQPQKXuu}
	The vector $X_{++}$ has non vanishing components in $\sH\cap\sP$, $\sH\cap\sK$, $\sQ\cap\sP$ and $\sQ\cap\sK$.
\end{corollary}

\begin{proof}
	Since $\ad(J_2)$ inverts the $\sH$ and $\sQ$-components of $X_{++}$ (equation \eqref{LONGSubEqsJdeuxXalphaneta}), they must be both non zero. In the same way $\ad(J_1)$ inverts the components $\sP$ and $\sK$ of vectors of $\sH$ and $\sQ$ (equations \eqref{LONGEqJUnXabPK}).
\end{proof}

\begin{lemma}		\label{LONGLEmDesZPP}
	We have	$(X^{k}_{0+})_{\sK\sQ}=(X^{k}_{0+})_{\sP\sH}=0$ and consequently, $(X^k_{0+})_{\sP}=(X^k_{0+})_{\sQ}$.
\end{lemma}

\begin{proof}
	Consider the decomposition of the equality $[J_1,X^{k}_{0+}]=0$ into components $\sP\sQ$, $\sP\sH$, $\sK\sQ$, $\sK\sH$. Since $J_1\in\sP\cap\sH$, the $\sK\sH$ and $\sP\sQ$ components are
	\begin{subequations}
		\begin{align}
			[J_1,(X^{k}_{0+})_{\sP\sH}]&=0\\
			[J_1,(X^{k}_{0+})_{\sK\sQ}]&=0.
		\end{align}
	\end{subequations}
	In the same way, using the fact that $J_2\in\sP\cap\sQ$, we have
	\begin{subequations}	\label{LONGEqDeuxJUDzpKq}
		\begin{align}
			[J_2,(X_{0+}^{k})_{\sP\sH}]&=(X_{0+}^{k})_{\sK\sQ}\\
			[J_2,(X^{k}_{0+})_{\sK\sQ}]&=(X^{k}_{0+})_{\sP\sH}.		\label{LONGsubEqDeuxJUDzpKqb}
		\end{align}
	\end{subequations}
	Since $\dim(\sK\cap\sQ)=1$, the component $(X^{k}_{0+})_{\sK\sQ}$ has to be a multiple of $(X_{++})_{\sK\sQ}$. Thus we have
	\begin{equation}
		0=[J_1,(X^{k}_{0+})_{\sK\sQ}]=\lambda [J_1,(X_{++})_{\sK\sQ}]=\lambda (X_{++})_{\sP\sQ},
	\end{equation}
	but $(X_{++})_{\sP\sQ}\neq 0$, thus $\lambda=0$ and we conclude that $(X_{0+}^{k})_{\sK\sQ}=0$. Now, equation \eqref{LONGsubEqDeuxJUDzpKqb} shows that $(X_{0+}^{k})_{\sP\sH}=0$.
\end{proof}

\begin{lemma}				\label{LONGLemSigmaXzpBien}
	We have $\sigma X^k_{0+}=X^k_{0-}$.
\end{lemma}

\begin{proof}
	Since we know that $\sigma\sG_{\alpha\beta}\subset\sG_{\alpha,-\beta}$, the work is to fix the sign in
	\begin{equation}		\label{LONGEqsigmaXzppmthetaXzp}
		\sigma X^k_{0+}=\pm X^k_{0-}=\pm\theta X^k_{0+}.
	\end{equation}
	Lemma \ref{LONGLEmDesZPP} states that $(X_{0+}^k)_{\sP}=(X_{0+}^k)_{\sQ}$. Thus the $\sQ$-component of $\theta X^k_{0+}$ is $-(X^k_{0+})_{\sQ}$, which is also equal to the $\sQ$-component of $\sigma(X^k_{0+})$. That  fixes the choice of sign in equation \eqref{LONGEqsigmaXzppmthetaXzp}.
\end{proof}

The following is an immediate corollary of lemma \ref{LONGLemSigmaXzpBien} and the fact that $\theta$ fixes $\sP$ and $\sK$ while $\sigma$ fixes $\sH$ and $\sQ$.

\begin{corollary}		\label{LONGCorXzpHQPKXzm}
	We have
	\begin{subequations}
		\begin{align}
			(X^k_{0+})_{\sH}	&=(X^k_{0-})_{\sH}\\
			(X^k_{0+})_{\sQ}	&=-(X^k_{0-})_{\sQ}\\
			(X^k_{0+})_{\sP}	&=-(X^k_{0-})_{\sP}\\
			(X^k_{0+})_{\sK}	&=(X^k_{0-})_{\sK}.
		\end{align}
	\end{subequations}
\end{corollary}

\begin{proof}
	Since $\sigma$ acts as the identity on $\sH$ and changes the sign on $\sQ$, we have
	\begin{equation}
		\sigma X^k_{0+}=\sigma\big(  (X^k_{0+})_{\sH}+(X^k_{0+})_{\sQ}  \big)=(X^k_{0+})_{\sH}-(X^k_{0+})_{\sQ},
	\end{equation}
	but lemma \ref{LONGLemSigmaXzpBien} states that $\sigma X^k_{0+}=X^k_{0-}=(X^k_{0-})_{\sH}+(X^k_{0-})_{\sQ}$. Equating the $\sH$ and $\sQ$-component of these two expressions of $\sigma X^k_{0+}$ brings the two first equalities.

	The two other are proven the same way. We know that $\theta X^k_{0+}=X^k_{0-}$, but
	\begin{equation}
		\theta X^k_{0+}=\theta\big( (X^k_{0+})_{\sP}+(X^k_{0+})_{\sK} \big)=-(X^k_{0+})_{\sP}+(X^k_{0+})_{\sK}.
	\end{equation}
	The two last relations follow.
\end{proof}

\subsubsection{An interesting basis of \texorpdfstring{$\sQ$}{Q}}
\label{LONGSubSubSecInterestingBasisQ}

Being the tangent space of $AdS$, the space $\sQ$ is of a particular importance. Let us now have a closer look at the vectors that we already mentioned in equations \eqref{LONGAlignPRedDefQQ}:
\begin{subequations}				\label{LONGEqBasQQzi}
	\begin{align}
		q_0&=(X_{++})_{\sK\sQ}\\
		q_1&=J_2\\
		q_2&=-[J_1,q_0]=-(X_{++})_{\sP\sQ}\\
		q_k&=(X^k_{0+})_{\sQ}		&\text{lemma \ref{LONGLEmDesZPP}}.
	\end{align}
\end{subequations}

By lemma \ref{LONGLEmDesZPP}, and the discussion about $\mZ(\sK)$ (equation \eqref{LONGEqDecompZKenDeuxSuivantDim}), we can express the elements $q_i$ without explicit references to $\sQ$ itself and each element corresponds to a particular space (once again, the choice of $q_0$ is not that simple in $AdS_3$):

\begin{subequations}		\label{LONGEqBasQQziPlusMieux}%
	\begin{align}
		q_0&=(X_{++})_{\mZ(\sK)}&\in \sK\cap\sQ\cap\tilde\sN_2					\\
		q_1&=J_2&\in\sQ\cap\sA\\							
		q_2&=-[J_1,q_0]	&\in\sP\cap\sQ\cap\tilde\sN_2	\label{LONGsubEqJUnQZQDeux}		\\
		q_k&=(X^k_{0+})_{\sP}&\in\sP\cap\sQ\cap\tilde\sN_k.
	\end{align}
\end{subequations}
These elements correspond to the expression \eqref{LONGEqDecQEspacesCools}.
The compact part isomorphic to $\so(2,l-1)$ is then generated by the elements%

\begin{remark}

The space $\mZ(\sK)$ is given by the structure of the compact part of $\so(2,n)$, the elements $(X_{0+}^k)_{\sP}$ are defined from the root space structure of $\so(2,n)$ and the Cartan involution. The elements $J_1$ and $J_2$ are a basis of $\sA$. However, we need to know $\sH$ in order to distinguish $J_1$ from $J_2$ that are respectively generators of $\sA_{\sH}$ and $\sA_{\sQ}$.

Thus the basis 
\eqref{LONGEqBasQQziPlusMieux} %
is given in a way almost independent of the choice of $\sH$.
\end{remark}
\begin{corollary}		\label{LONGCorQdansPetK}
	We have $q_0\in\sK$ and $q_i\in\sP$ if $i\neq 0$ and the set $\{ q_0,q_1,\ldots,q_l \}$ is a basis of $\sQ$. Moreover, we have $\sQ\cap\tilde\sN_k=\langle q_k\rangle$.
\end{corollary}
\begin{proof}
	The first claim is a direct consequence of the expressions \eqref{LONGEqBasQQziPlusMieux}. Linear independence is a direct consequence of the spaces to which each vector belongs. A dimensional counting shows that it is a basis of $\sQ$.
\end{proof}

\subsubsection{Magic intertwining elements}

The vectors $q_i$ have the property to be intertwined by some elements of $\sH$. Namely, the adjoint action of the elements 
\begin{subequations}\label{LONGAlignDefMagicIntert}
	\begin{align}
		X_1=p_1&=-[J_2,q_0]	&\in\sP\cap\sH\cap\tilde\sN_2\\
		X_2=s_1&=[J_1,X_1]	&\in\sK\cap\sH\cap\tilde\sN_2\\
		X_k=-r_k&=-[J_2,q_k]	&\in\sK\cap\sH\cap\tilde\sN_k				\label{LONGEqDefXkCommeComm}
	\end{align}
\end{subequations}%
intertwines the $q_i$'s in the sense of the following proposition.
\begin{proposition}[Intertwining properties]			\label{LONGXUnALaTwistingSuperCool}
	The elements defined by equation \eqref{LONGAlignDefMagicIntert} satisfy
	\begin{multicols}{2}
		\begin{subequations}				\label{LONGEqCalculBBBJUnUnNirme}
			\begin{align}
				\ad(J_1)q_0&=-q_2		\label{LONGEqCalculBBBJUnUnNirmeA}\\
				\ad(J_1)q_2&=-q_0.		\label{LONGEqCalculBBBJUnUnNirmeB}
			\end{align}
		\end{subequations}
		\begin{subequations}				\label{LONGEqSubEqbXUnqZero}
			\begin{align}
				\ad(X_1)q_1&=q_0		\label{LONGSubEqbXZeroqUn}\\
				\ad(X_1)q_0&=q_1,		\label{LONGSubEqbXUnqZero}
			\end{align}
		\end{subequations}
		\begin{subequations}
			\begin{align}
				\ad(X_2)q_2&=q_1		\label{LONGSubEqXdeuxQdeuxa}\\
				\ad(X_2)q_1&=-q_2		\label{LONGSubEqXdeuxQun}
			\end{align}
		\end{subequations}
		\begin{subequations}				\label{LONGEqSubEqbXkqZero}
			\begin{align}
				\ad(X_k)q_k&=-q_1.		\label{LONGSubEqbXIMoinsqZero}\\
				\ad(X_k)q_1&=q_k		\label{LONGSubEqbXkQunQk}
			\end{align}
		\end{subequations}
	\end{multicols}
\end{proposition}

\begin{proof}

	Equation \eqref{LONGEqCalculBBBJUnUnNirmeA} is by definition while equation \eqref{LONGEqCalculBBBJUnUnNirmeB} follows from the first one and the fact that $\ad(J_1)^2$ acts as the identity on $\tilde\sN_2$. 

	The equality \eqref{LONGSubEqbXZeroqUn} is a direct consequence of the fact that $\ad(J_2)^2$ is the identity on $\tilde\sN_2$, so that
	\begin{equation}
		[X_1,q_1]=-\big[ [J_2,q_0],q_1 \big]=\ad(J_2)^2q_0=q_0.
	\end{equation}
	
	For the relation \eqref{LONGSubEqbXUnqZero}, first remark that, since $q_0=(X_{++})_{\sK\sQ}$, we have
	\begin{equation}
		X_1=-(X_{++})_{\sP\sH}
	\end{equation}
	and we have to compute
	\begin{equation}
		[X_1,q_0]=-\big[ (X_{++})_{\sP\sH},(X_{++})_{\sK\sQ} \big]
	\end{equation}
	Using the projections \eqref{LONGEqProjHQPKsigmatheta}, we have
	\begin{equation}
		\begin{aligned}[]
			(X_{++})_{\sP\sH}&=\frac{1}{ 4 }(X_{++}+\sigma X_{++}-\theta X_{++}-\sigma\theta X_{++})\\
			(X_{++})_{\sK\sQ}&=\frac{1}{ 4 }(X_{++}-\sigma X_{++}+\theta X_{++}-\sigma\theta X_{++})
		\end{aligned}
	\end{equation}
	We compute the commutator taking into account the facts that $\sigma$ is an automorphism and that, for example, $[X_{++},\sigma X_{++}]=0$ because $\sigma X_{++}\in\sG_{(+-)}$. What we find is
	\begin{equation}
		\big[ (X_{++})_{\sP\sH},(X_{++})_{\sK\sP} \big]=\frac{1}{ 4 }\frac{ 1 }{2}\Big( [X_{++},\theta X_{++}]-\sigma [X_{++},\theta X_{++}] \Big)=\frac{1}{ 4 }[X_{++},\theta X_{++}]_{\sQ}.
	\end{equation}
	Since $[X_{++},X_{--}]=-4(J_1+J_2)$, we have $[X_1,q_0]=J_2=q_1$ as expected.

	For equation \eqref{LONGSubEqXdeuxQdeuxa} we use the Jacobi relation and the relation \eqref{LONGEqCalculBBBJUnUnNirmeB}.
	\begin{equation}
		\begin{aligned}[]
			[q_2,X_2]&=\big[ q_2,[J_1,p_1] \big]\\
				&=-\big[ J_1,[p_1,q_2] \big]-\big[ p_1,[q_2,J_1] \big]\\
				&=-[p_1,q_0]\\
				&=-q_1
		\end{aligned}
	\end{equation}
	For equation \eqref{LONGSubEqXdeuxQun}, we use the definition of $X_2$, the Jacobi identity and the facts that $[p_1,J_2]=q_0$ and $[J_1,q_0]=-q_2$.

	We pass now to the fourth pair of intertwining relations. By definition, $q_k=(X^k_{0+})_{\sP}$, but taking into account the fact that $J_2\in\sP$ we can decompose the relation $[J_2,X_{0+}]=X_{0+}$ into
	\begin{subequations}
		\begin{align}
			[J_2,(X_{0+})_{\sP}]&=(X_{0+})_{\sK}		\label{LONGEqJdeuxXzpsPsK}\\
			[J_2,(X_{0+})_{\sK}]&=(X_{0+})_{\sP}.
		\end{align}
	\end{subequations}
	Equation \eqref{LONGEqJdeuxXzpsPsK} told us that 
	\begin{equation}
		X_k=-(X_{0+})_{\sK}.
	\end{equation}
	Now we have to compute $[X_k,q_k]=-\big[ (X_{0+})_{\sK},(X_{0+})_{\sP} \big]$. We know that $[X_{0+},X_{0-}]=-2J_2\in\sP$. Thus corollary \ref{LONGCorXzpHQPKXzm} brings
	\begin{equation}
		-2J_2=\big[ (X_{0+})_{\sK},(X_{0-})_{\sP} \big]+\big[ (X_{0+})_{\sP},(X_{0-})_{\sK} \big]=-2\big[ (X_{0+})_{\sK},(X_{0+})_{\sP} \big]=2[X_k,q_k],
	\end{equation}
	and the result follows.

	For equation \eqref{LONGSubEqbXkQunQk}, we have to compute $[X_k,q_1]=[J_2,(X_{+0}^k)_{\sK}]$. The $\sP$-component of $[J_2,X_{0+}^k]=X_{0+}^k$ is exactly
	\begin{equation}
		[J_2,(X_{0+}^k)_{\sK}]=(X_{0+}^k)_{\sP}=q_k.
	\end{equation}
	
\end{proof}

These intertwining relations will be widely used in computing the norm of the vectors $q_i$ in proposition \ref{LONGPropBaseQOrtho} as well as in some other occasions.

Let us now give a few words about the existence and unicity of these elements. The fact that there exists an element $X_1$ such that $\ad(J_2)X_1=q_0$ comes from the decomposition \eqref{LONGEqsDecopmQXpmpm} and the fact that each $X_{\pm\pm}$ is an eigenvector of $\ad(J_2)$. It is thus sufficient to adapt the signs in order to manage a combination of $X_{++}$, $X_{+-}$, $X_{-+}$ and $X_{--}$ on which the adjoint action of $J_2$ creates $q_0$. However, the fact that this element has in the same time the ``symmetric'' property $\ad(X_1)q_0=q_1$ could seem a miracle. See theorem \ref{LONGThoAdSqIouZero}.

\begin{lemma}
	An element $X_1$ such that $\ad(X_1)q_1=q_0$ can be chosen in $\sP\cap\sH\cap\tilde\sN_2$. Moreover, this choice is unique up to normalisation.
\end{lemma}

\begin{proof}
	Unicity is nothing else than the fact that $\dim(\sP\cap\sH\cap\tilde\sN_2)=1$. Indeed, since $\sG=\sA\oplus\tilde\sN\oplus\mZ_{\sK}(\sA)$ and $\sA\subset\sP$, we have $\sK\subset\tilde\sN$. Dimension counting shows that $\dim(\tilde\sN_2\cap\sH)=2$ (because $\dim(\tilde\sN_2)=4$ and $q_0,q_2\in\sQ\cap\tilde\sN_2$). As we are looking in $\tilde\sN_2$, we are limited to elements in $\so(2,2)$ (not the higher dimensional slices), so that we can consider $\sK=\so(2)\oplus\so(2)$. One of these two $\so(2)$ factors belongs to $\sH$, so that $\dim(\sK\cap\sH\cap\tilde\sN_2)=1$  and finally $\dim(\sP\cap\sH\cap\tilde\sN_2)=1$.

	Let now $X_1$ be such that $[X_1,q_1]=q_0$. If $X_1$ has a component in $\sQ$, that component has to commute with $q_1$ (if not, the commutator $[X_1,q_1]$ would have a $\sH$-component). So we can redefine $X_1$ in order to have $X_1\in\sH$.

	In the same way, a $\sA$-component has to be $J_1$ (because $J_2\in\sQ$) which commutes with $q_1$. We redefine $X_1$ in order to remove its $J_1$-component. We remove a component in $\tilde\sN_k$ because $[\tilde\sN_2,\tilde\sN_k]\subset\tilde\sN_k$, and a $\sK$-component can also be removed since its commutator with $q_1$ would produce a $\sP$-component. We showed that $X_1\in\sP\cap\sH\cap\tilde\sN_2$.
\end{proof}

\begin{lemma}		\label{LONGLemChoixDeXk}
	An element $X_k$ such that $\ad(X_k)q_1=q_k$ can be chosen in $\sK\cap\sH\cap\tilde\sN_k$.
\end{lemma}

\begin{proof}
	The proof is elementary in tree steps using the fact that $q_1\in\sP\cap\sQ\cap\sA$:
	\begin{enumerate}
		\item
			A $\sP$-component can be annihilated because $[\sP,\sP]\subset\sK$ while $q_k\in\sP$,
		\item
			a $\sQ$-component can be annihilated because $[\sQ,\sQ]\subset\sH$ while $q_k\in\sQ$,
		\item
			if $k'\neq k$, a $\tilde\sN_{k'}$-component can be annihilated because $[\tilde\sN_{k'},\sA]\subset\tilde\sN_{k'}$ while $q_k\in\tilde\sN_k$.
	\end{enumerate}
\end{proof}

\subsubsection{Killing form and orthogonality}

We \emph{define} the norm of an element in $\sG=\so(2,n)$ as
\begin{equation}	\label{LONGEqDefNormeKillingSix}
	\| X \|=-\frac{ 1 }{2n}B(X,X).
\end{equation}
Notice that $q_0$ belongs to the compact part of $\sG$, so that its Killing norm is negative, so that $\| q_0 \|$ is positive.

\begin{proposition}		\label{LONGPropBaseQOrtho}
	We have $\| q_0 \|=1$ and $\| q_i \|=-1$ ($i\neq 0$). As a consequence, the space $G/H$ is Lorentzian. 
\end{proposition}

\begin{proof}
	We begin by computing the norm of $q_1=J_2$. The Killing form $B(J_2,J_2)=\tr\big( \ad(J_2)\circ\ad(J_2) \big)$ is the easiest to compute in the basis $\mZ_{\sK}(\sA)\oplus\sA\oplus\sN\oplus\bar\sN$ of eigenvectors of $J_2$. If we look at the matrix of $\ad(J_2)\circ\ad(J_2)$, we have one $1$ at each of the positions of $X_{++}$, $X_{+-}$, $X_{-+}$ and $X_{--}$. Moreover, for each higher dimensional slice, we get additional $2$ because of $X_{0+}^k$ and $X_{+0}^k$. When one look at $\so(2,n)$ we have $n-2$ higher dimensional slices, so that 
	\begin{equation}
		B(J_2,J_2)=4+2(n-2)=2n.
	\end{equation}
	
	The result is that $B(q_1,q_1)=6$, so that $\| q_1 \|=-1$.

	We are going to propagate that result to other elements of the basis, using the``magic'' intertwining elements $X_1$, $X_k$ and $J_1$. 

	Using left invariance of the Killing form, we find
	\begin{equation}		\label{LONGEqCalculBBBqZeroqUnNirme}
		B(q_0,q_0)=B\big( q_0,-\ad(J_1)q_2 \big)=B\big( \ad(J_1)q_0,q_2 \big)=-B(q_2,q_2),
	\end{equation}
	so that $\| q_0 \| =-\| q_2 \|$.
	
	Now, the same computation with $X_1$ and $X_k$ instead of $J_1$ show that $\| q_0 \|=-\| q_1 \|$ and $\| q_1 \|=\| q_k \|$. 
\end{proof}

\begin{remark}		\label{LONGRemBProdScal}
	Using the fact that the basis $\{ q_i\}$ is orthonormal, we can decompose an element of $\sQ$ by the Killing form. One only has to be careful on the sign: if $X=aq_0+\sum_{i>0}b_iq_i$, we have
	\begin{equation}		\label{LONGEqsabKillProjComp}
		\begin{aligned}[]
			a&=B(X,q_0)\\
			b_i&=-B(X,q_i).
		\end{aligned}
	\end{equation}
\end{remark}

\begin{remark}	\label{LONGRemOrdreNilpotentQ}
	As a consequence, a light like direction reads, up to normalization, $E(w)=q_0+\sum_{i=1}^{l-1}w_iq_i$ with $w\in S^{l-2}$.
\end{remark}

\subsubsection{Other properties}

\begin{lemma}		\label{LONGLemXZUAHetQ}
	We have $\sigma X_{\alpha\beta}\in\sG_{(\alpha,-\beta)}$. In particular, $X^k_{0+}$ has non vanishing components in $\sH$ and in $\sQ$.
\end{lemma}

\begin{proof}
	If one applies $\sigma$ to the equality $[J_2,X_{\alpha\beta}]=\beta X_{\alpha\beta}$, we see that $\sigma X_{\alpha\beta}$ is an eigenvector of $\ad(J_2)$ with eigenvalue $-\beta$. The same with $\ad(J_1)$ shows that $\sigma X_{\alpha\beta}$ has $+1$ as eigenvalue. Thus $\sigma X_{\alpha\beta}\in\sG_{(\alpha,-\beta)}$. 
	
	In particular, $\sigma X^k_{0+}\neq \pm X^k_{0+}$ so that it does not belongs to $\sH$ nor to $\sQ$.
\end{proof}

Notice that, as corollary, we have 
\begin{equation}
	\sigma X_{\alpha,\beta}=\pm X_{\alpha,-\beta}.
\end{equation}

\begin{lemma}				\label{LONGLemSigmaXppEgalXPm}
	We have $(X_{++})_{\sQ}=(X_{+-})_{\sQ}$ or, equivalently, $\sigma X_{++}=-X_{+-}$.
\end{lemma}

\begin{proof}
	Since $q_1=J_2\in\sA$ and $q_k\in\tilde\sN_k$, the $\sQ$-component of $X_{++}$ and $X_{+-}$ are only made of $q_0$ and $q_2$.
	We are	going to prove the following three equalities.
	\begin{enumerate}

		\item\label{LONGItemBpmqDeux}
			$B(X_{+-},q_2)=B(X_{+-},q_0)$
		\item
			$B(X_{++},q_2)=B(X_{++},q_0)$
		\item
			$B(X_{++},q_0)=B(X_{+-},q_0)$
	\end{enumerate}
	
	The first point is proved using the fact that $q_2=[q_0,J_1]$ and the $\ad$-invariance of the Killing form:
	\begin{equation}
		B(X_{+-},q_2)=-B\big( X_{+-},\ad(J_1)q_0 \big)=B\big( \ad(J_1)X_{+-},q_0 \big)=B(X_{+-},q_0).
	\end{equation}
	One checks the second point in the same way. For the third equality, we know from decomposition \eqref{LONGEqDecompZKenDeuxSuivantDim} that $q_0$ is a multiple of $X_{++}+X_{--}+X_{+-}+X_{-+}$. If the multiple is $\lambda$, $B(X_{++},q_0)=\lambda B(X_{++},X_{--})$ and $B(X_{+-},q_0)=\lambda B(X_{+-},X_{-+})$. Thus we have to prove that the traces of the operators
	\begin{equation}
		\begin{aligned}[]
			\gamma_1&=\ad(X_{++})\circ\ad(X_{--})\\
			\gamma_2&=\ad(X_{+-})\circ\ad(X_{-+})
		\end{aligned}
	\end{equation}
	are the same. That trace is straightforward to compute on the natural basis of $\sG=\mZ_{\sK}(\sA)\oplus\sA\oplus\sN\oplus\bar\sN$. 
	The only elements on which $\ad(X_{--})$ is not zero are $\sA$, $X^k_{0+}$, $X^k_{+0}$ and $X^k_{++}$, while for $\ad(X_{-+})$, the only non vanishing elements are $\sA$, $X^k_{0-}$, $X^k_{+0}$ and $X_{+-}$. From equation \eqref{LONGsubeqrXpmpm}, we have $\gamma_1(r_{ij})=\gamma_2(r_{ij})=0$. Using the commutation relations, we find
	\begin{subequations}
		\begin{align}
			\gamma_1J_1&=[X_{++},X_{--}]=-4(J_1+J_2) \\
			\gamma_1J_2&=[X_{++},X_{--}]=-4(J_1+J_2) \\
			\gamma_1X^k_{0+}&=2[X_{++},X^k_{-0}]=-4X^k_{0+} \\
			\gamma_1X^k_{+0}&=-2[X_{++},X^k_{0-}]=-4X^k_{+0} \\
			\gamma_1X_{++}&=[X_{++},4(J_1+J_2)]=-8X_{++}.
		\end{align}
	\end{subequations}
	Thus $\tr(\gamma_1)=-24$. The same computations bring
	\begin{subequations}
		\begin{align}
			\gamma_2J_1&=[X_{+-},X_{-+}]=-4(J_1-J_2)\\
			\gamma_2J_2&=[X_{+-},X_{-+}]=4(J_1-J_2)\\
			\gamma_2X^k_{0-}&=-2[X_{+-},X^k_{-0}]=-4X^k_{0-}\\
			\gamma_2X_{+-}&=[X_{+-},4(J_1-J_2)]=-8X_{+-}\\
			\gamma_2X^k_{+0}&=2[X_{+-},X^k_{0+}]=-2X^k_{+0},
		\end{align}
	\end{subequations}
	and $\tr(\gamma_2)=-24$. Thus we have 
	\begin{equation}
		\pr_{\mZ(\sK)}(X_{++})=\pr_{\mZ(\sK)}(X_{+-}).
	\end{equation}
\end{proof}

Notice that the lemma is trivial if we consider that $X_{++}-X_{+-}$ belongs to $\sH$ by definition of $\sH$. From a $AdS$ point of view, we define $AdS=G/H$ and we have to define $H$, so from that point of view, lemma \ref{LONGLemSigmaXppEgalXPm} is by definition. However, the direction we have in mind is to use the more generic tools as possible. From that point of view, the fact to set $\mZ(\sK)\subset\sQ$ is more intrinsic than to set $X_{++}-X_{+-}\in\sH$.

\begin{proposition}		\label{LONGPropXmpXppqq}
	We have $(X_{++})_{\sQ}=(X_{+-})_{\sQ}=q_0-q_2$.
\end{proposition}

\begin{proof}
	Using the remark \ref{LONGRemBProdScal}, the three Killing forms computed in the proof of lemma \ref{LONGLemSigmaXppEgalXPm} are expressed under the form
	\begin{subequations}
		\begin{align}
			(X_{+-})_{q_0}&=-(X_{+-})_{q_2}\\
			(X_{++})_{q_0}&=-(X_{++})_{q_2}\\
			(X_{++})_{q_0}&=(X_{+-})_{q_2}.
		\end{align}
	\end{subequations}
	Consequently, we have $(X_{++})_{\sQ}=\lambda(q_0-q_2)$ and $(X_{+-})_{\sQ}=\lambda(q_0-q_2)$ for a constant $\lambda$ to be fixed. It is fixed to be $1$ by the facts that, by definition, $q_0=(X_{++})_{\sK\sQ}$ and $q_2\in\sP$.
\end{proof}

\begin{lemma}		\label{LONGLemComJDeuxQ}
	We have
	\begin{equation}
		\begin{aligned}[]
			[J_2,q_0]&=(X_{++})_{\sH\sP}\neq 0\\
			[J_2,q_1]&=0\\
			[J_2,q_2]&=(X_{++})_{\sH\sK}\neq 0\\
			[J_2,q_k]&=(X^k_{0+})_{\sH}\neq 0\\
			[J_1,p_1]&=-(X_{++})_{\sK\sH}\neq 0
		\end{aligned}
	\end{equation}
	where $k\geq 3$. 
\end{lemma}

\begin{proof}
	Using the fact that $J_2\in\sQ\cap\sP$ and that $X_{++}$ has non vanishing components ``everywhere'' (corollary \ref{LONGCorHPHKQPQKXuu}), we have
	\begin{equation}
		\begin{aligned}[]
			[J_2,q_0]&=[J_2,(X_{++})_{\sK\sQ}]=(X_{++})_{\sP\sH}\neq 0\\
			[J_2,q_2]&=[J_2,(X_{++})_{\sP\sQ}]=(X_{++})_{\sK\sH}\neq 0\\
			[J_1,p_1]&=[J_1,(X_{++})_{\sP\sH}]=-(X_{++})_{\sK\sH}\neq 0\\
			[J_2,q_k]&=[J_2,(X_{0+}^k)_{\sQ}]=(X^k_{0+})_{\sH}\neq 0		&\text{lemma \ref{LONGLemXZUAHetQ}}\\
		\end{aligned}
	\end{equation}
\end{proof}

\begin{lemma}		\label{LONGLemNonHXaz}
	If $\alpha\neq 0$, then $X_{\alpha 0}^k\in\sH$.
\end{lemma}

\begin{proof}
	The element $\pr_{\sQ}X^k_{\alpha 0}$ is a combination of $q_i$. Since $\ad(J_2)\pr_{\sQ}X^k_{\alpha 0}=0$, we must have $(X^k_{\alpha 0})_{\sQ}=\lambda J_2$ by lemma \ref{LONGLemComJDeuxQ}. Using the fact that $J_1\in\sH$, the $\sQ$-component of the equality $[J_1,X^k_{\alpha 0}]=\alpha X^k_{\alpha 0}$ becomes
	\begin{equation}
		[J_1,\lambda J_2]=\alpha\lambda J_2.
	\end{equation}
	The left-hand side is obviously zero, so that $\lambda=0$ which proves that $X^k_{\alpha 0}\in\sH$.
\end{proof}

Applying successively the projections \eqref{LONGEqProjHQPKsigmatheta}, and lemma \ref{LONGLemSigmaXppEgalXPm}, we write the basis elements of $\sQ$ in the decomposition $\sG=\sG_0\oplus\sN\oplus\bar\sN$ :
\begin{subequations}			\label{LONGEqsDecopmQXpmpm}
	\begin{align}
		q_0&=\frac{1}{ 4 }(X_{++}+X_{+-}+X_{-+}+X_{--}),		\\
		q_1&=J_2,							\\
		q_2&=\frac{1}{ 4 }(-X_{++}-X_{+-}+X_{-+}+X_{--}),		\\
		q_k&=\frac{ 1 }{2}(X^{k}_{0+}-X^{k}_{0-})
	\end{align}
\end{subequations}
with $k\geq 3$. Notice that none of them has component in $\mZ_{\sK(\sA)}$.

These decompositions allow us to compute the commutators $[q_i,q_j]$ and $[q_i,J_p]$. Instead of listing here every commutation relations, we will only write the ones we use when we need them.

\begin{lemma}		\label{LONGLemQzQdeuxJun}
	We have $[q_0,q_2]=-J_1$.
\end{lemma}

\begin{proof}
	The proof is exactly the same as the one of equation \eqref{LONGSubEqbXUnqZero} in lemma \ref{LONGXUnALaTwistingSuperCool}. Here we use
	\begin{equation}
		(X_{++})_{\sP\sQ}=\frac{1}{ 4 }\big( X_{++}-\sigma X_{++}-\theta X_{++}+\sigma\theta X_{++} \big)
	\end{equation}
	and we find
	\begin{equation}
		[q_0,q_2]=-\big[ (X_{++})_{\sK \sQ},(X_{++})_{\sP\sQ} \big]=-\frac{1}{ 4 }[\theta X_{++},X_{++}]_{\sH}=-J_1.
	\end{equation}
\end{proof}

\begin{lemma}
We have
	\begin{equation}		\label{LONGEqXunQdeuxcommutent}
		[X_1,q_2]=[X_1,q_k]=0
	\end{equation}
	for $k\geq 3$.
\end{lemma}
\begin{proof}
	The proof is elementary:
	\begin{equation}
		\begin{aligned}[]
			[X_1,q_2]&\in[\sP\cap\sH\cap\tilde\sN_2,\sP\cap\sQ\cap\tilde\sN_2]\subset\sK\cap\sQ\cap\sA=\{ 0 \}\\
			[X_1,q_k]&\in[\sP\cap\sH\cap\tilde\sN_2,\sP\cap\sQ\cap\tilde\sN_k]\subset\sK\cap\sQ\cap\tilde\sN_k=\{ 0 \}.
		\end{aligned}
	\end{equation}
\end{proof}

The following is a first step in the proof of theorem \ref{LONGThoAdESqqq}.
\begin{corollary}			\label{LONGCorAdQUncarreqi}
	We have $\ad(J_1)|_{\tilde\sN_2}^2=\ad(J_2)|_{\tilde\sN_2}^2=\id$ and $\ad(q_1)^2q_i=q_i$.
\end{corollary}

\begin{proof}
	The action of $\ad(q_1)^2$ is to change two times the sign of the components $X_{\alpha -}$. Thus $\ad(q_1)^2=\id$ on $\tilde\sN_2$. The result is now proved for $i=0,1,2$. For the higher dimensions, we use the fact that $J_2=q_1$ and we find
	\begin{equation}
		q_k=[X_k,q_1]=-\big[ [q_1,X_k],q_1 \big]=\ad(q_1)^2q_k
	\end{equation}
	as claimed.

	Moreover, the elements of $\tilde\sN_2$ are build of elements of the form $X_{\pm\pm}$, so that $\ad(J_1)^2$ changes at most twice the sign.

\end{proof}

\begin{lemma}		\label{LONGLemXkqzerozero}\label{LONGLemXkJunzero}
	We have
	\begin{subequations}
		\begin{align}
			[X_k,q_0]=[X_k,J_1]=[X_k,q_2]=0\\
			[J_1,q_k]=0.		\label{LONGEqJUnqkzero}
		\end{align}
	\end{subequations}
\end{lemma}

\begin{proof}
	The first claim is proved in a very standard way:
	\begin{equation}
		[X_k,q_0]\in[\sK\cap\sH\cap\tilde\sN_k,\sK\cap\sQ\cap\tilde\sN_2]\subset\sK\cap\sQ\cap\tilde\sN_k=\{ 0 \}.
	\end{equation}

	For the second commutator, we use the Jacobi identity and the definition $X_k=-[J_2,q_k]$:
	\begin{equation}		\label{LONGEqJ1J2qkzero}
		\big[ J_1,[J_2,q_k] \big]=-\big[ J_2,[q_k,J_1] \big]-\big[ q_k,\underbrace{[J_1,J_2]}_{=0} \big],
	\end{equation}
	while
	\begin{equation}
		[q_k,J_1]\in[\sP\cap\sQ\cap\tilde\sN_k,\sP\cap\sH\cap\sA]\subset\sK\cap\sQ\cap\tilde\sN_k=\{ 0 \}.
	\end{equation}
	That proves \eqref{LONGEqJUnqkzero} in the same time.

	For the third commutator, remark that, since $q_2=[q_0,J_1]$, we have
	\begin{equation}
		[X_k,q_2]=-\big[ q_0,[J_1,X_k] \big]-\big[ J_1,[X_k,q_0] \big].
	\end{equation}
	which is zero by the two first claims.
\end{proof}

\begin{proposition}			\label{LONGEtOrdreDeux}
	If $E$ is nilpotent in $\sQ$, then $\ad(E)^3=0$.
\end{proposition}

\begin{proof}
	If $E_1$ is a nilpotent element of $\sQ$, then every nilpotent elements in $\sQ$ are of the form $\lambda\Ad(k)E_1$ for some $k\in K$ and $\lambda\in\eR$\cite{These}. It is then sufficient to prove that one of them is of order two. The element
	\begin{equation}		\label{LONGEqDecqzmoinsqDeux}
		q_0-q_2=\frac{ 1 }{2}(X_{++}+X_{+-}),
	\end{equation}
	is obviously of order two because the eigenvalue for $\ad(J_1)$ increases by one unit at each iteration of $\ad(q_0-q_2)$.
\end{proof}

\begin{lemma}
	We have $[J_1,q_k]=0$.
\end{lemma}

\begin{proof}
	The proof is standard :
	\begin{equation}
		[J_1,q_k]\in[\sH\cap\sP\cap\sA,\sQ\cap\sP\cap\tilde\sN_k]\subset\sQ\cap\sK\cap\tilde\sN_k=\{ 0 \}.
	\end{equation}
\end{proof}

The following theorem, which relies on the preceding lemmas, will be central in computing the Killing form which appears in the characterization of theorem \ref{LONGThosSequivJzero}.
\begin{theorem}			\label{LONGThoAdESqqq}
	We have
	\begin{equation}
		\ad(q_i)^2q_j=q_j
	\end{equation}
	if $i\neq j$ and $i\neq 0$. If $i=0$, we have
	\begin{equation}
		\ad(q_0)^2q_j=-q_j.
	\end{equation}
\end{theorem}

\begin{proof}
	The case $i=1$ is already done in corollary \ref{LONGCorAdQUncarreqi}. 

	We propagate that result to the other $\ad(q_i)^2$ with the intertwining elements $J_1$, $X_1$ and $X_k$.

	Let us compute $\ad(q_0)^2q_i=\ad\big([X_1,q_1]\big)^2q_i$ using twice the Jacobi identity and the properties of $X_1$ (in order to be more readable, we write $XY$ for $[X,Y]$)
	\begin{equation}		\label{LONGEqAdqZsqqi}
		\begin{aligned}[]
			\ad(q_0)^2q_i	&=	(X_1q_1)\Big( (X_1q_1)q_i \Big)\\
					&=	-(X_1q_1)\Big( (q_1q_I)X_1+(q_iX_1)q_1 \Big)\\
					&=(q_1q_i)\big( X_1(X_1q_1) \big)+(q_iX_1)\big( q_1(X_1q_1) \big)\\
					&\quad + X_1\big( (X_1q_1)(q_1q_i) \big)+q_1\big( (X_1q_1)(q_iX_1) \big)\\
					&=(q_1q_i)q_1-\ad(X_1)^2q_i+X_1\big( q_0(q_1q_i) \big)+q_1\big( q_0(q_iX_1) \big).
		\end{aligned}
	\end{equation}
	
	If $i=1$, the only non vanishing term is $-\ad(X_1)^2q_1=-q_1$.  Thus $\ad(q_0)^2q_1=-q_1$.
	
	If $i=2$, the relation \eqref{LONGEqXunQdeuxcommutent} annihilates the second and fourth terms while $[q_1,q_2]$ commutes with $q_0$ because $q_0\in\mZ(\sK)$. We are thus left with the term $-q_2$, and $\ad(q_0)^2q_2=-q_2$.

	If $i=k\geq 3$, we find
	\begin{equation}
		\ad(q_0)^2q_k=-\ad(q_1)^2q_k-\ad(X_1)^2q_k+X_1\big( q_0(q_1q_k) \big)+q_1\big( q_0(q_kX_1) \big).
	\end{equation}
	Since $[q_1,q_k]\in\sK$, it commutes with $q_0$. Using the fact that $[X_1,q_k]=0$, we get $\ad(q_0)^2q_k=-q_k$.

	Let us perform the same computations  as in \eqref{LONGEqAdqZsqqi} with $q_k$ ($k\geq 3$) instead of $q_0$ and $X_k$ (equations \eqref{LONGEqSubEqbXkqZero}) instead of $X_1$. What we get is
	\begin{equation}
		\ad(q_k)^2q_i=\ad(q_1)^2q_i-\ad(X_k)^2q_i+X_k\big( q_k(q_1q_i) \big)+q_1\big( q_k(q_iX_k) \big).
	\end{equation}

	If we set $i=0$, taking into account the commutator $[X_k,q_0]=0$, we have
	\begin{equation}		\label{LONGEqkdeuxqzerointer}
		\ad(q_k)^2q_0=\ad(q_1)^2q_0+X_k\big( q_k(q_1q_0) \big).
	\end{equation}
	As already proved, the first term is $q_0$. Now,
	\begin{equation}
		\big[ q_k,[q_1,q_0] \big]\in\sK\cap\sQ\cap\tilde\sN_k=\{ 0 \},
	\end{equation}
	so that the second term in \eqref{LONGEqkdeuxqzerointer} is zero. Thus we proved that $\ad(q_k)^2q_0=q_0$.

	If we set $i=1$, taking into account the relations \eqref{LONGEqSubEqbXkqZero}, we find
	\begin{equation}
		\ad(q_k)^2q_1=-\ad(X_k)^2q_1+q_1\big( q_k(q_1X_k) \big)=q_1.
	\end{equation}

	If we set $i=2$ and using the fact that $[X_k,q_2]=0$, we find
	\begin{equation}
		\ad(q_k)^2q_2=q_2-q_k\big( X_k(q_1q_2) \big).
	\end{equation}
	Using once again the Jacobi identity inside the big parenthesis, we find $2q_2-\ad(q_k)^2q_2$. This proves that $\ad(q_k)^2q_2=q_2$.

	We turn now our attention to $\ad(q_2)^2q_i$. We perform the same computation, using the intertwining property \eqref{LONGEqCalculBBBJUnUnNirme} of $J_1$. What we get is
	\begin{equation}		\label{LONGEqadqqqDeuxIGene}
		\ad(q_2)^2q_i=(J_1q_i)(q_0q_2)-\ad(q_0)^2q_i+q_0\big( q_2(J_1q_i) \big)+J_1\big( q_2(q_iq_0) \big).
	\end{equation}

	If we set $i=1$, we use the already proved property $\ad(q_0)^2q_1=-q_1$, and we obtain
	\begin{equation}
		\ad(q_2)^2q_1=(J_1q_1)(q_0q_2)+q_1+q_0\big( q_2(J_1q_1) \big)+J_1\big( q_2(q_1q_0) \big).
	\end{equation}
	We claim that all of these terms are zero except of $q_1$. First, $\big[ q_2,[q_1,k_k] \big]\in\big[ \tilde\sN_2,[\sA,\tilde\sN_2] \big]\subset\sA$.	Thus the last term vanishes .The commutator $[J_1,q_1]$ vanishes because $q_1=J_2$. We are done with $\ad(q_2)^2q_1=q_1$.

	If we set $i=k$ ($k\geq 3$) in \eqref{LONGEqadqqqDeuxIGene}, we use $\ad(q_0)^2q_k=-q_k$ and what we find is
	\begin{equation}		\label{LONGEqadqqqdeuxkk}
		\ad(q_2)^2q_k=(J_1q_k)(q_0q_2)+q_k+q_0\big( q_2(J_1q_k) \big)+J_1\big( q_2(q_kq_0) \big).
	\end{equation}
	We already know that $[J_1,q_k]=0$. We have $\big[ q_2,[q_k,q_0] \big]=0$ because
	\begin{equation}
		\begin{aligned}[]
			\big[ q_2,[q_k,q_0] \big]&\in\big[ \sP\cap\sQ\cap\tilde\sN_2,[\sP\cap\sQ\cap\tilde\sN_k,\sK\cap\sQ\cap\tilde\sN_2] \big]\\
						&\subset[\sP\cap\sQ\cap\tilde\sN_2,\sP\cap\sH\cap\tilde\sN_k]\\
						&\subset\sK\cap\sQ\cap\tilde\sN_k=\{ 0 \}.
		\end{aligned}
	\end{equation}
	The remaining terms in \eqref{LONGEqadqqqdeuxkk} are $\ad(q_2)^2q_k=q_k$.

	In order to compute $\ad(q_2)^2q_0$, we write $q_0=\ad(X_1)q_1$. Using twice the Jacobi identity, we get
	\begin{equation}
		\ad(q_2)^2q_0=X_1\big( (q_1q_2)q_2 \big)+q_1\big( (X_1q_2)q_2 \big)+(q_1q_2)(X_1q_2)+(X_1q_2)(q_2q_1).
	\end{equation}
	Using the fact that $[X_1,q_2]=0$, we are left with
	\begin{equation}
		\ad(q_2)^2q_0=X_1\big( \ad(q_2)^2q_1 \big)=[X_1,q_1]=q_0
	\end{equation}
	as desired.
\end{proof}

\subsection{A convenient basis for the root spaces and computations}
\label{LONGSubSecMOreConvBasisBlbla}

The most natural basis of $\tilde\sN_2$ is
\begin{equation}
	\tilde\sN_2=\langle X_{++},X_{+-},X_{-+},X_{--}\rangle,
\end{equation}
but the multiple commutators of these elements with $q_0$ reveal to require some work. 

We provide in this section an other basis for $\tilde\sN$ that corresponds to the decomposition $\sK\oplus\sP$. Since $q_0$ is central in $\sK$, the exponential $e^{xq_0}X$ is trivial when $X\in\sK$ and, since $q_0\in\sK$, the commutator $[q_0,X]$ remains in $\sP$ when $X\in\sP$.

Here is the new basis:
\begin{equation}		\label{LONGEqSuperBaseeB}
	\eB=\{J_1,J_2, q_0,q_2,p_1,s_1,q_k,p_k,r_k,s_k \}_{k=3,\ldots,l-1}
\end{equation}
where
\begin{subequations}		\label{LONGSubEqsBaseSuperTsNTrois}\label{LONGEqDecompqprskXX}
	\begin{align}
		q_0&		&&=\frac{1}{ 4 }(X_{++}+X_{+-}+X_{-+}+X_{--})	&\in\sK\cap\sQ\cap\tilde\sN_2\\
		q_2&=[q_0,J_1]	&&=\frac{1}{ 4 }(-X_{++}-X_{+-}+X_{-+}+X_{--})	&\in\sP\cap\sQ\cap\tilde\sN_2\\
		p_1&=[q_0,q_1]	&&=\frac{1}{ 4 }(-X_{++}+X_{+-}-X_{-+}+X_{--})	&\in\sP\cap\sH\cap\tilde\sN_2	\label{LONGEqqzqupuDef}\\
		s_1&=[J_1,p_1]	&&=\frac{1}{ 4 }(-X_{++}+X_{+-}+X_{-+}-X_{--})	&\in\sK\cap\sH\cap\tilde\sN_2\\
		q_k&		&&=\frac{ 1 }{2}(X_{0+}^k-X_{0-}^k)	&\in\sP\cap\sQ\cap\tilde\sN_k	\\
		p_k&=[q_0,q_k]	&&=\frac{ 1 }{ 2 }(X_{-0}^k-X_{+0}^k)	&\in\sP\cap\sH\cap\tilde\sN_k	\\
		r_k&=[J_2,q_k]	&&=\frac{ 1 }{2}(X^k_{0+}+X^k_{0-})	&\in\sK\cap\sH\cap\tilde\sN_k	\\
		s_k&=[J_1,p_k]	&&=-\frac{ 1 }{2}(X^k_{-0}+X^k_{+0})	&\in\sK\cap\sH\cap\tilde\sN_k	\\
		r_{ij}&=[q_i,q_j]&&					&\in\sK\cap\sH\cap\sG_0.
	\end{align}
and
	\begin{align}
		J_1&\in\sP\cap\sH\cap\sA\\
		q_1=J_2&\in\sP\cap\sQ\cap\sA
	\end{align}
\end{subequations}
Notice that the elements $p_1$ and $s_1$ are non vanishing by lemma \ref{LONGLemComJDeuxQ}. We have 
\begin{subequations}
	\begin{align}
		\sP\cap\sQ\cap\tilde\sN_2&=\langle q_2\rangle	&	\sK\cap\sQ\cap\tilde\sN_2&=\langle q_0\rangle\\
		\sP\cap\sQ\cap\tilde\sN_k&=\langle q_k\rangle	&	\sK\cap\sQ\cap\tilde\sN_k&=\emptyset\\
		\sP\cap\sQ\cap\sA&=\langle J_2\rangle	&	\sK\cap\sQ\cap\sA&= \emptyset\\
		\sP\cap\sH\cap\tilde\sN_2&=\langle p_1\rangle	&	\sK\cap\sH\cap\tilde\sN_2&=\langle s_1\rangle\\
		\sP\cap\sH\cap\tilde\sN_k&=\langle p_k\rangle	&	\sK\cap\sH\cap\tilde\sN_k&=\langle r_k,s_k\rangle\\
		\sP\cap\sH\cap\sA&=\langle J_1\rangle	&	\sK\cap\sH\cap\sA&=\emptyset
	\end{align}
\end{subequations}

The decomposition of $\tilde\sN_2$ into $\sK\oplus\sP$ is
\begin{equation}		\label{LONGEqDecomptsNTroisKP}
	\tilde\sN_2=\langle q_0,s_1\rangle\oplus \langle q_1,p_1\rangle.
\end{equation}

The decomposition of $\tilde\sN_k$ into $\sK\oplus\sP$ is
\begin{equation}		\label{LONGEqDecomptsNkKP}
	\tilde\sN_k =\langle r_k,s_k,r_{ij}\rangle_{k,i,j\geq 3}\oplus \langle q_k,p_k\rangle.
\end{equation}

We are now going to compute all the Killing form and commutators in this basis.
\begin{proposition}		\label{LONGPropBprsk}
	We have
	\begin{equation}
		\begin{aligned}[]
			B(p_k,p_k)&=-B(q_0,q_0)\\
			B(r_k,r_k)&=B(q_0,q_0)\\
			B(s_k,s_k)&=B(q_0,q_0),
		\end{aligned}
	\end{equation}
	and then
	\begin{equation}	\label{LONGeqNormInHigherDimensionalSlices}	
		-\| p_k \|^2=\| r_k \|^2=\| s_k \|^2=1.
	\end{equation}
\end{proposition}

\begin{proof}
	Using the definitions \eqref{LONGEqDecompqprskXX} and the theorem \ref{LONGThoAdESqqq}, we have
	\begin{subequations}		\label{LONGSubEqsBpprrssk}
	\begin{equation}	
			B(p_k,p_k)=B\big( \ad(q_k)q_0,\ad(q_k)q_0 \big)
					=-B\big( \ad(q_k)^2q_0,q_0 \big)
					=-B(q_0,q_0),
	\end{equation}
	and
	\begin{equation}
		B(r_k,r_k)=B\big( \ad(q_1)q_k,\ad(q_1)q_k \big)=-B(q_k,q_k)=B(q_0,q_0).
	\end{equation}
	and
	\begin{equation}
		B(s_k,s_k)=-B\big( \ad(J_1)^2p_k,p_k \big)=-B(p_k,p_k)=B(q_0,q_0).
	\end{equation}
	\end{subequations}
\end{proof}
Notice that we are not surprised by the positivity of the norms of $r_k$ and $s_k$ because they belong to the compact part of the algebra.

\begin{proposition}
	The Killing norm in the space $\mZ_{\sK}(\sA)$ are given by
	\begin{equation}
		B(r_{ij},r_{kl})=\begin{cases}
			B(q_0,q_0)	&	\text{if $\{ i,j \}=\{ k,l \}$}\\
			0	&	 \text{\text{otherwise}}.
		\end{cases}
	\end{equation}
\end{proposition}

\begin{proof}
	First, we have
	\begin{equation}
		\begin{aligned}[]
			B(r_{ij},r_{ij})&=B\big( [q_i,q_j],[q_i,q_j] \big)\\
					&=-B\big( \ad(q_i)^2q_j,q_j \big)\\
					&=-B(q_j,q_j)\\
					&=B(q_0,q_0).
		\end{aligned}
	\end{equation}
	For the mixed case we have
	\begin{equation}	
			B(r_{ij},r_{ik})=-B\big( \ad(q_i)^2q_j,q_k \big)=0.
	\end{equation}
	We suppose now that $i$, $j$, $k$ and $l$ are four different numbers. The action of $\ad(r_{kl})$ on $\sA$ is zero because $r_{kl}\in\mZ(\sA)$. From \eqref{LONGEqComsRRN}, the action of $\ad(r_{ij})\circ\ad(r_{kl})$ on $\sN$ is zero. Since the elements $r_{ij}$ satisfy the algebra of $\so(n)$, we have $\ad(r_{ij})\circ\ad(r_{kl})r_{mn}=0$ when $i$, $j$, $k$ and $l$ are four different numbers. Finally we have
	\begin{equation}
		B(r_{ij},r_{kl})=0.
	\end{equation}
\end{proof}

\begin{lemma}		\label{LONGLempunjdeuxqzero}
	We have $[p_1,J_2]=q_0$.
\end{lemma}

\begin{proof}
	The lemma comes from theorem \ref{LONGThoAdESqqq} because
	\begin{equation}
		[p_1,J_2]=-[q_1,p_1]=\ad(q_1)^2q_0=q_0.
	\end{equation}
\end{proof}

\begin{lemma}				 \label{LONGlemJDeuxqDeuxsUn}
	We have $s_1=[J_2,q_2]$.
\end{lemma}

\begin{proof}
	We use the definition $p_1=[q_0,q_1]$ and the Jacobi identity:
	\begin{equation}
		[J_1,p_1]=\big[ J_1,[q_0,q_1] \big]=-\big[ q_0,[q_1,J_1] \big]-\big[ q_1,[J_1,q_0] \big].
	\end{equation}
	The first terms vanishes because $q_1\in\sA$ while $[J_1,q_0]=-q_2$ by definition.
\end{proof}

\begin{proposition}		\label{LONGPropBJpsun}
	We have
	\begin{equation}
		\begin{aligned}[]
			B(J_1,J_1)&=-B(q_0,q_0)\\
			B(p_1,p_1)&=-B(q_0,q_0)\\
			B(s_1,s_1)&=B(q_0,q_0),
		\end{aligned}
	\end{equation}
	and then 
	\begin{equation}
		-\| J_1 \|^2=-\| p_1 \|^2=\| s_1 \|^2=1.
	\end{equation}
\end{proposition}

\begin{proof}
	Using $J_1=[q_0,q_2]$ (lemma \ref{LONGLemQzQdeuxJun}) we find
	\begin{equation}
		B(J_1,J_1)=B\big( \ad(q_2)q_0,\ad(q_2)q_0 \big)=-B\big( \ad(q_2)^2q_0,q_0 \big)=-B(q_0,q_0).
	\end{equation}
	In much the same way, using the definition of $p_1$ and $s_1=[q_1,q_2]$ (lemma \ref{LONGlemJDeuxqDeuxsUn}), we find $B(p_1,p_1)=B(q_1,q_1)$ and $B(s_1,s_1)=-B(q_2,q_2)$.

\end{proof}

\begin{lemma}		\label{LONGLemJunrkzero}
	We have $[J_1,r_k]=0$.
\end{lemma}

\begin{proof}
	Using the definition of $r_k$ and the Jacobi identity,
	\begin{equation}
		[J_1,r_k]= \big[ J_1,[J_2,q_k] \big]=-\big[ J_2,[q_k,J_1] \big]-\big[ q_k,[J_1,J_2] \big]=0
	\end{equation}
	because of equation \eqref{LONGEqJUnqkzero} and the fact that $\sA$ is abelian.
\end{proof}

Now, the Killing norms of the basis $\eB$ can be computed.
\begin{theorem}		\label{LONGThoBaisXXorthoigher}
	The basis
	\begin{equation}
		\eB=\{ q_0,q_2,p_1,s_1,q_k,p_k,r_k,s_k,J_1,J_2 \}.
	\end{equation}
	given by the definitions \eqref{LONGEqSuperBaseeB}  is orthonormal and
	\begin{equation}
		\begin{aligned}[]
			\| J_1 \|^2=\| q_1 \|^2=\| q_2 \|^2=\| p_1 \|^2=\| q_k \|^2=\| p_k \|^2&=-1\\
			\| r_{ij} \|=\| q_0 \|^2=\| s_1 \|^2=\| r_k \|^2=\| s_k \|^2&=1
			\end{aligned}
	\end{equation}
\end{theorem}

\begin{proof}
	The norms are given by the propositions \ref{LONGPropBprsk}, \ref{LONGPropBJpsun}, and \ref{LONGPropBaseQOrtho}.

	For the orthogonality, we know that $\sP\perp\sK$, $\sQ\perp\sH$ (from general theory) as well as $\tilde\sN_2\perp\tilde\sN_k$, $\tilde\sN_2\perp\sA$ and $\tilde\sN_k\perp\sA$ (equations \eqref{LONGEqtsnkperprtsnkp}, \eqref{LONGEqAperpNTrois}, \eqref{LONGEqAperpNk} and \eqref{LONGEqNTtroisperpNk}). Thus, among the elements of the basis $\eB$, the two only ones that could not orthogonal are $r_k$ and $s_k$. However, we have $B(r_k,s_k)=0$ because
	\begin{equation}
		B(r_k,s_k)=B\big( r_k,\ad(J_1)p_k \big)=-B\big( \ad(J_1)r_k,p_k \big)=0
	\end{equation}
	by lemma \ref{LONGLemJunrkzero}. 

\end{proof}

It turns out that we are able to compute all the commutators using the following techniques
\begin{enumerate}

	\item
		the orthonormality of the basis, theorem \ref{LONGThoBaisXXorthoigher} among with the $\ad$-invariance of the Killing form
	\item
		the Jacobi identity
	\item
		the theorem \ref{LONGThoAdESqqq} and the commutators of lemmas \ref{LONGLempunjdeuxqzero}, \ref{LONGlemJDeuxqDeuxsUn} and \ref{LONGLemJunrkzero}.

\end{enumerate}
Computing all the commutators that way is quite long and very few interesting. The interesting point is that it is possible. You can immediately jump to subsection \ref{LONGSubSecSomeExpo}.

\begin{enumerate}
	\item$\ad(q_0)J_1=q_2$. Definition.
	\item$\ad(q_0)J_2=p_1$. Definition.
	\item$\ad(p_1)J_1=-s_1$. Definition.
	\item$\ad(q_0)q_k=p_k$\label{LONGItemComqzpk}. Definition.
	\item$\ad(q_k)J_2=-r_k$. Definition.
	\item$\ad(p_k)J_1=-s_k$\label{LONGItemCompkJun}. Definition.
	\item$\ad(p_1)J_2=q_0$\label{LONGItemCompunJdeux}. Lemma \ref{LONGLempunjdeuxqzero}.
	\item$\ad(q_2)J_2=-s_1$. Lemma \ref{LONGlemJDeuxqDeuxsUn}.
	\item$\ad(J_1)r_k=0$\label{LONGItemComJunrk}. Lemma \ref{LONGLemJunrkzero}.
	
	\item$\ad(q_0)p_1=-J_2$\label{LONGItemComqzpun}. We have $\ad(q_0)p_1=\ad(q_0)^2q_1=-q_1$.
	\item$\ad(q_0)s_1=0$. By the usual techniques, we get $[q_0,s_1]\in\sK\cap\sQ\cap\sA=\{ 0 \}$. The following few are obtained in the same way.
	\item$\ad(q_0)r_k=0$\label{LONGItemComqzrk}.
	\item$\ad(q_0)s_k=0$.
	\item$\ad(q_2)p_1=0$.
	\item$\ad(q_k)J_1=0$.
	\item$\ad(q_2)p_k=0$.
	\item$\ad(q_2)p_k=0$.
	\item$\ad(p_1)q_k=0$\label{LONGItemCompunqk}.
	\item$\ad(J_2)p_k=0$.
	\item$\ad(q_2)s_1=-J_2$. We know that $[q_2,s_1]\in\sP\cap\sQ\cap\sA=\langle J_2\rangle$. Thus $[q_0,s_1]$ is a multiple of $J_2$. The coefficient is given by$B\big( [q_2,s_1],J_2 \big)/B(J_2,J_2)$. Using the $\ad$-invariance of the Killing form, we are left to compute $B(s_1,\ad(q_2)J_2)$. Lemma \ref{LONGlemJDeuxqDeuxsUn} shows then that the coefficient we are searching if $B(s_1,s_1)/B(J_2,J_2)=1$.
	\item\label{LONGItemComqdeuxqk}$\ad(q_2)q_k=-s_k$. The spaces show that $[q_2,q_k]\in\sK\cap\sH\cap\tilde\sN_k=\langle r_k,s_k\rangle$. Thus we have to check the two possible components. First, $B\big( [q_2,q_k],r_k \big)=-B\big( q_2,\ad(q_k)^2J_2 \big)=0$ by theorem \ref{LONGThoAdESqqq}. For the second, we use the definition of $s_k$ and the Jacobi identity:
		\begin{equation}
			\begin{aligned}[]
				B\big( [q_2,q_k],s_k \big)&=B\big( [q_2,q_k],\ad(J_1)[q_0,q_k] \big)\\
				&=B\big( [q_2,q_k],-\ad(q_0)\underbrace{[q_k,J_1]}_{=0}-\ad(q_k)\underbrace{[J_1,q_0]}_{=-q_2} \big)\\
				&=-B\big( \ad(q_2)q_k,\ad(q_2)q_k \big)\\
				&=B\big( q_k,\ad(q_2)^2q_k \big)\\
				&=B(q_k,q_k).
			\end{aligned}
		\end{equation}
		Finally, what we have is
		\begin{equation}
			[q_2,q_k]=\frac{ B(q_2,q_2) }{ B(s_k,s_k) }s_k=-s_k.
		\end{equation}
	\item$\ad(q_2)r_k=0$. From the spaces, $[q_2,r_k]\in\langle q_k\rangle$, but
		\begin{equation}
			B\big( \ad(q_2)r_k,q_k \big)=-B\big( r_k,[q_2,q_k] \big)=B(r_k,s_k)=0.
		\end{equation}
	\item$\ad(q_2)s_k=-q_k$. From the spaces, $[q_2,s_k]\in\sP\cap\sQ\cap\tilde\sN_k=\langle q_k\rangle$, so we compute
		\begin{equation}
			B\big( [q_2,s_k],q_k \big)=-B(s_k,[q_2,s_k])=B(s_k,s_k),
		\end{equation}
		and 
		\begin{equation}
			[q_2,s_k]=\frac{ B(s_k,s_k) }{ B(q_k,q_k) }q_k=-q_k.
		\end{equation}
	\item$\ad(q_2)r_k=0$. From the spaces, $[q_k,q_2]\in\sP\cap\sQ\cap\tilde\sN_k=\langle q_k\rangle$. Thus we compute
		\begin{equation}
			B\big( [q_2,r_k],q_k \big)=-B\big( r_k,[q_2,q_k] \big)=B(r_k,s_k)=0
		\end{equation}
		where we used the item \ref{LONGItemComqdeuxqk}.
	\item\label{LONGItemComjunpun}$\ad(J_1)p_1=s_1$. From the spaces, $[J_1,p_1]\in\langle s_1\rangle$. We have
		\begin{equation}
			B\big( [J_1,p_1],s_1 \big)=-B\big( p_1,[J_1,s_1] \big)=-B\big( p_1,\ad(J_1)^2p_1 \big)=-B(p_1,p_1),
		\end{equation}
		thus
		\begin{equation}
			[J_1,p_1]=-\frac{ B(p_1,p_1) }{ B(s_1,s_1) }s_1=s_1.
		\end{equation}
	\item\label{LONGItemComskqk}$\ad(s_k)q_k=-q_2$. From the spaces, $[s_k,q_k]\in\langle J_2,q_2\rangle$. Thus we have two Killing forms to compute. The first is
		\begin{equation}
			B\big( [s_k,q_k],J_2 \big)=B\big( s_k,[q_k,J_2] \big)=B(s_k,r_k)=0
		\end{equation}
		were we used the definition of $r_k$. The second is
		\begin{equation}
			B\big( [s_k,q_k],q_2 \big)=B\big( s_k,[q_k,q_2] \big)=B(s_k,s_k).
		\end{equation}
		where we used item \ref{LONGItemComqdeuxqk}. Thus we have
		\begin{equation}
			[s_k,q_k]=\frac{ B(s_k,s_k) }{ B(q_2,q_2) }q_2=-q_2
		\end{equation}
	\item$\ad(s_k)J_2=0$. From the spaces, $[s_k,J_2]\in\langle q_k\rangle$, but
		\begin{equation}
			B\big( \ad(s_k)J_2,q_k \big)=-B\big( J_2,[s_k,q_k] \big)=-B(J_2,q_2)=0
		\end{equation}
		where we used item \ref{LONGItemComskqk}.
	\item$\ad(q_2)s_k=-q_k$. From the spaces, $[q_2,s_k]\in\langle q_k\rangle$. We have
		\begin{equation}
			B\big( \ad(q_2)s_k,q_k \big)=-B\big( s_k,\ad(q_2)q_k \big)=B(s_k,s_k)
		\end{equation}
		where we used item \ref{LONGItemComqdeuxqk}. Thus
		\begin{equation}
			[q_2,s_k]=\frac{ B(s_k,s_k) }{ B(q_k,q_k) }q_k=-q_k.
		\end{equation}
	\item$\ad(p_1)s_1=-J_1$. From the spaces, $[p_1,s_1]\in\langle J_1\rangle$. We have
		\begin{equation}
			B\big( [p_1,s_1],J_1 \big)=-B\big( s_1,[p_1,J_1] \big)=B(s_1,s_1),
		\end{equation}
		where we used item \ref{LONGItemComjunpun}. Thus,
		\begin{equation}
			[p_1,s_1]=\frac{ B(s_1,s_1) }{ B(J_1,J_1) }J_1=-J_1.
		\end{equation}
	\item$\ad(p_1)r_k=p_k$\label{LONGItemCompunrk}. We use the definition of $r_k$ and Jacobi:
		\begin{equation}
			[p_1,r_k]=\big[ p_1,[J_2,q_k] \big]=-\big[ J_2,\underbrace{[q_k,p_1]}_{=0} \big]-\big[ q_k,\underbrace{[p_1,J_2]}_{=q_0} \big]=-[q_k,q_0]=p_k.
		\end{equation}
		where we used items \ref{LONGItemCompunqk}, \ref{LONGItemCompunJdeux} and \ref{LONGItemComqzpk}.
	\item$\ad(q_k)J_2=-r_k$\label{LONGItemComkJdeux}. From the spaces, $[q_k,J_2]\in\langle r_k,s_k\rangle$. First, we have
		\begin{equation}
			B\big( [q_k,q_2],s_k \big)=-B\big( J_2,\underbrace{[q_k,s_k]}_{=-q_2} \big)=0
		\end{equation}
		where we used item \ref{LONGItemComskqk}. Secondly we have
		\begin{equation}
			B\big( [q_k,J_2],r_k \big)=B\big( q_k,\underbrace{[J_2,r_k]}_{\ad(J_2)^2q_k} \big)=B(q_k,q_k)
		\end{equation}
		Thus
		\begin{equation}
			[q_k,J_2]=\frac{ B(q_k,q_k) }{ B(r_k,r_k) }r_k=-r_k.
		\end{equation}
	\item$\ad(p_1)p_k=r_k$\label{LONGItemCompunpk}. We use the definition of $p_k$ and the Jacobi identity:
		\begin{equation}
			\begin{aligned}[]
				[p_1,p_k]&=-\big[ q_0,\underbrace{[q_k,p_1]}_{=0} \big]-\big[ q_k,\underbrace{[p_1,q_0]}_{=J_2} \big]\\
					&=-[q_k,J_2]\\
					&=r_k
			\end{aligned}
		\end{equation}
		where we used items \ref{LONGItemCompunqk}, \ref{LONGItemComqzpun} and \ref{LONGItemComkJdeux}.
	\item$\ad(p_1)s_k=0$. From the spaces, $[p_1,s_k]\in\langle p_k\rangle$. We have
		\begin{equation}
			B\big( [p_1,s_k],p_k \big)=-B\big( s_k,\underbrace{[p_1,p_k]}_{=0} \big)=0
		\end{equation}
		where we used the item \ref{LONGItemCompunpk}.
	\item$\ad(s_1)q_k=0$. From the spaces, $[s_1,q_k]\in\langle q_k\rangle$. We have
		\begin{equation}
			B\big( [s_1,q_k],q_k \big)=B\big( s_1,[q_k,q_k] \big)=0.
		\end{equation}
	\item$\ad(s_1)p_k=0$\label{LONGItemComsunpk}. From the spaces, $[s_1,p_k]\in\langle p_k\rangle$. We have
		\begin{equation}
			B\big( [s_1,p_k],p_k \big)=B\big( s_1,[p_k,p_k] \big)=0.
		\end{equation}
	\item$\ad(s_1)r_k=s_k$. From the spaces, $[s_1,q_k]\in\langle r_k,s_k\rangle$. Using Jacobi,
		\begin{equation}
			\begin{aligned}[]
				[s_1,r_k]&=\big[ [J_1,p_1],r_k \big]\\
				&=-\big[ \underbrace{[p_1,r_k]}_{=p_k},J_1 \big]-\big[ \underbrace{[r_k,J_1]}_{=0},p_1 \big]\\
					&=-[p_k,J_1]\\
					&=s_k
			\end{aligned}
		\end{equation}
		where we used items \ref{LONGItemCompunrk}, \ref{LONGItemComJunrk} and the definition of $s_k$.
	\item$\ad(s_1)J_1=-p_1$\label{LONGItemComsunJun}. We have $[s_1,J_1]=-\ad(J_1)^2p_1=-p_1$ by theorem \ref{LONGThoAdESqqq}.
	\item$\ad(s_1)s_k=-r_k$\label{LONGItemComsunsk}. We use the definition of $s_k$ and Jacobi:
		\begin{equation}
			\begin{aligned}[]
				[s_1,s_k]&=\big[ s_1,[J_1,p_k] \big]\\
				&=-\big[ J_1,\underbrace{[p_k,s_1]}_{=0} \big]-\big[ p_k,\underbrace{[s_1,J_1]}_{=-p_1} \big]\\
					&=[p_k,p_1]\\
					&=-r_k
			\end{aligned}
		\end{equation}
		where we used items \ref{LONGItemComsunpk}, \ref{LONGItemComsunJun} and \ref{LONGItemCompunpk}.
	\item$\ad(J_2)s_k=0$. From the spaces, $[J_2,s_k]\in\langle q_k\rangle$. We have
		\begin{equation}
			B\big( [J_2,s_k],q_k \big)=B\big( J_2,\underbrace{[s_k,q_k]}_{=-q_2} \big)=0
		\end{equation}
		where we used item \ref{LONGItemComskqk}.
	\item$\ad(q_k)p_k=-q_0$. Using the definition of $p_k$ and the theorem \ref{LONGThoAdESqqq},
		\begin{equation}
			[q_k,p_k]=\big[ q_k[q_0,q_k] \big]=-\ad(q_k)^2q_0=-q_0.
		\end{equation}
	\item$\ad(q_k)r_k=-J_2$\label{LONGItemComqkrk}. Using the definition of $r_k$ and theorem \ref{LONGThoAdESqqq}, we have
		\begin{equation}
			\big[ q_k,[J_2,q_k] \big]=-\ad(q_k)^2J_2=-J_2.
		\end{equation}
	\item$\ad(p_k)r_k=-p_1$. Using the definition of $r_k$ and Jacobi,
		\begin{equation}
			\begin{aligned}[]
				[p_k,r_k]&=-\big[ pk,[q_0,q_k] \big]\\
				&=\big[ q_0,\underbrace{[q_k,r_k]}_{=-J_2} \big]+\big[ q_k,\underbrace{[r_k,q_0]}_{=0} \big]\\
					&=[J_2,q_0]\\
					&=-p_1
			\end{aligned}
		\end{equation}
		where we used the items \ref{LONGItemComqkrk}, \ref{LONGItemComqzrk} and the definition of $p_1$.
	\item$\ad(p_k)s_k=-J_1$. The spaces show that $[p_k,s_k]\in\langle J_1,p_1\rangle$. We have
		\begin{equation}
			B\big( [p_k,s_k],J_1 \big)=-B\big( s_k,\underbrace{[p_k,s_1]}_{=-s_k} \big)=B(s_k,s_k)
		\end{equation}
		and
		\begin{equation}
			B\big( [p_k,s_k],p_1 \big)=-B\big( s_k,\underbrace{[p_k,p_1]}_{=-r_k} \big)=0
		\end{equation}
		where we used the items \ref{LONGItemCompkJun} and \ref{LONGItemCompunpk}. Thus
		\begin{equation}
			[p_k,s_k]=\frac{ B(s_k,s_k) }{ B(J_1,J_1) }J_1=-J_1.
		\end{equation}
	\item$\ad(r_k)s_k=s_1$. From the spaces, $[r_k,s_k]\in\langle s_1\rangle$. We have
		\begin{equation}
			B\big( [r_k,s_k],s_1 \big)=B\big( r_k,\underbrace{[s_k,s_1]}_{=r_k} \big)=B(r_k,r_k)
		\end{equation}
		where we used item \ref{LONGItemComsunsk}. Thus
		\begin{equation}
			[r_k,s_k]=\frac{ B(r_k,r_k) }{ B(s_1,s_1) }s_1=s_1.
		\end{equation}
	\item$\ad(J_1)q_2=-q_0$. Using the definition of $q_2$ and theorem \ref{LONGThoAdESqqq},
		\begin{equation}
			[J_1,q_2]=-\ad(J_1)^2q_0=-q_0.
		\end{equation}
	\item$\ad(J_1)s_k=p_k$. From the spaces, $[J_1,s_k]\in\langle p_k\rangle$. We have
		\begin{equation}
			B\big( [J_1,s_k],p_k \big)=-B\big( s_k,\underbrace{[J_1,p_k]}_{=s_k} \big)=-B(s_k,s_k)
		\end{equation}
		where we used the definition of $s_k$.
	\item$\ad(J_2)p_1=-q_0$\label{LONGItemComJdeuxpun}. Using the definition of $p_1$ and theorem \ref{LONGThoAdESqqq}, $[J_2,p_1]=-\ad(q_1)^2q_0=-q_0$.
	\item$\ad(J_1)s_1=q_2$. We use the definition of $s_1$ and Jacobi:
		\begin{equation}
			\begin{aligned}[]
				[J_2,s_1]&=\big[ J_2,[J_1,p_1] \big]\\
				&=-\big[ J_1,\underbrace{[p_1,J_2]}_{=q_0} \big]-\big[ p_1,\underbrace{[J_2,q_1]}_{=0} \big]\\
				&=-[J_1,q_0]\\
				&=q_2
			\end{aligned}
		\end{equation}
		where we used item \ref{LONGItemComJdeuxpun} and the definition of $q_2$.
	\item$\ad(J_2)r_k=q_k$. Using the definition of $r_k$ and theorem \ref{LONGThoAdESqqq}, $[J_2,r_k]=\ad(q_1)^2q_k=q_k$.
\end{enumerate}
 
\begin{proposition}
	We have $\sH=[\sQ,\sQ]$.
\end{proposition}

\begin{proof}
	The inclusion $[\sQ,\sQ]\subset\sH$ is by construction. Now every elements in the basis \eqref{LONGAlignPremDefHH} can be expressed in terms of commutators in $\sQ$ because
	\begin{subequations}
		\begin{align}
			J_1&=[q_0,q_2]\\
			s_k&=[q_k,q_2]\\
			s_1&=[J_2,q_2]	\label{LONGEqJdqdsu}
		\end{align}
	\end{subequations}
\end{proof}

\subsection{Properties of the basis}

The basis $\eB$ is motivated by the fact that $\ad(q_0)^2q_k=-q_k$, so that $ e^{\ad(xq_0)}$ is easy to compute on $q_k$ and $p_k$. Moreover, $r_k$ and $s_k$ belong to $\sK$, so that $[q_0,r_k]=[q_0,s_k]=0$. The drawback of that decomposition is that the basis elements do not belong to $\sN$ or $\bar\sN$ while it will be useful to have basis elements in $\sN$ and $\bar\sN$, among other for theorem \ref{LONGThoOrbitesOuverttes}. 

At a certain point, we are going to compute the exponentials $ e^{\ad(xq_0)}X$ when $X$ runs over $\sN$ and $\bar\sN$. We are going to extensively use the commutation relations listed in \eqref{LONGsubEqsGenPySO}, \eqref{LONGSubEqsPlusPresPySO} and \eqref{LONGSubEqsThethaPySO}. A particular attention will be devoted to the projection over $\sQ$ which will be central in determining the open and closed orbits of $AN$ in $G/H$.

\begin{lemma}		\label{LONGLemDecomptsNDanseB}
	The decomposition of $\tilde\sN_2$ with respect to $\eB$ is
	\begin{equation}		\label{LONGSubeqsDecompXqps}
		\begin{aligned}[]
			X_{++}&=q_0-q_2-p_1-s_1\\
			X_{+-}&=q_0-q_2+p_1+s_1
		\end{aligned}
	\end{equation}
	and the decomposition of $\tilde\sN_k$ with respect to $\eB$ is
	\begin{subequations}
		\begin{align}
			X_{0+}^k & = q_k+r_k	&    X_{0-}^k&=-q_k+r_k		\label{LONGsubEqqrkXzpegal}	\\
			X_{+0}^k & = -p_k-s_k	&	X_{-0}^k & = p_k-s_k
		\end{align}
	\end{subequations}
	
\end{lemma}

\begin{proof}

	Using known commutator and the fact that $[\ad(J_1),\ad(J_2)]=0$ on $\tilde\sN_2$, we find the following commutators:
	\begin{subequations}			
		\begin{align}
			[J_1,q_0]&=-q_2	&&	&[J_2,q_0]&=-p_1			\label{LONGsubEqJunqzJdeuxqzmoinspun}\\
			[J_1,q_2]&=-q_0	&&	&[J_2,q_2]&=s_1\label{LONGsubEqJunqzJdeuxqzmoinspdeux}\\
			[J_1,p_1]&=s_1	&&	&[J_2,p_1]&=-q_0\\
			[J_1,s_1]&=p_1	&&	&[J_2,s_1]&=q_2.
		\end{align}
	\end{subequations}
	From these properties, we deduce that $q_0-q_2-p_1-s_1$ is proportional to $X_{++}$. Since, by definition, $q_0$ is the $\sK\sQ$-component of $X_{++}$, the proportionality factor is $1$. We also know that  $X_{+-}$ is proportional to $q_0-q_2+p_1+s1$. Since $q_0-q_2=(X_{++})_{\sQ}=(X_{+-})_{\sQ}$ (proposition \ref{LONGPropXmpXppqq}), the proportionality coefficient is $1$. Thus we have

	For the elements of $\tilde\sN_k$, the commutation relations give
	\begin{subequations}
		\begin{align}
			[J_1,q_k+r_k]&=0		&[J_1,p_k-s_k]&=s_k-p_k\\
			[J_2,q_k+r_k]&=q_k+r_k		&[J_2,p_k-s_k]&=0
		\end{align}
	\end{subequations}
	so that
	\begin{equation}
		\begin{aligned}[]
			q_k+r_k&\propto X_{0+}^k\in\sN\\
			s_k+p_k&\propto X^k_{+0}\in\sN\\
			p_k-s_k&\propto X^k_{-0}\in\bar\sN
		\end{aligned}
	\end{equation}

	We have $r_k=[J_2,q_k]\in\sK\cap\sH$, so that the $\sP$-component of $q_k+r_k$ is $q_k$. But $q_k=(X_{0+}^k)_{\sP}$ is the $\sP$-component of $X_{0+}^k$. The proportionality between $q_k+r_k$ and $X_{0+}^k$ together with the equality of their $\sP$-component provide the equality \eqref{LONGsubEqqrkXzpegal}.

	For the two other, let us suppose that
	\begin{subequations}
		\begin{align}
			X_{+0}^k&=a(p_k+s_k)\\
			X_{-0}^k&=b(p_k-s_k).
		\end{align}
	\end{subequations}
	In this case, we have 
	\begin{equation}
			(X_{+0}^k)_{\sP}=\frac{ 1 }{2}(X_{+0}^k-\theta X_{+0}^k)
				=\frac{ 1 }{2}\big( (a-b)p_k+(a+b)s_k \big),
	\end{equation}
	so that $a=-b$ because $s_k\in \sK$. Now let us look at the $\sK\sQ$-component of the equality $[X_{+0}^k,X_{0+}^k]=-X_{++}$ taking into account the fact that $X_{+0}^k\in\sH$ and $(X_{0+}^k)_{\sK\sQ}=0$. What we have is $\big[ (X_{+0}^k)_{\sP\sH},(X_{0+}^k)_{\sP\sQ} \big]=-q_0$, but $(X_{0+}^k)_{\sP}=q_k$ and $(X_{+0}^k)_{\sP}=ap_k$, so that $[ap_k,q_k]=-q_0$. If we replace $p_k$ by its definition $[q_0,q_k]$, we get
	\begin{equation}
		a\big[ [q_0,q_k],q_k \big]=a\ad(q_k)^2q_0=-q_0,
	\end{equation}
	so that $a=-1$. 

	The last point comes from $X_{0-}^k=\theta X_{0+}^k$.

\end{proof}
Notice that this result was already obvious from the decompositions given in \eqref{LONGEqDecompqprskXX}.

\subsection{Some exponentials}
\label{LONGSubSecSomeExpo}

It will be important to compute the element $ e^{\ad(xq_0)}X$ when $X$ runs over the vectors of $\eB$. 
The action of $e^{xq_0}$ on $\sA$ is
\begin{subequations}			\label{LONGSubEqsAdxqzJJ}
	\begin{align}
		e^{\ad(xq_0)}J_1=\cos(x)J_1+\sin(x)q_2\\
		e^{\ad(xq_0)}J_2=\sin(x)p_1+\cos(x)q_1,
	\end{align}
\end{subequations}
on $\tilde\sN_2$ we have
\begin{subequations}					\label{LONGSubEqsexpxzqpsdzuu}
	\begin{align}
		e^{\ad(xq_0)}q_2&=\cos(x)q_2-\sin(x)J_1\\
		e^{\ad(xq_0)}q_0&=q_0\\
		e^{\ad(xq_0)}p_1&=\cos(x)p_1-\sin(x)q_1 		\label{LONGEqAdqzpUn}\\
		e^{\ad(xq_0)}s_1&=s_1,
	\end{align}
\end{subequations}
and the action on $\tilde\sN_k$ is
\begin{subequations}		\label{LONGEqExpAdqkpk}
	\begin{align}
		e^{\ad(xq_0)}q_k&=\cos(x)q_k+\sin(x)p_k\\
		e^{\ad(xq_0)}p_k&=\cos(x)p_k-\sin(x)q_k\\
		e^{\ad(xq_0)}p_k&=\cos(x)p_k-\sin(x)q_k.
	\end{align}
\end{subequations}

Combining with lemma \ref{LONGLemDecomptsNDanseB},
\begin{subequations}			\label{LONGEqExpoQzSurNk}
	\begin{align}
		e^{\ad(xq_0)}X_{0+}^k=e^{\ad(xq_0)}(q_k+r_k)&=r_k+\cos(x)q_k+\sin(x)p_k		\label{LONGEqExpoQzSurNka}\\
		e^{\ad(xq_0)}X_{+0}^k=-e^{\ad(xq_0)}(s_k+p_k)&=-s_k-\cos(x)p_k+\sin(x)q_k.\label{LONGEqExpoQzSurNkb}
	\end{align}
\end{subequations}
The same way,
\begin{subequations}			\label{LONGEqExpoQzN}
	\begin{align}
		e^{\ad(xq_0)}X_{++}&=q_0+\sin(x)q_1-\cos(x)q_2+\sin(x)J_1-\cos(x)p_1\\
		e^{\ad(xq_0)}X_{+-}&=q_0-\sin(x)q_1-\cos(x)q_2+\sin(x)J_1+\cos(x)p_1.
	\end{align}
\end{subequations}
The projections on $\sQ$ of all these combinations are immediate.

\subsection{Classification of the basis by the spaces}

The basis \eqref{LONGEqDecompqprskXX} allows to  generalize the theorem \ref{LONGThoAdESqqq}.

By very definition, we have $\ad(\sP)^2\sP\subset\sP$, $\ad(\sP)^2\sK\subset\sK$ and the same for the couple $(\sH,\sQ)$. The relations \eqref{LONGEqsCommWithtsNDeuxkA} say that the same is true with the triple $(\tilde\sN_2,\tilde\sN_k,\sA)$, i.e. $\ad(\sA)^2\tilde\sN_k\subset\tilde\sN_k$ and $\tilde\sN_k(\tilde\sN_2)\subset\tilde\sN_2$.

\begin{theorem}		\label{LONGThoAdSqIouZero}
	The basis $\eB$ has the property to be stable under the commutators: $[X,Y]\in\{0,\pm\eB\}$ when $X,Y\in\eB$. Moreover, we have
	\begin{equation}		\label{LONGEqadXsqYzYmY}
		\ad(X)^2Y=\begin{cases}
			0	&	\text{if $\ad(X)Y=0$}\\
			Y	&	\text{if $X\in\sP$}\\
			-Y	&	 \text{if $X\in\sK$}.
		\end{cases}
	\end{equation}
\end{theorem}

\begin{proof}
	This theorem can be immediately checked using the commutators. It is however instructive to see that equation \eqref{LONGEqadXsqYzYmY} can be checked from few considerations. First, remark that, if we look at the decompositions $\sG=\sQ\oplus\sH=\sP\oplus\sK=\mZ_{\sK}(\sA)\oplus\sA\oplus\tilde\sN_2\oplus\tilde\sN_k$, the commutation relations \eqref{LONGEqsCommWithtsNDeuxkA}, we find
	\begin{equation}
		\ad(\tilde\sN_2)\circ\ad(\tilde\sN_2)\colon
		\begin{cases}
			\tilde\sN_2\to\sA\to\tilde\sN_2\\
			\tilde\sN_k\to\tilde\sN_k\to\tilde\sN_k\\
			\sA\to\tilde\sN_2\to\sA,
		\end{cases}	
	\end{equation}
	\begin{equation}
		\ad(\tilde\sN_k)\circ\ad(\tilde\sN_k)\colon
		\begin{cases}
			\tilde\sN_2\to\tilde\sN_k\to\sA\oplus\tilde\sN_2\\
			\tilde\sN_k\to(\sA\oplus\tilde\sN_2)\to\sA\oplus\tilde\sN_k\\
			\sA\to\tilde\sN_k\to\sA\oplus\tilde\sN_2,
		\end{cases}
	\end{equation}
	\begin{equation}
		\ad(\sA)\circ\ad(\sA)\colon
		\begin{cases}
			\tilde\sN_2\to\tilde\sN_2\to\tilde\sN_2\\
			\tilde\sN_k\to\tilde\sN_k\to\tilde\sN_k\\
			\sA\to 0.
		\end{cases}
	\end{equation}
	Thus we have $\ad(X)^2Y=\lambda Y$ whenever we are not in the cases $(X,Y)\in(\tilde\sN_k,\tilde\sN_2)$ and $(X,Y)\in(\tilde\sN_k,\tilde\sN_k)$. We should check that these cases cannot bring a $\sA$-component. We should also check that, since $\sK\cap\sH\cap\tilde\sN_k$, is two dimensional, $\ad(X)^2r_k$ has no $s_k$-component as well as $\ad(X)^2s_k$ has no $r_k$-component.

	Let us suppose that the spaces fit. We have $\ad(X)^2Y=\lambda Y$ and we still have to check the values of $\lambda$. Since $\eB$ is orthonormal, we have
	\begin{equation}
		\lambda=\frac{ B\big( \ad(X)^2Y,Y \big) }{ B(Y,Y) }=-\frac{ B\big( \ad(X)Y,\ad(X)Y \big) }{ B(Y,Y) }.
	\end{equation}
	First, remark that this is zero if and only if $\ad(X)Y=0$. Now, there are $4$ possibilities following that $X,Y$ belong to $\sP$ or $\sK$ because we know that $B|_{\sK}<0$ and $B|_{\sP}>0$. The result is that $\lambda$ is positive when $X\in\sP$ and negative when $X\in\sK$. Now, the fact that $\eB$ is stable under the commutators implies that $\ad(X)^2Y\in\{ 0,\pm \eB \}$, so that $\lambda\in\{ 0,1,-1 \}$.
\end{proof}
\section{The causally singular structure}
\label{LONGSecBlacHole}

\subsection{Closed orbits}

The singularity in $AdS_l$ is defined as the closed orbits of $AN$ and $A\bar N$ in $G/H$. This subsection is intended to identify them.

\begin{proposition}		\label{LONGPropCartanExtExpo}
	The Cartan involution $\theta\colon \sG\to \sG$ is an inner automorphism, namely it is given by $\theta=\Ad(k_{\theta})$ where $k_{\theta}= e^{\pi q_0}$.
\end{proposition}

\begin{proof}
	The operator $\Ad(k_{\theta})$ acts as the identity on $\sK$ because $q_0$ is central in $\sK$ by definition.  Looking at the decompositions \eqref{LONGEqDecomptsNkKP} and \eqref{LONGEqDecomptsNTroisKP}, and taking into account that the result is already guaranteed on $\sK$, we have to check the action of $\Ad(k_{\theta})$ on $J_1$, $J_2$, $q_k$, $p_k$ and $p_1$. It is done in setting $x=\pi$ in equations
	\eqref{LONGSubEqsAdxqzJJ},
	\eqref{LONGEqExpAdqkpk}
	and \eqref{LONGEqAdqzpUn}.
	What we get is that $\Ad(k_{\theta})$ changes the sign on $\sP$.
	
\end{proof}

\begin{proposition}		\label{LONGPropKanUnicAbarN}
	For each $an\in AN$, there exists one and only one $k\in K$ such that $kan\in A\bar N$. There also exists one and only one $k\in K$ such that $ank\in A\bar N$.
\end{proposition}

\begin{proof}
	For unicity, let $an\in AN$ and suppose that $k_1^{-1}an$ and $k_2^{-1}an$ both belong to $A\bar N$. Then there exist $a_1$, $a_2$, $\bar n_1$ and $\bar n_2$ such that $k_1^{-1}an=a_1\bar n_1$ and $k_2^{-1}an=a_2\bar n_2$ and we have
	\begin{equation}
		an=k_1a_1\bar n_1=k_2a_2\bar n_2.
	\end{equation}
	By unicity of the decomposition $KA\bar N$, we conclude that $k_1=k_2$.

	For the existence, let $an\in AN$ and consider the $KAN$ decomposition $\theta(an)=ka'n'$. We claim that $k^{-1}$ answers the question. Indeed, $\theta$ is the identity on $K$, so that $an=k\theta(a'n')$, and then
	\begin{equation}
		k^{-1}an=\theta(a'n')\in A\bar N.
	\end{equation}
	
	One checks the statement about $ank\in A\bar N$ in much the same way.
\end{proof}%
In the following results, we use the fact that the group $K$ splits into the commuting product $K=\SO(2)\times\SO(l-1)$.
\begin{lemma}		\label{LONGLemExistxTqansAbarN}
	For every $an\in AN$, there exists $x\in\mathopen[ 0 , 2\pi [$ such that $[an e^{xq_0}]\in[A\bar N]$.
\end{lemma}

\begin{proof}
	Let $k\in K$ such that $ank\in A\bar N$. The element $k$ decomposes into $k=st$ with $s= e^{xq_0}\in \SO(2)$ and $t\in\SO(n)\subset H$. Thus $[ans]\in[A\bar N]$.
\end{proof}

\begin{lemma}		\label{LONGLemansse}
	If $[an]=[s]$ with $s\in\SO(2)$, then $s=e$.
\end{lemma}

\begin{proof}
	The assumption implies that there exists a $h\in H$ such that $an=sh$. Using the $KAN$ decomposition of $H$, such a $h$ can be written under the form $h=ta'n'$ with $t\in\SO(n)$. Thus we have $an=sta'n'$. By unicity of the decomposition $kan$, we must have $st=e$, and then $s=e$.
\end{proof}

\begin{lemma}		\label{LONGLemANksk}
	If $an\in AN$ and if $[ank_{\theta}]=[s]$ with $s\in\SO(2)$, then $s=k_{\theta}$.
\end{lemma}

The hypothessis implies that there exists a $h\in H$ such that $ank_{\theta}=sh$. The element $h$ can be decomposed as
\begin{equation}
	h=h_Kh_{AN}k_{\theta}
\end{equation}
with $h_K\in K$ and $h_{AN}\in AN$ but with no warranty that $h_K$ or $h_{AN}$ belong to $H$. We have
\begin{equation}
	an=\underbrace{sh_K}_{\in K}h_{AN},
\end{equation}
so that $sh_K=e$ by unicity of the $KAN$ decomposition. In particular $h_K\in\SO(2)$.

Now we want to prove that $h_K=k_{\theta}$ because it would implies $s=k_{\theta}$ by the relation $sh_K=e$. If $h_K= e^{yq_0}k_{\theta}$ we have
\begin{equation}
	h=h_Kh_{AN}k_{\theta}= e^{yq_0}k_{\theta}h_{AN}k_{\theta}= e^{yq_0}h_{A\bar N}.
\end{equation}
where $h_{A\bar N}=\AD(k_{\theta})h_{AN}$. By unicity of the $KA\bar N$ decomposition in $H$, the latter relation shows that $ e^{yq_0}$ has to be the compact part of $H$ and then belong to $K_H$ which is only possible when $y=0$.

\begin{theorem}		\label{LONGThoOrbitesOuverttes}
	The closed orbits of $AN$ in $AdS_l$ are $[AN]$ and $[AN k_{\theta}]$ where $k_{\theta}$ is the element of $K$ such that $\theta=\Ad(k_{\theta})$. The closed orbits of $A\bar N$ are $[A\bar N]$ and $[A\bar N k_{\theta}]$. The other orbits are open.
\end{theorem}

\begin{proof}
	Let us deal with the $AN$-orbits in order to fix the ideas. First, remark that each orbit of $AN$ pass trough $[SO(2)]$. Indeed, each $[ank]$ is in the same orbit as $[k]$ with $k\in K=\SO(2)\otimes\SO(n)$. Since $\SO(n)\subset H$, we have $[k]=[s]$ for some $s\in\SO(2)$.

	We are thus going to study openness of the $AN$-orbit of elements of the form $[e^{x q_0}]$ because these elements are ``classifying'' the orbits. Using the isomorphism $  dL_{g^{-1}}\colon T_{[g]}(G/H)\to \sQ$, we know that a set $\{ X_1,\ldots X_l \}$ of vectors in $T_{[ e^{x q_0}]}AdS_l$ is a basis if and only if the set $\{ dL_{ e^{-xq_0}}X_i \}_{i=1,\ldots l}$ is a basis of $\sQ$. We are thus going to study the elements 
	\begin{equation}
		\begin{aligned}[]
			dL_{ e^{-xq_0}}X^*_{[ e^{xq_0}]}	&=dL_{ e^{-xq_0}}\Dsdd{ \pi\big(  e^{-t X} e^{xq_0} \big) }{t}{0}\\
								&=\Dsdd{ \pi\big(  \AD( e^{-xq_0} ) e^{-tX} \big) }{t}{0}\\
								&=-\pr_{\sQ} e^{-\ad(xq_0)}X
		\end{aligned}
	\end{equation}
	when $X$ runs over the elements of $\sA\oplus \sN$. 
	The projections on $\sQ$ of equations \eqref{LONGSubEqsAdxqzJJ}, \eqref{LONGEqExpoQzSurNk} and \eqref{LONGEqExpoQzN} are
	\begin{subequations}		\label{LONGSubEqsExpAdxqDivers}
		\begin{align}
			\pr_{\sQ}\Big( e^{\ad(xq_0)}J_1\Big)&=\sin(x)q_2\\
			\pr_{\sQ}\Big( e^{\ad(xq_0)}J_2\Big)&=\cos(x)q_1\\
			\pr_{\sQ}\Big(  e^{xq_0}X_{++} \Big)&=q_0+\sin(x)q_1-\cos(x)q_2\\
			\pr_{\sQ} \Big(   e^{\ad(xq_0)}X_{+-} \Big)&=q_0-\sin(x)q_1-\cos(x)q_2\\
			\pr_{\sQ}\Big(  e^{\ad(xq_0)}(s_k-p_k) \Big)&= \sin(x)q_k\\
			\pr_{\sQ}\Big(  e^{\ad(xq_0)}(q_k+r_k) \Big)&=\cos(x)q_k.
		\end{align}
	\end{subequations}
	It is immediately visible
	that an orbit trough $[ e^{xq_0}]$ is open if and only if $\sin(x)\neq 0$. It remains to study the orbits of $[ e^{\pi q_0}]$ and $[e]$. Lemma \ref{LONGLemansse} shows that these two orbits are disjoint.

	Let us now prove that $[AN]$ is closed. A point outside $\pi(AN)$ reads $\pi(ans)$ where $s$ is an elements of $\SO(2)$ which is not the identity. Let $\mO$ be an open neighborhood of $ans$ in $G$ such that every element of $\mO$ read $a'n's't'$ with $s'\neq e$. The set $\pi(\mO)$ is then an open neighborhood of $\pi(ans)$ which does not intersect $[AN]$. This proves that the complementary of $[AN]$ is open. The same holds for the orbit $[A\bar N]$.
	
	The orbit $[ANk_{\theta}]$ and $[A\bar Nk_{\theta}]$ are closed too because $ANk_{\theta}=k_{\theta}A\bar N$.

\end{proof}

\begin{lemma}
	We have $[AN]\cap[ANk_{\theta}]=\emptyset$.
\end{lemma}

\begin{proof}
	A representative of an element in $[AN]\cap[ANk_{\theta}]$ can be written $an=a'n'k_{\theta}h$, and then we have $k_{\theta}h\in AN$. Decomposing $h$ into its components $h_K\in K$ and $h_{AN}\in AN$, we see that $k_{\theta}h_K\in AN$, which is impossible because $k_{\theta}\in\SO(2)$ while $K_H=\SO(l-1)$.
\end{proof}

\begin{proposition}		\label{LONGPropUniquexxxxANANbarktheta}
	Let $an\in AN$. Then there exist unique $x_0$, $x_1$, $x_2$, $x_3$ in $\mathopen[ 0 , 2\pi [$ such that
	\begin{enumerate}
		\item
			$[an e^{x_0q_0}]\in[AN]$, and $an e^{x_0q_0}=an$,
		\item
			$[an e^{x_1q_0}]\in[ANk_{\theta}]$, and $an e^{x_1q_0}=ank_{\theta}$,
		\item
			$[an e^{x_2q_0}]\in[A\bar N]$, and $an e^{x_2q_0}=a'\bar n h_K$ (lemma \ref{LONGLemExistxTqansAbarN}),
		\item
			$[an e^{x_3q_0}]\in[A\bar Nk_{\theta}]$, and $an e^{x_3q_0}=a'\bar nk_{\theta}h_K$.
	\end{enumerate}
	These numbers satisfy $x_0=0$, $x_1=\pi$ and $x_3=x_2+\pi$ modulo $2\pi$. 

	Moreover if $P$ does not belong to $[AN]\cap[A\bar N]$ nor to $[AN]\cap[A\bar Nk_{\theta}]$, these are four different numbers.
\end{proposition}

\begin{proof}
	 Existence comes from lemma \ref{LONGLemExistxTqansAbarN}. Now we discuss the unicity.
	 \begin{enumerate}

		 \item
			 Let $x'_0$ such that $[an e^{x'_0q_0}]\in[AN]$. In that case, there exist $a'n'\in AN$ and $h\in H$ such that $ e^{x'_0q_0}=a'n'$. If we decompose $h=h_{AN}h_K$, unicity of the $ANK$ decomposition show that $ e^{x'_0q_0}=h_K$, which is only possible when $x'_0=0$ and $h_K=e$.
		\item
			If $an e^{x'_1q_0}=a'n'k_{\theta}h$, then there exist $a''n''$ such that $ e^{x'_1q_0}=a''n''k_{\theta}h$. Taking the class, $[a''n''k_{\theta}]=[ e^{x'_1q_0}]\in[\SO(2)]$. Then $ e^{x'_1q_0}=k_{\theta}$ by lemma \ref{LONGLemANksk}.
		\item
			Let $x_2$ such that $an e^{x_2q_0}t=a'\bar n$ with $t\in\SO(l-1)$ (proposition \ref{LONGPropKanUnicAbarN}), and suppose that we have an other $x'_2$ such that $[an e^{x'_2}]\in[A\bar N]$. We have the system
			\begin{subequations}
				\begin{numcases}{}
					an e^{x_2q_0}t=a'\bar n\\
					an e^{x'_2q_0}=a''\bar n' h.
				\end{numcases}
			\end{subequations}
			We extract $an$ from the first equation and we put the result in the second one. Taking into account the fact that $t$ commutes with $ e^{xq_0}$ and that $t\in H$ and renaming $h\to ht$, if we decompose $h=h_{A\bar N}h_K$ we find
			\begin{equation}
				a'\bar n e^{(x'_2-x_2)q_0}=a''\bar n'h_{A\bar N}h_K.
			\end{equation}
			By unicity of the decomposition $A\bar NK$ we find $h_K= e^{(x'_2-x_2)q_0}$. The left hand side belongs to $\SO(l-1)$ and the right hand side belongs to $\SO(2)$, so that one has to have $x'_2=x_2$.
		\item
			If we consider $x_3=x_2+\pi$ we have $an e^{x_3q_0}t=a'\bar nk_{\theta}$ (with $t\in\SO(l-1)$) and if we consider an other $x'_3$ such that $[an e^{x'_3q_0}]\in[A\bar Nk_{\theta}]$, we have the system
			\begin{subequations}
				\begin{numcases}{}
					an e^{(x_2+\pi)q_0}t=a'\bar nk_{\theta}\\
					an e^{x'_3q_0}=a'' \bar n'k_{\theta}h.
				\end{numcases}
			\end{subequations}
			Some manipulations including a redefinition of $h$ to include $t$ yield
			\begin{equation}
				e^{(x'_3-x_3)q_0}=k_{\theta}a_0\bar n_0k_{\theta}h,
			\end{equation}
			but $k_{\theta}A\bar Nk_{\theta}\subset AN$, thus
			\begin{equation}
				e^{(x'_3-x_3)q_0}=a_1n_1h_{AN}h_K
			\end{equation}
			so that unicity of the decomposition $ANK$ implies, $ e^{(x'_3-x_3)q_0=h_K}$, so that $x'_3=x_3$. 
	\end{enumerate}
	Let us now prove the second part. We suppose that $[P]$ does not belong to $[AN]\cap[A\bar N]$. If $x_0=x_2=0$, we have $[P]\in[AN]\cap[A\bar N]$. If $x_0=x_3$, then we have $an=a'\bar nk_{\theta}h_K$, and then $[an]\in[AN]\cap[A\bar Nk_{\theta}]$ since $h_K$ commutes with $k_{\theta}$.

	 If $x_1=x_2$, then $an=a'\bar nk_{\theta}h_K$, so that $[an]\in[AN]\cap[A\bar Nk_{\theta}]$. If $x_1=x_3$, then $ank_{\theta}=a'\bar nk_{\theta}h_K$, so that $an\in[AN]\cap[A\bar N]$.
\end{proof}

\subsection{Vanishing norm criterion}

In the preceding section, we defined the singularity by means of the action of an Iwasawa group. We are now going to give an alternative way of describing the singularity, by means of the norm of a fundamental vector of the action. This ``new'' way of describing the singularity is, in fact, much more similar to the original BTZ black hole where the singularity was created by identifications along the integral curves of a Killing vector field\cite{these_Detournay}. The vector $J_1$ in theorem \ref{LONGThosSequivJzero} plays here the role of that ``old'' Killing vector field.

Discrete identifications along the integral curves of $J_1$ would produce the causally singular space which is at the basis of our black hole.

What we will prove is the
\begin{theorem}		\label{LONGThosSequivJzero}
	We have $\sS\equiv \| J_1^* \|_{[g]}=\| \pr_{\sQ}\Ad(g^{-1})J_1 \|=0$.
\end{theorem}

Thanks to this theorem, our strategy will be to compute $\| \pr_{\sQ}\Ad(g^{-1})J_1 \|$ in order to determine if $[g]$ belongs to the singularity or not. The proof will be decomposed in three steps. The first step is to obtain a manageable expression for $\| J_1^* \|$.

\begin{lemma}		\label{LONGLemExpressionCoolNormJUn}
	We have $\| (J_1^*)_{[g]} \|=\| \pr_{\sQ}\Ad(g^{-1})J_1 \|$ for every $[g]\in AdS_l$.
\end{lemma}

\begin{proof}
	By definition, 
	\begin{equation}
		(J_1^*)_{[g]}=\Dsdd{ \pi( e^{-tJ_1}g) }{t}{0}=-d\pi dR_g J_1.
	\end{equation}
	The norm of this vector is the norm induced from the Killing form on $\sG$\cite{Kerin}. First we have to put $dR_g J_1$ under the form $dL_g X$ with $X\in \lG$. One obviously has $dR_g J_1=dL_g\Ad(g^{-1})J_1$, and the norm to be computed is
	\begin{equation}
		\begin{aligned}[]
			\| J_1^* \|_{[g]}=\| d\pi_gdL_g\Ad(g^{-1})J_1\|_{[g]}&=  \|  d\pi_gdL_g\pr_{\sQ}\Ad(g^{-1})J_1\|_{[g]}\\ 
			&= \| dL_g\pr_{\sQ}\Ad(g^{-1})J_1 \|_{g} \\
			&=\| \pr_{\sQ}\Ad(g^{-1})J_1 \|_e
		\end{aligned}
	\end{equation}
\end{proof}

\begin{proposition}		\label{LONGPropPtpsSjzero}
	If $p\in\hS$, then $\| J_1^* \|_p=0$.
\end{proposition}

\begin{proof}
	We are going to prove that $\pr_{\sQ}\Ad(g^{-1})J_1$ is a light like vector in $\sQ$ when $g$ belongs to $[AN]$ or $[A\bar N]$. A general element of $AN$ reads $g=a^{-1}n^{-1}$ with $a\in A$ and $n\in N$. Since $\Ad(a)J_1=J_1$, we have $\Ad(g^{-1})J_1=\Ad(n)J_1$. We are going to study the development
	\begin{equation}
		\Ad( e^{Z})J_1= e^{\ad(Z)}J_1=J_1+\ad(Z)J_1+\frac{ 1 }{2}\ad(Z)^2J_1+\ldots
	\end{equation}
	where $Z=\ln(n)\in\sN$. The series is finite because $Z$ is nilpotent (see theorem \ref{LONGEtOrdreDeux} for more informations) and begins by $J_1$ while all other terms belong to $\sN$. Notice that the same remains true if one replace $\sN$ by $\bar \sN$ everywhere. 
	
	Moreover, $\Ad( e^{Z})J_1$ has no $X_{0+}$-component (no $X_{0-}$-component in the case of $Z\in\bar\sN$) because $[X_{0+},J_1]=0$, so that the term $[Z,J_1]$ is a combination of $X_{+0}$, $X_{++}$ and $X_{+-}$.  Since the action of $\ad(X_{+\pm})$ on such a combination is always zero, the next terms are produced by action of $\ad(X_{0+})$ on a combination of $X_{+0}$, $X_{++}$ and $X_{+-}$. Thus we have
\begin{equation}		\label{LONGEqAdanJUnabck}
	\Ad( e^{Z})J_1=J_1+aX_{++}+bX_{+-}+c_kX_{+0}^k
\end{equation}
for some\footnote{One can show that every combinations of these elements are possible, but that point is of no importance here.} constants $a$, $b$ and $c_k$.

The projection of $\Ad( e^{Z})J_1$ on $\sQ$ is made of a combination of the projections of $X_{+0}$, $X_{++}$ and $X_{+-}$. From the definitions \eqref{LONGEqBasQQzi}, we have $\pr_{\sQ}X_{++}=q_0+q_2$, lemma \ref{LONGLemNonHXaz} implies $\pr_{\sQ}X_{+0}=0$ and lemma \ref{LONGLemSigmaXppEgalXPm} yields $\pr_{\sQ}X_{+-}=-\sigma\pr_{\sQ}X_{++}=q_0+q_2$. The conclusion is that $\pr_{\sQ}\big( e^{\ad(Z)}J_1\big)$ is a multiple of $q_0+q_2$, which is light like. The conclusion still holds with $\bar\sN$, but we get a multiple of $q_0-q_2$ instead of $q_0+q_2$.

	Now we have $\Ad(k_{\theta})J_1=J_1$ and $\Ad(k_{\theta})(q_0\pm q_2)=-(q_0\pm q_2)$, so that the same proof holds for the closed orbits $[ANk_{\theta}]$ and $[A\bar N k_{\theta}]$.
\end{proof}

\begin{remark}		\label{LONGRemANANbarYapas}
	The coefficients $a$, $b$ and $c_k$ in equation \eqref{LONGEqAdanJUnabck} are continuous functions of the starting point $an\in AN$. More precisely, they are polynomials in the coefficients of $X_{++}$, $X_{+-}$, $X_{0+}$ and $X_{0+}$ in $Z$. The vector $\pr_{\sQ}\Ad(g^{-1})J_1=(a+b)(q_0+q_1)$ is thus a continuous function of the point $[g]\in[AN]$.

	If $[g]\in[AN]\cap[A\bar N]$, then $\pr_{\sQ}\Ad(g^{-1})J_1$ has to vanish as it is a multiple of $q_0+q_1$ and of $q_0-q_1$ in the same time. We conclude that in each neighborhood in $[AN]$ of an element of $[AN]$, there is an element which does not belong to $[A\bar N]$.
\end{remark}

\begin{proposition}
	If $\| J_1^* \|_p=0$, then $p\in\sS$.
\end{proposition}

\begin{proof}
	As before we are looking at a point $[g]=[(an)^{-1}s^{-1}]$ with $s= e^{xq_0}$. The norm $\| J_1^* \|$ vanishes if
	\begin{equation}
		\| \pr_{\sQ} \Ad( e^{xq_0})\Ad(an)J_1 \|=0.
	\end{equation}
	We already argued in the proof of proposition \ref{LONGPropPtpsSjzero} that $\Ad(an)J_1$ is equal to $J_1$ plus a linear combination of $X_{++}$, $X_{+-}$ and $X_{+0}$. Using the relations \eqref{LONGSubEqsExpAdxqDivers}, we see that
	\begin{equation}		\label{LONGEqprQexpxqzXanroots}
		\begin{aligned}[]
			\pr_{\sQ} e^{\ad(xq_0)}(J_1&+aX_{++}+bX_{+-}+\sum_k c_kX^k_{+0})\\
								&=(a+b)q_0+(a-b)\sin(x)q_1+\big( \sin(x)-(a+b)\cos(x) \big)q_2+\sum _k c_k\sin(x)q_k.
		\end{aligned}
	\end{equation}
	The norm of this vector, as function of $x$, is given by
	\begin{equation}
		n(x)=(a+b)\sin(2x)+(4ab-C^2-1)\big( 1-\cos(2x) \big),
	\end{equation}
	where $C^2=\sum_kc_k^2$. Using the variables with $u=a+b$ and $v=(1+c^2-4ab)/2$,
	\begin{equation}
		n(x)=u\sin(2x)+v\cos(2x)-v.
	\end{equation}
	Following $u=0$ or $u\neq 0$, the graph of that function vanishes two or four times between $0$ and $2\pi$, see figure \ref{LONGLabelFigGraphn}. Points of $[AN]$ are divided into two parts: the \emph{red points} which correspond to $u\neq 0$, and the \emph{blue points} which correspond to $u=0$. By continuity, the red part is open.

		\newcommand{\CaptionFigGraphn}{In red, the function $n(x)$ with $u\neq 0$ and in blue, the function with $u=0$.}
\begin{figure}[ht]
\centering
\includegraphics{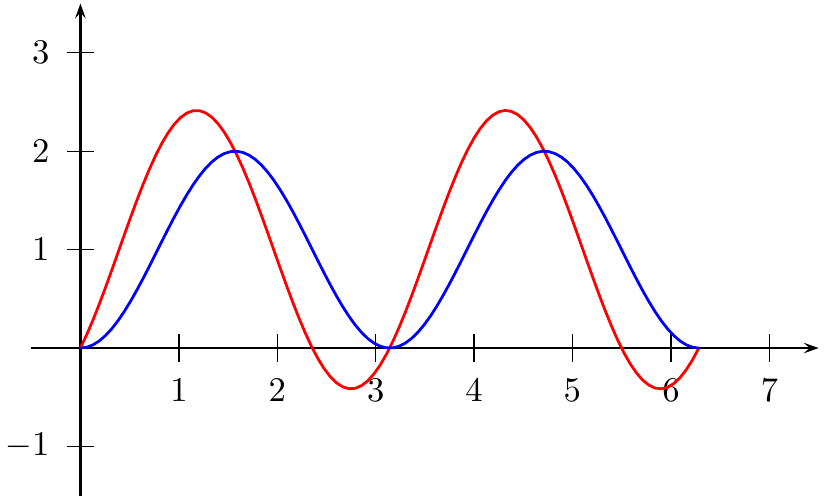}
\caption{\CaptionFigGraphn}\label{LONGLabelFigGraphn}
			\end{figure}

	Let $P=an\in AN$. 
	By proposition \ref{LONGPropUniquexxxxANANbarktheta}, we consider the unique  $x_0$, $x_1$, $x_2$ and $x_3$ in $\mathopen[ 0 , 2\pi [$ such that 
\begin{subequations}
	\begin{align}
		[P e^{x_0q_0}]&\in[AN]\\
		[P e^{x_1q_0}]&\in[ANk_{\theta}]\\
		[P e^{x_2q_0}]&\in[A\bar N]\\
		[P e^{x_3q_0}]&\in[A\bar Nk_{\theta}].
	\end{align}
\end{subequations}
They satisfy $x_0=0$, $x_1=\pi$ and $x_3=x_2+\pi$ modulo $2\pi$. Now, we divide $[AN]$ into two parts. The elements of $[AN]\cap [A\bar N]$ and $[AN]\cap[A\bar Nk_{\theta}]$ are said to be of \emph{type I}, while the other are said to be of \emph{type II}. We are going to prove that type I points are exactly blue points, while type II points are the red ones.

If $P$ is a point of type II, we know that the $x_i$ are four different numbers so that the norm function $n_P(x)$ vanishes \emph{at least} four times on the interval $\mathopen[ 0 , 2\pi [$, each of them corresponding to a point in the singularity. But our division of $[AN]$ into red and blue points shows that $n_P(x)$ can vanish \emph{at most} four times. We conclude that a point of type II is automatically red, and that the four roots of $n_P(x)$ correspond to the four values $x_i$ for which $P e^{x_iq_0}$ belongs to the singularity. The proposition is thus proved for points of type II.

Let now $[P]$ be of type I (say $P\in [AN]\cap [A\bar N]$) and let us show that $P$ is blue. We consider a sequence of points $P_k$ of type II which converges to $P$ (see remark \ref{LONGRemANANbarYapas}). We already argued that $P_k$ is red, so that $x_0(P_k)\neq x_2(P_k)$ and $x_1(P_k)\neq x_3(P_k)$, but
\begin{subequations}
	\begin{align}
		x_0(P_k)-x_2(P_k)\to 0\\
		x_1(P_k)-x_3(P_k)\to 0.
	\end{align}
\end{subequations}
The continuity of $n_Q(x)$ with respect to both $x\in\mathopen[ 0 , 2\pi [$ and $Q\in[AN]$ implies that $P$ has to be blue, and then $n_P(x)$ vanishes for exactly two values of $x$ which correspond to $P e^{xq_0}\in\sS$.

Let us now prove that everything is done. We begin by points of type I. Let $P$ be of type $I$ and say $P\in[AN]\cap[A\bar N]$. The curve $n_P(x)$ vanishes exactly two times in $\mathopen[ 0 , 2\pi [$. Now, if $P e^{x_1 q_0}\in[ANk_{\theta}]$, thus $x_1 = \pi$ and we also have $P e^{x_1q_0}\in[A\bar Nk_{\theta}]$, but $P$ does not belong to $[ANk_{\theta}]$, which proves that $n_P(x)$ vanishes \emph{at least} two times which correspond to the points $P e^{xq_0}$ that are in the singularity. Since the curve vanishes in fact exactly two times, we conclude that $n_P(x)$ vanishes if and only if $P e^{xq_0}$ belongs to the singularity.

If we consider a point $P$ of type II, we know that the values of $x_i$ are four different numbers, so that the curve $n_P(x)$ vanishes \emph{at least} four times, corresponding to the points $P e^{xq_0}$ in the singularity. Since the curve is in fact red, it vanishes \emph{exactly} four times in $\mathopen[ 0 , 2\pi [$ and we conclude that the curve $n_P(x)$ vanishes if and only if $P e^{xq_0}$ belongs to the singularity.

The conclusion follows from the fact that 
\begin{equation}
	AdS_l=\Big\{ [P e^{xq_0}] \tq \text{$P$ is of type I or II and $x\in\mathopen[ 0 , 2\pi [$} \Big\}.
\end{equation}

\end{proof}
Proof of theorem \ref{LONGThosSequivJzero} is now complete.

\subsection{Existence of the black hole}
\label{LONGSubSecExistenceTrouNoir}

We know that the geodesic trough $[g]$ in the direction $X$ is given by
\begin{equation}
	\pi\big( g e^{sX} \big)
\end{equation}
and that a geodesics is light-like when the \defe{direction}{Direction} $X$ is given by a nilpotent element in $\sQ$\cite{lcTNAdS}.
\newcommand{\CaptionFigGeodExistence}{We are looking at a geodesics issued from one point of the line $[\SO(2)]=\{ e^{xq_0}\}_{x\in\mathopen[ 0 , 2\pi [}$. Here, $E(w)=q_0+w_1q_1+w_2q_2+\sum_{k\geq 3}w_kq_k$ with $\sum_{k}w_k^2=1$.}

			\makeatletter
			\@ifundefined{lengthOffiguresbtzpy}{\newlength{\lengthOffiguresbtzpy}}
			\makeatother
			
\setlength{\lengthOffiguresbtzpy}{\totalheightof{$o=[\mtu]$}}
 \makeatletter 
				\@ifundefined{writeOffiguresbtzpy}			
				{\newwrite{\writeOffiguresbtzpy}
				\immediate\openout\writeOffiguresbtzpy=writeOffiguresbtzpy.pstricks.aux
				}
				\makeatother
\immediate\write\writeOffiguresbtzpy{totalheightofGeodExistencebaab:\the\lengthOffiguresbtzpy:}
\setlength{\lengthOffiguresbtzpy}{\widthof{$o=[\mtu]$}}
\immediate\write\writeOffiguresbtzpy{widthofGeodExistencebaac:\the\lengthOffiguresbtzpy:}
\setlength{\lengthOffiguresbtzpy}{\totalheightof{$[\SO(2)]$}}
\immediate\write\writeOffiguresbtzpy{totalheightofGeodExistencebaad:\the\lengthOffiguresbtzpy:}
\setlength{\lengthOffiguresbtzpy}{\widthof{$[\SO(2)]$}}
\immediate\write\writeOffiguresbtzpy{widthofGeodExistencebaae:\the\lengthOffiguresbtzpy:}
\setlength{\lengthOffiguresbtzpy}{\totalheightof{$[ e^{xq_0}]$}}
\immediate\write\writeOffiguresbtzpy{totalheightofGeodExistencebaaf:\the\lengthOffiguresbtzpy:}
\setlength{\lengthOffiguresbtzpy}{\widthof{$[ e^{xq_0}]$}}
\immediate\write\writeOffiguresbtzpy{widthofGeodExistencebaag:\the\lengthOffiguresbtzpy:}

\begin{figure}[ht]
\centering
\includegraphics{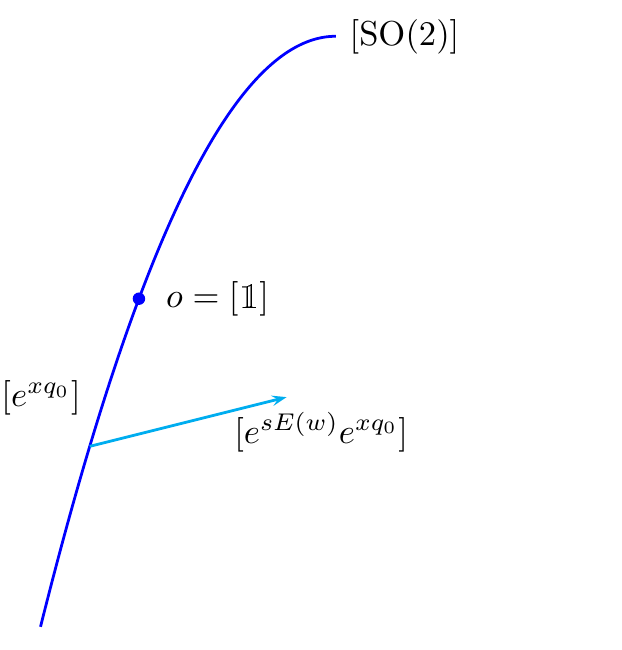}
\caption{\CaptionFigGeodExistence}\label{LONGLabelFigGeodExistence}
			\end{figure}

Let us study the geodesics issued from the point $[ e^{-xq_0}]$, see figure \ref{LONGLabelFigGeodExistence}. They are given by
\begin{equation}
	l^w_x(s)=\pi\big(    e^{-xq_0} e^{sE(w)} \big)
\end{equation}
where $E(w)=q_0+\sum_iw_iq_i$ with $\| w \|=1$. According to our previous work, the point $l^w_x(s)$ belongs to the singularity if and only if 
\begin{equation}		\label{LONGEqNormAFaireZeroOuPas}
	n_x^w(s)=\left\|   \pr_{\sQ} e^{-\ad(sE(w))} e^{\ad(xq_0)}J_1  \right\|^2=0.
\end{equation}
We already computed that $ e^{\ad(xq_0)}J_1=\cos(x)J_1+\sin(x)q_2$. By construction, $E(w)$ is nilpotent and $\ad(E)^3=0$ by proposition \ref{LONGEtOrdreDeux}. Using the fact that $[\sQ,\sH]\subset\sQ$ and $[\sQ,\sQ]\subset\sH$, we collect the terms in $\sQ$ in the development of the exponential. The $\sQ$ component of
\begin{equation}
	e^{-s\ad(E)}\big( \cos(x)J_1+\sin(x)q_2 \big)
\end{equation}
is 
\begin{equation}
	\ell=\frac{ s^2 }{ 2 }\sin(x)\ad(E)^2q_2-s\cos(x)\ad(E)J_1+\sin(x)q_2.
\end{equation}
The square norm of that expression is \emph{a priori} a polynomial of order $4$. Hopefully, the coefficient of $s^4$ contains
\begin{equation}
	B\big( \ad(E)^2q_2,\ad(E)^2q_2 \big),
\end{equation}
while the coefficient of $s^3$ is given by
\begin{equation}
	B\big( \ad(E)J_2,\ad(E)^2q_2 \big).
\end{equation}
Both of these two expressions are zero because the $\ad$-invariance of the Killing form makes appear $\ad(E)^3$. Equation \eqref{LONGEqNormAFaireZeroOuPas} is thus the second order polynomial given by
\begin{equation}		\label{LONGEqnwxBBB}
	\begin{aligned}[]
		n_x^w(s)	&=s^2\sin^2(x)B\big( \ad(E)^2q_2,q_2 \big)\\
				&\quad+s^2\cos^2(x)B\big( \ad(E)J_1,\ad(E) J_1\big)\\
				&\quad-2s\cos(x)\sin(x)B\big( \ad(E)J_1,q_2 \big)\\
				&\quad+\sin^2(x)B(q_2,q_2).
	\end{aligned}
\end{equation}
The problem now reduces to the evaluation of the three Killing products in this expression. %
Let us begin with $B\big( \ad(E)^2q_2,q_2 \big)$. For this one, we need to know the $q_2$-component of $\ad(E)^2q_2$. We have to review all the possibilities $\ad(q_i)\ad(q_j)q_2$ and determine which one(s) have a $q_2$-component.

In this optic, let us recall that $q_2$ is characterised by
\begin{equation}
	q_2\in\sP\cap\sQ\cap\tilde\sN_2.
\end{equation}

All the combinations $\ad(q_i)^2$ work. Since $\sG=\sK\oplus\sP$ is reductive and since $q_0$ is the only basis element of $\sQ$ to belong to $\sK$, none of the combinations $\ad(q_i)\ad(q_0)$ or $\ad(q_0)\ad(q_i)$ work. Using the relations
\begin{equation}
	\begin{aligned}[]
		q_1&\in\sP\cap\sQ\cap\sA	,	&&	[\tilde\sN_k,\tilde\sN_2]\subset\tilde\sN_k\\
		q_2&\in\sP\cap\sQ\cap\tilde\sN_2,	&&	[\tilde\sN_2,\tilde\sN_2]\subset\sA\\
		q_k&\in\sP\cap\sQ\cap\tilde\sN_k,	&&	[\tilde\sN_k,\tilde\sN_k]\subset\sA\oplus\tilde\sN_2,
	\end{aligned}
\end{equation}%
we check that the only working combinations are $\ad(q_i)^2$.

Thus, the only elements $\ad(q_i)\ad(q_j)q_2$ which have a $q_2$-component are $\ad(q_i)^2q_2$, while theorem \ref{LONGThoAdESqqq} says that this component is $q_2$ for $2\neq i\neq 0$ and $-q_2$ for $i=0$. Therefore, the $q_2$-component of $\ad(E)^2q_2$ is
\begin{equation}
	\ad(q_0)^2q_2+w_1^2\ad(q_1)^2q_2+\sum_{k\geq 3}w_k^2\ad(q_k)^2q_2=-w_2^2q_2
\end{equation}
where we used the fact that $\sum_i w_i^2=1$. Thus we have
\begin{equation}		\label{LONGEqBeDeuxqqwDeux}
	B\big( \ad(E)^2q_2,q_2 \big)=-w_2^2B(q_2,q_2).
\end{equation}

Let us now search for the $q_2$-component of $\ad(E)J_1$. We have $[q_1,J_1]\in[\sA,\sA]=0$, $[q_k,J_1]=0$ (equation \eqref{LONGEqJUnqkzero}), and $[q_2,J_1]=-q_0$, $[q_0,J_1]=-q_2$ (equation \eqref{LONGEqCalculBBBJUnUnNirme}). Then, we have
\begin{equation}
	\ad(E)J_1=w_2q_0+q_2.
\end{equation}
That implies
\begin{equation}
	B\big( \ad(E)J_1,q_2 \big)=B(q_2,q_2),
\end{equation}
and
\begin{equation}
	B\big( \ad(E)J_1,\ad(E)J_1 \big)=B(q_2,q_2)+w_2^2B(q_0,q_0).
\end{equation}
Equation \eqref{LONGEqnwxBBB} now reads
\begin{equation}
	\begin{aligned}[]
		\frac{ n_x^w(s) }{ B(q_2,q_2) }=\big( \cos^2(x)-w^2_2 \big)s^2-2\cos(x)\sin(x)s+\sin^2(x).
	\end{aligned}
\end{equation}
We have $n_x^w(s)=0$ when $s$ equals
\begin{equation}
	s_{\pm}=\frac{ \cos(x)\sin(x)\pm| w_2\sin(x) | }{ \cos^2(x)-w_2^2 }.
\end{equation}
If $w_2\sin(x)\geq 0$, we have\footnote{The solutions \eqref{LONGEqRacinesellTN} were already deduced in \cite{lcTNAdS} in a quite different way.}
\begin{equation}				\label{LONGEqRacinesellTN}
	\begin{aligned}[]
		s_+	&=\frac{ \sin(x) }{ \cos(x)-w_2 }&\text{and}&&
		s_-	&=\frac{ \sin(x) }{ \cos(x)+w_2 },
	\end{aligned}
\end{equation}
and if $w_2\sin(x)<0$, we have to exchange $s_+$ with $s_-$.

If we consider a point $ e^{xq_0}$ with $\sin(x)>0$ and $\cos(x)<0$, the directions $w$ with $| w_2 |<| \cos(x) |$ escape the singularity as the two roots \eqref{LONGEqRacinesellTN} are simultaneously negative. Such a point does not belong to the black hole. That proves that the black hole is not the whole space.

If we consider a point $ e^{xq_0}$ with $\sin(x)>0$ and $\cos(x)>0$, we see that for every $w_2$, we have $s_+>0$ or $s_->0$ (or both). That shows that for such a point, every direction intersect the singularity. Thus the black hole is actually larger than only the singularity itself.

The two points with $\sin(x)=0$ belong to the singularity. At the points $\cos(x)=0$, $\sin(x)=\pm1$, we have $s_+=-1/w_2$ and $s_-=1/w_2$. A direction $w$ escapes the singularity only if $w_2=0$ (which is a closed set in the set of $\| w \|=1$).

\newcommand{\CaptionFigCercleK}{Points in $\pi(K)$ are classified by their angle in $\SO(2)$. Red points are part of the singularity, points in the black zone belong to the black hole and points in the green zone are free. The upper and lower boundaries belong to the horizon.}

			\makeatletter
			\@ifundefined{lengthOffiguresbtzpy}{\newlength{\lengthOffiguresbtzpy}}
			\makeatother
			
\setlength{\lengthOffiguresbtzpy}{\totalheightof{$K_H$}}
 \makeatletter 
				\@ifundefined{writeOffiguresbtzpy}			
				{\newwrite{\writeOffiguresbtzpy}
				\immediate\openout\writeOffiguresbtzpy=writeOffiguresbtzpy.pstricks.aux
				}
				\makeatother
\immediate\write\writeOffiguresbtzpy{totalheightofCercleKbaab:\the\lengthOffiguresbtzpy:}
\setlength{\lengthOffiguresbtzpy}{\widthof{$K_H$}}
\immediate\write\writeOffiguresbtzpy{widthofCercleKbaac:\the\lengthOffiguresbtzpy:}

\begin{figure}[ht]
\centering
\includegraphics{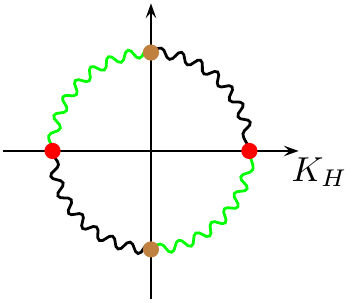}
\caption{\CaptionFigCercleK}\label{LONGLabelFigCercleK}
			\end{figure}

\section{Some more computations}
\label{LONGSecMoreComputations}

As we saw, the use of theorem \ref{LONGThosSequivJzero} leads us to study the function
\begin{equation}			\label{LONGEqprSqbignor}
	n_{[g]}^w(s)=\|   \pr_{\sQ} \Ad\big(  e^{-s E(w)} \big)X  \|^2=0
\end{equation}
where $E(w)=q_0+w_1q_1+\ldots+w_{l-1}q_{l-1}$ ($w\in S^{l-2}$) and (see equation \eqref{LONGEqAdanJUnabck})
\begin{equation}
	X=e^{\ad(xq_0)}\Ad(na)J_1=e^{\ad(xq_0)}\big(  J_1+aX_{++}+bX_{+-}+\sum_{k=3}^{l-1}c_kX_{+0}^k\big).
\end{equation}

From proposition \ref{LONGEtOrdreDeux} we have $\ad(E)^3=0$ and the exponential in equation \eqref{LONGEqprSqbignor} contains only three terms. Using the fact that $\pr_{\sQ}\ad(E)X=\ad(E)X_{\sH}$ and  we are then lead to study the norm of
\begin{equation}
	X_{\sQ}-s\ad(E)X_{\sH}+\frac{ s^2 }{2}\ad(E)^2X_{\sQ}.
\end{equation}
Notice that, since $\sQ$ is Killing-orthogonal to $\sH$, we have $B(X_{\sQ},Y_{\sQ})=B(X,Y_{\sQ})$. Thus we have
\begin{equation}
	n_{[g]}^w(s)=\|   \pr_{\sQ} \Ad\big(  e^{-s E(w)} \big)X  \|^2=a(E)s^2+b(E)s+c
\end{equation}
where
\begin{subequations}		\label{LONGSubEqsabcEBBB}
	\begin{align}
		a(E)&=-B\big( \ad(E)X,\sigma\ad(E)X \big)		\label{LONGEqCoefaEBX}\\
		b(E)&=-2B\big( X_{\sQ},\ad(E)X_{\sH} \big)\\
		c&=B(X_{\sQ},X_{\sQ}).
	\end{align}
\end{subequations}
It is convenient to decompose $X$ into $X=X_2+\sum_kc_kX_k$ with $X_2\in\tilde\sN_2$ and $X_k\in\tilde\sN_k$ as well as $E=E_2+\sum_kE_k$. Using the decompositions \eqref{LONGSubeqsDecompXqps}, the first part is
\begin{equation}
	\begin{aligned}[]
		X_2&=\Ad( e^{xq_0})\big( a(q_0-q_2)+b(p_1+s_1) \big)	\label{LONGEqXsiAdSTrois}\\
		X_k&=\Ad( e^{xq_0})(p_k+s_k)
	\end{aligned}
\end{equation}
where we have renamed $a$, $b$ and $c_k$. Using the exponentials \eqref{LONGSubEqsexpxzqpsdzuu} and \eqref{LONGEqExpoQzSurNkb} we have
\begin{subequations}
	\begin{align}
		X_2 &= aq_0-b\sin(x)q_1-a\cos(x)q_2 +a\sin(x)J_1+b s_1+b\cos(x)p_1	\label{LONGsubEqXtroisdonne}\\
		X_k &= -s_k-\cos(x)p_k+\sin(x)q_k					\label{LONGsubEqXkdonne}
	\end{align}
\end{subequations}

Let us first look at the case we encounter in $AdS_3$, i.e. with $X_k=0$ and $E_k=0$. In this case, a general direction is given by
\begin{equation}
	E_2(w)=E(\theta)=q_0+\cos(\theta)q_1+\sin(\theta)q_2.
\end{equation}
We achieve the computation of $ e^{\ad(E(\theta))}X_2$ using the known commutators:
\begin{subequations}		\label{LONGSubEqsEXtrois}
\begin{equation}
	\begin{aligned}[]
		\ad(E_2)X_{\sQ}&=J_1\big( a\sin(\theta)+a\cos(x) \big)\\
				&\quad+p_1\big( -a\cos(\theta)-b\sin(x) \big)\\
				&\quad+s_1\big( b\sin(x)\sin(\theta)-a\cos(x)\cos(\theta) \big).
	\end{aligned}
\end{equation}
and
\begin{equation}
	\begin{aligned}[]
		\ad(E_2)X_{\sH}&=q_0\big( a\sin(x)\sin(\theta)-b\cos(x)\cos(\theta) \big)\\
				&\quad-q_1 b\big( \sin(\theta)+\cos(x) \big)\\
				&\quad+q_2\big( a\sin(x)+b\cos(\theta) \big).
	\end{aligned}
\end{equation}
\end{subequations}
Putting together and writing $\big( \cos(\theta),\sin(\theta) \big)=(w_1,w_2)$,
\begin{equation}\label{LONGEqadEtroissurXtrois}
	\begin{aligned}[]
		\ad(E_2)X_2 &= J_1\big( a\cos(x)+aw_2 \big)\\
			&\quad p_1\big( -b\sin(x)-aw_1 \big)\\
			&\quad s_1\big( -a\cos(x)w_1+b\sin(x)w_2 \big)\\
			&\quad q_0\big( -b\cos(x)w_1+a\sin(x)w_2 \big)\\
			&\quad q_1\big( -bw_2-b\cos(x) \big)\\
			&\quad q_2\big( bw_1+a\sin(x) \big)
	\end{aligned}
\end{equation}
Using the expression \eqref{LONGEqadEtroissurXtrois} among with the norms and collecting the terms with respect to the dependence in $\theta$, we have
\begin{equation}			\label{LONGEqKillingEXQEXQ}
	\begin{aligned}[]
		\frac{ B\big( \ad(E)X_{\sQ},\ad(E)X_{\sQ} \big)}{B(q_0,q_0)}&=a^2\cos^2(x)-a^2\\
				&\quad +\sin(\theta)\big( -2a^2\cos(x) \big)\\
				&\quad +\cos(\theta)\big( -2ab\sin(x) \big)\\
				&\quad +\cos^2(\theta)\big( a^2\cos^2(x)-b^2\sin^2(x) \big)\\
				&\quad +\sin(\theta)\cos(\theta)\big( -2ab\sin(x)\cos(x) \big),
	\end{aligned}
\end{equation}
and
\begin{equation}
	\begin{aligned}[]
		\frac{ B\big( \ad(E)X_{\sH},\ad(E)X_{\sH} \big)}{B(q_0,q_0)}&= -b^2(1+\cos^2(x))\\
							&\quad +\sin(\theta)\big( -2b^2\cos(x) \big)\\
							&\quad +\cos(\theta)\big( -2ab\sin(x) \big)\\
							&\quad +\cos^2(\theta)\big( b^2\cos^2(x)-a^2\sin^2(x) \big)\\
							&\quad +\sin(\theta)\cos(\theta)\big( -2ab\sin(x)\cos(x) \big)
	\end{aligned}
\end{equation}
and finally, in the case of $AdS_3$ we have
\begin{equation}		\label{LONGEqaEdansAdsTrois}
	\begin{aligned}[]
		a(E)&=B\big( \ad(E)X_{\sH},\ad(E)X_{\sH} \big)-B\big( \ad(E)X_{\sQ},\ad(E)X_{\sQ} \big)\\
			&=(a^2-b^2)\big( \cos^2(x)+\sin(\theta)\cos(x)+\sin^2(\theta) \big).
	\end{aligned}
\end{equation}

For investigating the higher dimensional cases, we decompose $\ad(E)X$ into the four parts $\ad(E_2)X_2$, $\ad(E_2)X_k$, $\ad(X_k)X_2$ and $\ad(E_k)X_k$. 

The action of $E_2$ on $\tilde\sN_k$ is given by
\begin{subequations}
	\begin{align}
		\ad(E_2)r_k&=w_1q_k\\
		\ad(E_2)p_k&=-q_k\\
		\ad(E_2)q_k&=p_k+w_1r_k-w_2s_k\\
		\ad(E_2)s_k&=-w_2q_k.
	\end{align}
\end{subequations}
Thus 
\begin{equation}\label{LONGadEsurXki}
		\ad(E_2)X_k	=	q_k\big( w_2+\cos(x) \big)
				 + p_k\big( \sin(x) \big)
				 + r_k \big( w_1\sin(x) \big)
				 + s_k \big( -w_2\sin(x) \big)
\end{equation}
The same way we find
\begin{equation}
	\begin{aligned}[]
		\ad(E_k)X_2 = \sum_{k=3}^{l-1}\Big[ & p_k(-aw_k)
						+ r_k\big( b\sin(x)w_k\big)
						+ s_k\big( -a\cos(x)w_k \big)
		\Big].
	\end{aligned}
\end{equation}
Using the relations
\begin{equation}
	\ad(q_k)X_k=-q_2+\cos(x)q_0.
\end{equation}
and
\begin{subequations}
	\begin{align}
		\ad(E)q_0&=-w_1p_1+w_2J_1-\sum_{k=1}^lw_kp_k\\
		\ad(E)q_1&=p_1-w_2s_1-\sum_{k=3}^lw_kr_k\\
		\ad(E)q_2&=-J_1+w_1s_1+\sum_{k=3}^ls_k\\
		\ad(E)q_k&=p_k+w_1r_k-w_2s_k,
	\end{align}
\end{subequations}
we find
\begin{equation}
	\ad(E_k)X_{k'}=\begin{cases}
		c_{k'}w_k\sin(x)	&	\text{if $k\neq k'$}\\
		-c_kw_k\cos(q_0)-c_kw_kq_2	&	 \text{if $k=k'$}
	\end{cases}
\end{equation}
Putting all together, we find
\begin{equation}
	\begin{aligned}[]
		\ad(E)X	& = J_1\big( a\cos(x)+aw_2 \big)\\
			&\quad +  p_1\big( -b\sin(x)-aw_1 \big)\\
			&\quad+  s_1\big( -a\cos(x)w_1+bw_2\sin(x) \big)\\
			&\quad  +  q_0\big( -bw_1\cos(x)+aw_2\sin(x)+\sum_{k}c_kw_k\cos(x) \big)\\
			&\quad +  q_1\big( -bw_2-b\cos(x) \big)\\
			&\quad +  q_2\big( bw_1+a\sin(x)-\sum_kc_kw_k \big)\\
			&\quad +  \sum_kq_kc_k\big( w_2+\cos(x) \big)\\
			&\quad +  \sum_kp_k\big( c_k\sin(x)-aw_k \big)\\
			&\quad +  \sum_k r_k\big( c_kw_1+bw_k \big)\sin(x)\\
			&\quad +  \sum s_k\big( -c_kw_2\sin(x)-aw_k\cos(x) \big)\\
			&\quad + \sum_{k\geq 3}\sum_{k'>k} r_{kk'}\big( c_kw_{k'}-c_{k'}w_k \big)\sin(x)
	\end{aligned}
\end{equation}
It is quite easy but long to compute $a(E)$, $b(E)$ and $c$ from that expression. The results are

\begin{subequations}\label{LONGEqCoefsabcBE}
	\begin{align}
		\frac{ a(E) }{ B(q_0,q_0) }&=M\Big( w_2^2+\cos(x)w_2+\cos^2(x)\Big)	   \label{LONGEqCoeffaE}\\
		\frac{ b(E) }{ B(q_0,q_0) }&=-2M\sin(x)\big( w_2+\cos(x)\big))\\
		\frac{ c }{ B(q_0,q_0) }&=M\sin^2(x)
	\end{align}
\end{subequations}
where $M=\big( a^2-b^2-\sum_kc_k^2 \big)$. The important point to notice is that these expressions only depend on the $w_2$-component of the direction. Notice that $c=0$ if and only if the point $[g]$ belongs to the singularity because $s=0$ is a solution of \eqref{LONGEqprSqbignor} only in the case $[g]\in\hS$.

We can avoid the computation of a certain number of terms by exploiting the properties of the decomposition $\sG=\mZ_{\sK}(\sA)\oplus\tilde\sN_2\oplus\tilde\sN_k\oplus\sA$. The dependence in $w_1^2$ of $a(E)$ is given by the term
\begin{equation}
	B\big( \ad(q_1)X,\sigma\ad(q_1)X \big)=B\big( \ad(q_1)^2X,\sigma X \big).
\end{equation}
This is easily computed using the theorems \ref{LONGThoAdSqIouZero} and \ref{LONGThoBaisXXorthoigher}. The result is that the coefficient of $w_1^2$ in $a(E)/B(q_0,q_0)$ is
\begin{equation}
	-\sin^2(x)(a^2-b^2- C^2)
\end{equation}
where $C^2=\sum_{k\geq 3}c_k^2$.

The term which does not depend on $w$ is 
\begin{equation}
	B\big( \ad(q_0)X,\sigma\ad(q_0)X \big)=B\big( \ad(q_0)^2X,\sigma X \big).
\end{equation}
The result is that the independent term in $a(E)/B(q_0,q_0)$ is
\begin{equation}
	\big( 1-2\cos^2(x) \big)(a^2-b^2-C^2).
\end{equation}
In the same way, the coefficient of $w_2^2$ is $B\big( \ad(q_2)^2X,\sigma X \big)$ and we find
\begin{equation}
	-\big( \sin^2(x)+1 \big)(a^2-b^2-C^2).
\end{equation}

\begin{remark}	\label{LONGRemImapoabcE}
	Importance of the coefficients \eqref{LONGEqCoefsabcBE}. If $v\in\hF_l$, there is a direction $E_0$ in $AdS_l$ which escapes the singularity from $v$. Thus the polynomial $a(E_0)s^2+b(E_0)s+c$ has only non positive roots. From the expressions \eqref{LONGEqCoefsabcBE}, we see that the polynomial corresponding to $\iota(v)$ is the same, so that the direction $E_0$ escapes the singularity from $\iota(v)$ as well.

	This is the main ingredient of the next section.
\end{remark}

\section{Description of the horizon}
\label{LONGSecHorizonSansMatrices}

\subsection{Induction on the dimension}

The horizon in $AdS_3$ is already well understood \cite{Keio,BTZ_horizon}. We are not going to discuss it again. We will study how does the causal structure (black hole, free part, horizon) of $AdS_{l}$ includes itself in $AdS_{l+1}$ by the inclusion map
\begin{equation}
	\iota\colon AdS_l\to AdS_{l+1}.
\end{equation}

\begin{lemma}		\label{LONGLemMemeQueLemQuatre}
	Let $[g]\in\iota(AdS_3)\subset AdS_l$ be outside the singularity. We suppose that there is an open set $\mO$ in $S^1$ of directions escaping the singularity from $[g]$. Then there exists an open set $\mO'$ in $S^{l-2}$ of directions escaping the singularity.
\end{lemma}

\begin{proof}
	The proof is a consideration about the coefficients \eqref{LONGSubEqsabcEBBB}. The hypothesis means that the points
\begin{equation}
	\pi\left( g e^{sE(\theta)} \right)
\end{equation}
do not belong to $\hS$ for $s\geq 0$ when $E(\theta)=q_0+\cos(\theta)q_1+\sin(\theta)q_2$ and $\theta$ belongs to the given open set $\mO\subset\mathopen[ 0 , 2\pi \mathclose]$. If $a(E_0)\neq 0$ for some $E_0\in\mO$, the solutions are given by
\begin{equation}
	s_{\pm}=\frac{ -b\pm\sqrt{b^2-4ac} }{ 2a }.
\end{equation}
In such a direction, there are two values, both outside\footnote{When we say ``outside'' of $\eR^+$, we include the case of complex solutions.} of $\eR^+$, of $s$ such that $[g e^{sE_0}]\in \hS$. By continuity, we can find a neighborhood of $E_0$ in $S^{l-2}$ such that $[g e^{sE}]$ belongs to the singularity only for non positive numbers.

A problem arises when $a(E)=0$ for every direction $E$ in the open set $\mO$. In that case the equation \eqref{LONGEqprSqbignor} has only one solution which is negative by hypothesis. But it could appear that in every neighborhood of $E$, a second solution, positive, appears. If we write $X=\Ad(g^{-1})J_1$, what we have to prove is that the quantity
\begin{equation}
	a(E)=B\big( \ad(E)X,\sigma\ad(E)X \big)
\end{equation}
is not constant when $E$ runs over $\mO$, in particular, there exists a direction $\theta_0\in\mO$ such that $a(\theta_0)\neq 0$. We supposed that $[g]\in \iota(AdS_3)$, so that $X=\Ad(an)J_1$ is given by $X_3$ of equation \eqref{LONGsubEqXtroisdonne}.

The function $a(E)$ \eqref{LONGEqaEdansAdsTrois} is analytic with respect to $\theta$, thus if it vanishes on an open set $\mO$, it has to vanish everywhere. This can only be achieved with $a=\pm b$. Now, simple computation show that
\begin{equation}
	c=a^2-b^2\sin^2(x)-a^2\cos^2(x)=(a^2-b^2)\sin^2(x)
\end{equation}
which vanishes when $a=\pm b$, so that $a(E)$ can only be constant with respect to $E$ on the singularity. Thus we conclude that $a(E)$ is not constant with respect to $E\in S^1$ outside the singularity.

This concludes the proof of lemma \ref{LONGLemMemeQueLemQuatre}.
\end{proof}

\begin{lemma}
	The direction $E_0$ in $AdS_l$ escapes the singularity from $v\in AdS_l$ if and only if it escapes the singularity from $\iota(v)$ in $AdS_{l+1}$.
\end{lemma}

\begin{proof}
	The fact for $v$ to escape the singularity in the direction $E_0$ means that the equation
	\begin{equation}
		a_{v}(E_0)s^2+b_{v}(E_0)+c_{v}=0
	\end{equation}
	where the coefficients are given by \eqref{LONGEqCoefsabcBE} has no positive solutions with respect to $s$. Since these coefficients are the same for $v$ and $\iota(v)$, the equation for $\iota(v)$ is in fact the same and has the same solutions.
\end{proof}

As a warm up, let us prove the following, which is a particular case of lemma \ref{LONGLemDueiINtlIntlpu}.
\begin{lemma}
	Let $[g]\in AdS_l$ be such that there exists an open set $\mO\in S^1$ of directions that escape the singularity. Then there is an open set in $S^{l-1}$ that escapes the singularity from $i[g]\in AdS_{l+1}$.
\end{lemma}

\begin{proof}

	With the notations of section \ref{LONGSecMoreComputations}, we only have to compute $a(E)$ when $E=E_2$ and $X$ is general. So we pick the expression \eqref{LONGEqCoeffaE} and we put $w_2=\cos(\theta)$ while $w_k=0$ for every $k\geq 3$. What we have is
	\begin{equation}
		a(E)=M\Big( \cos^2(\theta)+\cos(x)\cos(\theta)+\cos^2(x) \Big)
	\end{equation}
	If $a(E)=0$ for every $\theta\in\mO$, then $M=0$ which is impossible since we suppose that the starting point does not belong to the singularity.
\end{proof}
\begin{lemma}		\label{LONGLemDueiINtlIntlpu}
	We have
	\begin{equation}
		\iota\big( \Int(\hF_l) \big)\subset\Int(\hF_{l+1})
	\end{equation}
	or, equivalently, 
	\begin{equation}
		\Adh(BH_{l+1})\cap\iota(AdS_l)\subset\iota\big( \Adh(BH_l) \big).
	\end{equation}
\end{lemma}

\begin{proof}
	Let $v'\in\Int(\hF_l)$ and $\mO$, an open set of directions in $AdS_l$ that escape the singularity. The coefficient $a_l(E)$ is not constant on $\mO$ because the coefficient $M=a^2-b^2-C^2$ is only zero on the singularity (see equation \eqref{LONGEqCoeffaE}). Thus we can choose $E_0\in\mO$ such that $a_l(E_0)\neq 0$. We consider $a_{l+1}(E_0)$, the coefficient of $s^2$ for the point $\iota(v')$ in the direction $E_0$. From the expression \eqref{LONGEqCoeffaE} we know that $a_{l+1}(E_0)=a_l(E_0)$. The coefficients $b(E_0)$ and $c$ are also the same for $v'$ and $\iota(v')$. 
	
	Since $a(E_0)\neq 0$ and $v'\in\Int(\hF_l)$, we have two solutions to the equation $a(E_0)s^2+b(E_0)s+c=0$ and both of these are outside $\eR^+_0$. This conclusion is valid for $v'\in AdS_l$ as well as for $\iota(v')\in AdS_{l+1}$. Then there is a neighborhood of $\iota(v')$ on which the two solutions keep outside $\eR^+_0$. That proves that $\iota(v')\in\Int(\hF_{l+1})$.

	For the second line, suppose that $v\in\iota(AdS_l)$ does not belong to $\iota\big( \Adh(BH_l) \big)$, thus $v\in\iota\big( \Int(\hF_l) \big)\subset\Int(\hF_{l+1})$. In that case $v$ does not belong to $\Adh(BH_{l+1})$.

\end{proof}

\begin{proposition}		\label{LONGProphFdanshF}
	We have
	\begin{equation}
		\hF_{l+1}\cap\iota(AdS_l)\subset\iota(\hF_l)
	\end{equation}
\end{proposition}

\begin{proof}
	If $v=\iota(v')\in\hF_{l+1}$, there is a direction $E_0$ in $AdS_{l+1}$ which escape the singularity from $v$. That direction is given by a vector $(w_1,\cdots v_l)\in S^{l}$. Since the coefficients $a(E)$, $b(E)$ and $c$ do only depend on $w_2$, a direction $(w'_1,\cdots,w'_{l-1},0)$ with $w'_2=w_2$ escapes the singularity from $v'$. This proves that $v'\in\hF_l$.
\end{proof}

\begin{lemma}		\label{LONGLemHiH}
	We have
	\begin{equation}
		\hH_{l+1}\cap\iota(AdS_l)\subset\iota(\hH_{l}).
	\end{equation}
\end{lemma}

\begin{proof}
	First,
	\begin{equation}
		v\in\hH_{l+1}\cap\iota(AdS_l)\subset\hF_{l+1}\cap\iota(AdS_l)\subset\iota(\hF_l)
	\end{equation}
	from proposition \ref{LONGProphFdanshF}. Now, let's take $v'\in\hF_l$ such that $v=\iota(v')$. We have to prove that $v'\in\hH_l$. Let us suppose that $v'\in\Int(\hF_l)$, so $v\in\Int(\hF_{l+1})$ because of lemma \ref{LONGLemDueiINtlIntlpu}. This is in contradiction with the fact that $v\in\hH_{l+1}$.
\end{proof}

\begin{corollary}		\label{LONGCorDeuxTrucsBHhH}
	We have
	\begin{multicols}{2}
		\begin{enumerate}

		\item
			$\iota(\hS_l)\subset\hS_{l+1}$,
		\item
			$\iota(\hF_l)\subset\hF_{l+1}$,
		\item
			$\iota(BH_l)\subset BH_{l+1}$,
		\item
			$\iota(\hH_l)\subset \hH_{l+1}$.

		\end{enumerate}
	\end{multicols}
\end{corollary}

\begin{proof}
	We have $\Ad\big( \iota(g^{-1}) \big)J_1=\Ad(g^{-1})J_1$, so that the condition of theorem \eqref{LONGThosSequivJzero} is invariant under $\iota$. Thus one immediately has $\iota(\hS_l)\subset\hS_{l+1}$ and $\iota(\hF_l)\subset\hF_{l+1}$.

	An element $v$ which does not belong to $BH_{l+1}$ belongs to $\hF_{l+1}$, but if $v$ belongs to $\iota(AdS_l)\cap\hF_{l+1}$, it belongs to $\iota(\hF_l)$ by proposition \ref{LONGProphFdanshF} and then does not belong to $\iota(BH_l)$. Thus $\iota(BH_{l})\subset BH_{l+1}$.

	Now if $v'\in\hH_l$, let us consider $\mO$, a neighborhood of $v=\iota(v')$ in $AdS_{l+1}$. The set $\iota^{-1}\big( \mO\cap\iota(AdS_l) \big)$ contains a neighborhood $\mO'$ of $v'$ in $AdS_{l}$. Since $v'\in\hH_l$, there is $\bar v\in\mO'$ such that $\bar v\in BH_l$. Thus $\iota(\bar v)\in\mO$ belongs to $BH_{l+1}$ by the first item.

\end{proof}

Let $X\in\sG_{l+1}$ such that $[X,J_1]=0$, and let $R$ be the group generated by $X$. The following results are intended to show that such a group can be used in order to transport the causal structure from $AdS_l$ to $AdS_{l+1}$

The key ingredient will be the fact that, since $X$ commutes with $J_1$, we have
\begin{equation}		\label{LONGEqAdOkSurJun}
	\Ad\big( (g e^{sE})^{-1} \big)J_1=\Ad\big( ( e^{\alpha X}g e^{sE})^{-1} \big)J_1.
\end{equation}

\begin{lemma}		\label{LONGLemRSsubsetS}
	A group $R$ as described above preserves the causal structure in the sense that
	\begin{multicols}{2}
	\begin{enumerate}

		\item
			$R\cdot\hS\subset\hS$
		\item
			$R\cdot BH\subset BH$
		\item
			$R\cdot \hF\subset\hF$
		\item
			$R\cdot \hH\subset\hH$.

	\end{enumerate}
	\end{multicols}
\end{lemma}

\begin{proof}
	From equation \eqref{LONGEqAdOkSurJun}, we deduce that a direction $E_0$ will escape the singularity from the point $[g]$ is and only if it escapes the singularity from the points $r[g]$ for every $r\in R$. The first three points follow.

	Now let $v\in\hH$ and $r\in R$ and let us prove that $r\cdot v\in BH$. By the third point, $r\cdot v\in\hF$. Let now $\mO$ be a neighborhood of $r\cdot v$. The set $r^{-1}\cdot \mO$ is a neighborhood of $v$ and we can consider $\bar v\in BH\cap r^{-1}\cdot\mO$. By the second point, $r\cdot \bar v$ is a point of the black hole in $\mO$.
\end{proof}

\begin{remark}		\label{LONGRemdqnqRSlsubsetSlpu}
	Combining corollary \ref{LONGCorDeuxTrucsBHhH} and lemma \ref{LONGLemRSsubsetS} we have
	\begin{multicols}{2}
	\begin{enumerate}
		\item
			$ R\cdot\iota(\hS_l)\subset\hS_{l+1}$
		\item
			$ R\cdot\iota(\hF_l)\subset\hF_{l+1}$
		\item
			$ R\cdot\iota(BH_{l})\subset BH_{l+1}$
		\item
			$ R\cdot\iota(\hH_{l})\subset\hH_{l+1}$.
	\end{enumerate}
	\end{multicols}
\end{remark}

\begin{theorem}		\label{LONGThoCausalPasseParR}
	If moreover the one parameter group $R$ has the property to generate $AdS_{l+1}$ (in the sense that $R\cdot\iota(AdS_l)=AdS_{l+1}$), then we have
	\begin{multicols}{2}
	\begin{enumerate}
		\item
			$ R\cdot\iota(\hS_l)=\hS_{l+1}$
		\item
			$ R\cdot\iota(\hF_l)=\hF_{l+1}$
		\item\label{LONGItemStrucalpb}
			$ R\cdot \iota(BH_{l})=BH_{l+1}$
		\item
			$ R\cdot\iota(\hH_{l})=\hH_{l+1}$.
	\end{enumerate}
	\end{multicols}
\end{theorem}

\begin{proof}
	The inclusions in the direct sense are already done in the remark \ref{LONGRemdqnqRSlsubsetSlpu}.

	Let $r= e^{\alpha X}$ be an element of $R$. Since, by assumption, we have $[X,J_1]=0$, the action of $r$ leaves invariant the condition of theorem \eqref{LONGThosSequivJzero}:
	\begin{equation}		\label{LONGEqalpharjnagitpas}
		\Ad\big( (g e^{sE})^{-1} \big)J_1=\Ad\big( ( e^{\alpha X}g e^{sE})^{-1} \big)J_1.
	\end{equation}
	\begin{enumerate}
		\item
			Let $[g]\in\hS_{l+1}$, there exists a $r\in R$ such that $r[g]\in \iota(AdS_l)$. There exists an element $g'\in G_l$ such that $rg=\iota(g')$. Now $[g']\in\hS_l$ because
			\begin{equation}
				\Ad(g'^{-1})J_1=\Ad(\iota(g'^{-1}))J_1=\Ad\big(  (rg)^{-1} \big)J_1=\Ad(g^{-1})J_1,
			\end{equation}
			but by assumption the norm of the projection on $\sQ$ of the right hand side is zero.

		\item
			If $v$ is free in $AdS_{l+1}$, there is a direction $E_0$ escaping the singularity from $v$ and an element $r\in R$ such that $v'=r\cdot v\in\iota(AdS_l)$. The point $v'$ is also free in $AdS_{l+1}$ as the direction $E_0$ works for $r\cdot v$ as well as for $v$. Thus by proposition \ref{LONGProphFdanshF} we have
			\begin{equation}
				v'\in\hF_{l+1}\cap\iota(AdS_l)\subset\iota(\hF_l)
			\end{equation}
			and $v\in R\cdot\iota(\hF_l)$.
		\item
			If $v\in BH_{l+1}$, the point $r\cdot v\in\iota(AdS_l)$ also belongs to $BH_{l+1}$. If $r\cdot v=\iota(v')$, then $v'\in BH_l$ from if $v'\in\hF_l$, then $\iota(v')\in\hF_{l+1}$.
		\item
			If $v\in\hH_{l+1}$, there exists a $r\in R$ such that $v'=r\cdot v\in\iota(AdS_l)$ and moreover, $v'$ belongs to the horizon in $AdS_{l+1}$ since the horizon is invariant under $R$. Thus $v'$ belongs to $\iota(AdS_l)\cap\hH_{l+1}\subset \iota(\hH_l)$ by lemma \ref{LONGLemHiH}.  Now, $v\in R\cdot\iota(\hH_l)$.

	\end{enumerate}
\end{proof}

\subsection{Examples of surjective groups}

Theorem \ref{LONGThoCausalPasseParR} describes the causal structure in $AdS_l$ by induction on the dimension provided that one knows a group $R$ such that $AdS_{l+1}=R\cdot \iota(AdS_l)$. Can one provide examples of such groups? The following proposition provides a one.

\begin{proposition}		\label{LONGPropSurjectif}
	If $R$ is the one parameter subgroup of $\SO(2,l)$ generated by $r_{l,l+1}$, then we have
	\begin{equation}
		R\cdot \iota(AdS_l)= AdS_{l+1}.
	\end{equation}
\end{proposition}

\begin{proof}
	If one realises $AdS_n$ as the set of vectors of length $1$ in $\eR^{2,n-1}$, $AdS_l$ is included in $AdS_{l+1}$ as the set of vectors with vanishing last component and the element $r_{l,l+1}$ is the rotation in the plane of the two last coordinates. In that case, we have to solve
	\begin{equation}
		e^{\alpha r_{l,l+1}}\begin{pmatrix}
			u'	\\ 
			t'	\\ 
			x'_1	\\ 
			\vdots	\\ 
			x'_{l-2}	\\ 
			0	
		\end{pmatrix}=
		\begin{pmatrix}
			u	\\ 
			t	\\ 
			x_1	\\ 
			\vdots	\\ 
			x_{l-2}	\\ 
			x_{l-1}	
		\end{pmatrix}
	\end{equation}
	with respect to $\alpha$, $u'$, $t'$ and $x'_i$. There are of course exactly two solution if $x_{l-2}^2+x_{l-1}^2\neq 0$.
\end{proof}

In fact, many others are available, as the one showed at the end of \cite{BTZ_horizon}. In fact, since, in the embedding of $AdS$ in $\eR^{2,n}$, the singularity is given by $t^2-y^2=0$, almost every group which leaves invariant the combination $t^2-y^2$ can be used to propagate the causal structure. One can found lot of them for example by looking at the matrices given in \cite{These}.

\subsection{Backward induction}

Using proposition \ref{LONGProphFdanshF}, lemma \ref{LONGLemHiH}, corollary \ref{LONGCorDeuxTrucsBHhH} and the fact that the norm of $J_1^*$ is the same in $AdS_l$ as in $\iota(AdS_l)\subset AdS_{l+1}$, we have
\begin{enumerate}
	\item
		$\iota(\hF_l)=\hF_{l+1}\cap\iota(AdS_l)$,
	\item
		$\iota(\hH_l)=\hH_{l+1}\cap\iota(AdS_l)$,
	\item 
		$\iota(\hS_l)=\hS_{l+1}\cap\iota(AdS_l)$.
\end{enumerate}
These equalities hold for $l\geq 3$. For $l=2$ we can take the latter as a definition and set
\begin{equation}
	\hS_2=\{ v\in AdS_2\tq \iota(v)\in\hS_3 \}.
\end{equation}
The Iwasawa decomposition of $\SO(2,1)$ is given by $\sA=\langle J_2\rangle$, $\sN=\langle X_+\rangle$, $\sK=\langle q_0\rangle$ where $X_+=p_1-q_0$. Notice that we \emph{do not} have $\iota(X_+)=X_{++}$. Instead we have $X_+=\frac{ 1 }{2}(X_{++}+X_{-+})$. Thus $\hS_2$ is not given by the closed orbits of $AN$ in $AdS_2$.

The light-like directions are given by the two vectors $E=q_0\pm q_1$. In order to determine if the point $[e^{-\alpha J_2}e^{-aX_+} e^{-xq_0}]$ belongs to the black hole, we follow the same way as in section \ref{LONGSubSecExistenceTrouNoir}: we compute the norm
\begin{equation}
	\big\| \pr_{\sQ}  e^{-s\ad(E)}e^{x\ad(q_0)} e^{a\ad(X_+)} e^{\alpha\ad(J_2)}J_1 \big\|
\end{equation}
and we see under which conditions it vanishes. 
\begin{equation}
	e^{a\ad(X_+)}J_1=J_1+a(q_2+s_1),
\end{equation}
and then
\begin{equation}
	\begin{aligned}[]
		X&= e^{x\ad(q_0)}\big( J_1+a(q_2+s_1) \big)\\
		&=J_1\big( \cos(x)-a\sin(x) \big)+q_2\big( \sin(x)+a\cos(x) \big)+as_1.
	\end{aligned}
\end{equation}
With $E=q_0+q_1$, we have
\begin{equation}
	\begin{aligned}[]
		e^{s\ad(E)}J_1&=\big( -\frac{ 1 }{2}s^2+1 \big)J_1+sq_2+\frac{ 1 }{2}s^2s_1\\
		e^{s\ad(E)}q_2&=-sJ_1+q_2+ss_1\\
		e^{s\ad(E)}s_1&=-\frac{ 1 }{2}s^2J_1+sq_2+(\frac{ 1 }{2}s^2+1)s_1.
	\end{aligned}
\end{equation}
Thus
\begin{equation}
	\pr_{\sQ} e^{s\ad(E)}X=\Big( \big( \cos(x)-a\sin(x)+a \big)s+(\sin(x)+a\cos(x)) \Big)q_2.
\end{equation}
Its norm vanishes for the value of $s$ given by
\begin{equation}
	s^+=\frac{ a\cos(x)+\sin(x) }{ \big( \sin(x)-1 \big)a-\cos(x) }.
\end{equation}
The same computation with $E=q_0-q_1$ provides the value
\begin{equation}
	s^-=\frac{ a\cos(x)+\sin(x) }{ \big( \sin(x)+1 \big)a-\cos(x) }.
\end{equation}

For small enough $a$, the sign of $s^+$ and $s^-$ are both given the sign of $-\tan(x)$ that can be either positive or negative. Thus there is an open set of points in $AdS_2$ which intersect the singularity in every direction and an open set of points which escape the singularity.

As a side note, the singularity $\hS_2$ described here is not given by the closed orbits of $AN$ or $A\bar N$. Indeed, we show that
\begin{equation}
	\| \pr_{\sQ}\Ad( e^{a\iota(X_+)})J_1 \|^2=\| \pr_{\sQ}(J_1+aq_2+as_1) \|=-4a^2\neq 0.
\end{equation}
Thus the points of $N$ are not part of the singularity.

\pagestyle{headings}

\bibliographystyle{unsrt}
\def\polhk#1{\setbox0=\hbox{#1}{\ooalign{\hidewidth
  \lower1.5ex\hbox{`}\hidewidth\crcr\unhbox0}}}

\end{document}